\documentclass[nohyperref]{article}

\usepackage{microtype}
\usepackage{graphicx}
\usepackage{subfigure}
\usepackage{booktabs} 

\usepackage{hyperref}

\usepackage[accepted]{icml2022}

\usepackage{amsmath}
\usepackage{amssymb}
\usepackage{mathtools}
\usepackage{amsthm}
\usepackage{titletoc}

\usepackage[capitalize,noabbrev]{cleveref}
\usepackage{amsmath}
\usepackage[utf8]{inputenc} 
\usepackage[T1]{fontenc}    
\usepackage{hyperref}       
\usepackage{url}            
\usepackage{booktabs}       
\usepackage{amsfonts}       
\usepackage{nicefrac}       
\usepackage{microtype}      
\usepackage{mathrsfs}
\usepackage{flushend}
\usepackage{wrapfig}
\usepackage{lipsum}
\usepackage{etoolbox}
\usepackage{graphicx}
\usepackage{subfigure}
\usepackage{bbm}
\usepackage{amssymb}
\usepackage{amsthm}
\usepackage{color}
\usepackage{xcolor}
\usepackage{mdframed}
\usepackage{enumitem}
\setlist[itemize]{noitemsep, topsep=0pt, label=$\blacktriangleright$, leftmargin=*}
\usepackage[suppress]{color-edits}
\addauthor{mh}{red}
\addauthor{st}{blue}

\def\bblambda{\overbar{\blambda}'}

\def\b1{\mathbf 1}

\def\bb{\mathbf b}
\def\bd{\mathbf d}

\def\bx{\mathbf x}

\def\bg{\mathbf g}

\def\bs{\mathbf s}

\def\bw{\mathbf w}

\def\by{\mathbf y}

\def\bv{\mathbf v}

\def\bz{\mathbf z}

\def\bB{\mathbf B}
\def\bC{\mathbf C}

\def\bS{\mathbf S}

\def\bW{\mathbf W}

\def\blambda{\boldsymbol \lambda}
\def\wblambda{\widetilde{\boldsymbol \lambda}}
\def\bmu{\boldsymbol \mu}

\def\cf{\emph{\textmd{F}}}

\def\bmu{\boldsymbol \mu}

\newcommand{\T}{\scriptscriptstyle T}

\def\leref#1{Lemma~\ref{#1}}
\def\prref#1{Proposition~\ref{#1}}
\def\thref#1{Theorem~\ref{#1}}

\def\figref#1{Figure~\ref{#1}}
\def\figtab#1{Table~\ref{#1}}

\def\asref#1{Assumption~\ref{#1}}

\def\secref#1{Section~\ref{#1}}
\def\bydef{\triangleq}

\def\remark{\addtocounter{remark}{1}\def\@currentlabel{\theremark}%
\emph{Remark~\theremark}. } \makeatother
\newcounter{remark}
\newtheorem{proposition}{Proposition}
\newtheorem{lemma}{Lemma}
\newtheorem{theorem}{Theorem}
\newtheorem{corollary}{Corollary}
\newtheorem{definition}{Definition}
\newtheorem{assumption}{Assumption}

\newcommand{\overbar}[1]{%
  \mkern 1.5mu\overline{\mkern-1.5mu#1\mkern-1.5mu}\mkern 1.5mu}

\usepackage{amsopn}

\usepackage{pifont}

\theoremstyle{plain}
\theoremstyle{definition}
\theoremstyle{remark}

\usepackage[textsize=tiny]{todonotes}

\icmltitlerunning{A Single-Loop Gradient Descent and Perturbed Ascent Algorithm for Nonconvex Functional Constrained Optimization}

\begin{document}

\twocolumn[
\icmltitle{A Single-Loop Gradient Descent and Perturbed Ascent Algorithm for Nonconvex Functional Constrained Optimization}




\begin{icmlauthorlist}
\icmlauthor{Songtao Lu}{ibm}
\end{icmlauthorlist}

\icmlaffiliation{ibm}{IBM Research AI, IBM Thomas J. Watson Research Center, Yorktown Heights, New York 10598, USA}

\icmlcorrespondingauthor{Songtao Lu}{songtao@ibm.com}

\icmlkeywords{Machine Learning, ICML}

\vskip 0.3in
]



\printAffiliationsAndNotice{}  

\begin{abstract}
Nonconvex constrained optimization problems can be used to model a number of machine learning problems, such as multi-class Neyman-Pearson classification and constrained Markov decision processes. However, such kinds of problems are challenging because both the objective and constraints are possibly nonconvex, so it is difficult to balance the reduction of the loss value and reduction of constraint violation. Although there are a few methods that solve this class of problems, all of them are double-loop or triple-loop algorithms, and they require oracles to solve some subproblems up to certain accuracy by tuning multiple hyperparameters at each iteration. In this paper, we propose a novel gradient descent and perturbed ascent (GDPA) algorithm to solve a class of smooth nonconvex inequality constrained problems. The GDPA is a primal-dual algorithm, which only exploits the first-order information of both the objective and constraint functions to update the primal and dual variables in an alternating way.  The key feature of the proposed algorithm is that it is a single-loop algorithm, where only two step-sizes need to be tuned. We show that under a mild regularity condition GDPA is able to find Karush-Kuhn-Tucker (KKT) points of nonconvex functional constrained problems with convergence rate guarantees. To the best of our knowledge, it is the first single-loop algorithm that can solve the general nonconvex smooth problems with nonconvex inequality constraints. Numerical results also showcase the superiority of GDPA  compared with the best-known algorithms (in terms of both stationarity measure and feasibility of the obtained solutions).
\end{abstract}

\section{Introduction}

In this work, we consider the following class of nonconvex optimization problems under smooth nonconvex constraints
\begin{equation}\label{eq.org}
\mathbb{P}1\quad \min_{\bx\in\mathcal{X}\subseteq\mathbb{R}^d} \quad f(\bx)\quad\textrm{subject to}  \quad g(\bx)\le 0
\end{equation}
where functions $f(\bx):\mathbb{R}^d\to\mathbb{R}$ and $g(\bx):\mathbb{R}^d\to\mathbb{R}^m$ are smooth (possibly) nonconvex, $\mathcal{X}$ denotes the feasible set, and $m$ is the total number of constraints. This class of constrained optimization problems has been found very useful in formulating practical learning tasks, as the requirements of enhancing the interpretability of neural nets or safety and fairness guarantees raise. When the machine learning models are applied in different domains, functional constraint $g(\bx)$ can be specialized to particular forms. For example, if the first-order logic is considered in modeling the reasoning behaviors among the inputs, logical constraints will be incorporated in the optimization process \cite{hlmrf17,fischer19a}. Also, in a safe reinforcement learning problem, safety-aware constraints, e.g., cumulative long-term rewards \cite{yuya19,kaiqing20} or expected probability of failures \cite{sesh20},  will be included in the policy improvement step. Besides, in the design of deep neural nets (DNN) architectures, the energy consumption budget will be formulated as constraints in the DNN compression models \cite{yang2019ecc}. Different from the projection-friendly constraints, these constraints are (possibly) functional ones.

\subsection{Motivating Examples}

To be more specific, we give the following problems that can be formulated by \eqref{eq.org} as the motivating examples of this work.

\noindent{\bf Multi-class Neyman-Pearson Classification (mNPC)}. mNPC is a classic multi-class pattern recognition problem. The previous works \cite{weston1998multi,crammer2002learnability} propose to formulate this problem as a constrained optimization problem and adopt the support vector machine method to find $m$ classifiers. To be specific, it considers that there are $m$ classes of data, where each of them contains a data set, denoted by $\xi_i, i\in[m]$. The goal of the problem is to learn $m$ classifiers, denoted by $\{\bx_{(i)}\}$, so that the loss of a prioritized class is minimized while the rest ones are below a certain threshold denoted by $r_i$, i.e.,
\begin{subequations}
\begin{align}
\min_{\{\bx_{(i)}\}} \quad&f_{\xi_1}(\{\bx_{(i)}\}),
\\
\textrm{subject to}\quad &f_{\xi_i}(\{\bx_{(i)}\})\le r_i, i=2,\ldots,m,
\end{align}
\end{subequations}
where $f_{\xi_i}$ denotes the (possible nonconvex) non-increasing loss function of each class, and here class 1 is set as the prioritized one.

\noindent{\bf Constrained Markov decision processes (CMDP).} A CMDP is described by a tuple $(\mathcal{S},\mathcal{A},P,R,G_i,\gamma)$, where $\mathcal{S}$ is the state space, $\mathcal{A}$ is the action space of an agent, $P$ denotes the state transition probability of the MDP, $R$ is the reward, $G_i,i\in[m]$ are budget-related rewards,  and $\gamma$ denotes the discounted factor \cite{sutton2018reinforcement}. A policy $\pi: \mathcal{S}\to\Delta_{\mathcal{A}}$ of the agent is a function, mapping from the state space to action space, which determines a probability simplex $\Delta_{\mathcal{A}}$. The goal of the policy improvement problem is to learn an optimal policy so that the cumulative rewards can be maximized at the agent under some expected utility constraints, which can be formulated as the following constrained problem with respect to $\pi$ \cite{kaiqing20}:
\begin{equation}\label{eq.cmdpo}
\max_{\pi\in\Delta_{\mathcal{A}}}\; \langle Q^{\pi},\pi \rangle, \quad \textrm{subject to}\quad \langle Q^{\pi}_i,\pi\rangle \ge b_i, \forall i
\end{equation}
where $Q^{\pi}\bydef(1-\gamma)\mathbb{E}[\sum^{\infty}_{t=0}\gamma^tR(S_t,A_t)|S_0=s,A_0=a]$ denotes the action value functions associated with policy $\pi$, and similarly $Q^{\pi}_i\bydef(1-\gamma)\mathbb{E}[\sum^{\infty}_{t=0}\gamma^tG_i(S_t,A_t)|S_0=s,A_0=a],\forall i$ are the action value functions related to the constraints, $b_i,\forall i$ stand for the thresholds of each budget. In practice, policy $\pi$ is parametrzied by a neural network. Again, the above CMDP problem is a nonconvex optimization problem with nonconvex functional constraints \eqref{eq.org}.

\noindent{\bf Deep neural networks training (DNN) under Energy Budget}. One efficient technique to reduce the complexity of DNNs is model compression. Consider a layer-wise weights sparsification problem \cite{yang2019ecc}. Let $\bW$ denote the stacks of weight tensors of all the layers, i.e., $\bW\bydef\{\bw_{(u)}\}$, where $u$ is the index of layers, and $\bS\bydef\{\bs_{(u)}\}$ denote the stacks of the non-sparse weights of all the layers. Then, the energy-constrained DNN training problem can be written as
\begin{subequations}
\begin{align}
\min_{\bW,\bS}\quad &\ell(\bW),
\\
\textrm{subject to}\quad &\phi(\bw_{(u)})\le s_{(u)}, \psi(\bS)\le E_{\textrm{budget}},\forall u
\end{align}
\end{subequations}
where $\ell(\bW)$ denotes the training loss, $s_{(u)}$ corresponds to the sparsity level of layer $u$, $\phi(\bw_{(u)})$ calculates the layer-wise sparsity, $\psi(\bS)$ represents the energy consumption of the DNN, and constant $E_{\textrm{budget}}$ is the threshold of the maximum energy budget. Here, both $\phi$ and $\psi$ are potentially nonconvex functions with respect to $\bw_{(u)}$ and $\bS$.

\begin{table*}[htp]
\caption{\small Comparison of algorithms to solve nonconvex optimization problems, where ncvx denotes ``nonconvex'' and const. represents ``constraint''.}
\label{tab.com}
\begin{center}
\begin{small}
\begin{tabular}{llllllll}
\toprule\footnotesize
Algorithm &  Framework &  Const. & Const. Type & Implementation & Complexity \\
\midrule
Proximal ADMM \cite{zhang2020proximal} & inexact & linear & equality & single-loop & $\mathcal{O}(1/\epsilon^2)$ \\
IALM \cite{sahin2019inexact}& inexact  & ncvx  &  equality &  double-loop & $\mathcal{O}(1/\epsilon^4)$ \\
IALM \cite{li21d} & inexact  & ncvx  &  equality &  triple-loop & $\mathcal{O}(1/\epsilon^3)$ \\
IPPP \cite{lin2019inexact} &  penalty  & ncvx  &  inequality  & triple-loop & $\mathcal{O}(1/\epsilon^3)$ \\
IQRC \cite{ma20d}   & primal  & ncvx & inequality & double-loop & $\mathcal{O}(1/\epsilon^3)$\\
\midrule
{\bf GDPA} (This work) &  primal-dual  & ncvx  &  inequality & {\bf single-loop}  &  $\mathcal{O}(1/\epsilon^3)$\\
\bottomrule
\end{tabular}\\
\end{small}
\end{center}
\vskip -0.1in
\end{table*}
\subsection{Related Work}

Solving the nonconvex problems is a long-standing question in machine learning as well as other fields. When the constraints are linear equality constraints or convex constraints, many existing works study different types of algorithms to solve nonconvex objective function optimization problems, such as primal-dual algorithms \cite{hong2016convergence,hajinezhad2019perturbed,zhang2020proximal,wotao21}, inexact proximal accelerated augmented
Lagrangian methods \cite{xu2019iteration,kong2019complexity,melo2020iteration}, trust-region approaches \cite{cartis2011evaluation}, to just name a few.

\noindent{\bf Penalty based Methods}. When the objective function and constraints are both nonconvex, there is a penalty based method, called penalty dual decomposition method (PDD) \cite{pdd}, which 1) directly penalizes the nonconvex constraint function to the objective; 2) solves the problem based on the new surrogate function; 3) then check the feasibility of the obtained solution; 4) if the constraint is not satisfied then increase the penalty parameter and go back to 2). Although it is shown that PDD can solve a wide range of nonconvex constrained problems,  the convergence complexity is not obtained.  Based on the  proximal point method and the quadratic penalty method, an inexact proximal-point penalty (IPPP) method \cite{lin2019inexact} is the first penalty-based method that can converge to the Karush-Kuhn-Tucker (KKT) points for nonconvex objective and constraints problems.

\noindent{\bf Inexact Augmented Lagrangian Method}. The inexact augmented lagrangian method (IALM) is one of the most popular ones that solve the nonconvex optimization problems with constraints. The main idea of this family of algorithms is to add a quadratic term (proximal term) to the augmented Lagrangian function so that the resulting surrogate of the loss function becomes strongly convex and can be solved efficiently by leveraging the existing computationally efficient accelerated first-order methods as an oracle or subroutine. For an equality nonlinear/nonconvex constrained problem, IALM is proposed in \cite{sahin2019inexact,xie2021complexity} with quantifiable convergence rate guarantees to find KKT points of nonlinear constrained problems. For inequality non-convex constrained problems,  inexact quadratically regularized constrained (IQRC) methods \cite{ma20d,gudi19} are proposed recently for solving problem \eqref{eq.org}, which are developed based on the Moreau envelope notation of stationary points/KKT points under a uniform Slater regularity condition.  Although these works are able to solve the nonconvex constrained problems, all of them require double or triple inner loops, which makes the implementation of these algorithms complicated and time-consuming.
A short summary of the existing algorithms for solving nonconvex problems is shown in \figtab{tab.com}.

\noindent{\bf Single-Loop Min-Max Algorithms}. As solving constrained optimization problems can be formulated by searching a minmax equilibrium point of  the augmented Lagrangian function, another line of work on solving nonconvex minmax problem is very related to this framework. Remarkably, the recent developed nonconvex minmax solvers are single-loop algorithms, e.g., gradient descent and ascent (GDA) \cite{lin2020gradient}, hybrid block successive approximation  (HiBSA) method \cite{hibsa}, and smoothed-GDA \cite{zhang2020single}. However, they need the compactness of the dual variable resulting that there is no result that can quantify the constraint violation of the iterates generated by those algorithms, since the dual bound  is essential to measure satisfaction of the solutions as KKT points \cite{jialuo20}.

\subsection{Main Contributions of This Work}
In this work, we propose a single-loop gradient descent and perturbed ascent algorithm (GDPA) by using the idea of designing single-loop nonconvex minmax algorithms to solve nonconvex objective optimization problems with nonconvex constraints. Inspired by the dual perturbation technique \cite{koshal2011multiuser,hajinezhad2019perturbed}, it is shown that GDPA is able to find the KKT points of problem \eqref{eq.org} with provable convergence rate guarantees under mild assumptions.

The main contributions of this work are highlighted as follows
\begin{itemize}
\item {\bf Single-Loop}. To the best of our knowledge, this is the first single-loop algorithm that can find KKT points of nonconvex optimization problems under nonconvex inequality constraints.

\item {\bf Convergence Analysis}. Under a mild regularity condition, we provide the theoretical convergence rate of GDPA to KKT points of problem \eqref{eq.org} in an order of $1/\epsilon^3$, matching the best known rate achieved by double- and triple-loop algorithms.

\item {\bf Applications}. We discuss several possible applications of this class of algorithms with applications to machine learning problems, and give the numerical experimental results to showcase the computational efficiency of the single-loop algorithm compared with the state-of-the-art double-loop or triple-loop methods.
\end{itemize}

\section{Gradient Descent and Perturbed Ascent Algorithm}
\label{submission}
First, we can write down the Lagrangian function of problem \eqref{eq.org} as \cite{nocedal2006numerical}
\begin{equation}\label{eq.lag}
\mathcal{L}(\bx,\blambda)\bydef f(\bx) + \langle g(\bx),\blambda\rangle
\end{equation}
where  non-negative $\blambda\in\mathbb{R}^m_+$ denotes the dual variable (Lagrangian multiplier).

Instead of designing an algorithm based on optimizing the original Lagrangian function, we propose to construct the following perturbed augmented Lagrangian function:

\begin{equation}\label{eq.defF}
F_{\beta}(\bx,\blambda)\!\bydef\! f(\bx) + \frac{\beta}{2}\left\|\!\left[g(\bx)\!+\!\frac{(1-\tau)\blambda}{\beta}\right]_+\!\right\|^2-\frac{\|(1-\tau)\blambda\|^2}{2\beta}
\end{equation}
where $[\bx]_+$ denotes the component-wise nonnegative part of vector $\bx$, $\beta>0$, and constant $\tau\in(0,1)$. Here, perturbation term $\tau$ plays the critical role of ensuring the convergence of the designed algorithm. (Please see \secref{sec.key} for more discussion.)

Next, we consider finding a stationary (quasi-Nash equilibrium) point  \cite{pang2011nonconvex} of the following problem to solve the nonconvex constrained problem \eqref{eq.org}

\begin{equation}\label{eq.minmax}
    \min_{\bx\in\mathcal{X}}\max_{\blambda\ge 0} F_{\beta}(\bx,\blambda).
\end{equation}

To this end, a single-loop GDPA algorithm is proposed as follows:
\begin{subequations}\label{eq.gpda}
\begin{align}
\bx_{r+1}\!=&{}\arg\min_{\bx\in\mathcal{X}}\!\big\langle\!\nabla_{\bx} F_{\beta_r}(\bx_r,\blambda_r), \bx-\bx_r\!\big\rangle\!+\!\frac{1}{2\alpha_r}\|\bx-\bx_r\|^2,\label{eq.updatex0}
\\\notag
\blambda_{r+1}\!=&{}\arg\max_{\blambda\ge 0}\left\langle \frac{1}{1-\tau}\nabla_{\blambda} F_{\beta_r}(\bx_{r+1},\blambda_r),\blambda-\blambda_r\right\rangle
\\
&\qquad\qquad -\frac{1-\tau}{2\beta_r}\|\blambda-\blambda_r\|^2-\frac{\gamma_r}{2}\|\blambda\|^2\label{eq.updatel0}
\end{align}
\end{subequations}
where
\begin{equation}\label{eq.deftau}
\beta_r\gamma_r\bydef\tau\in(0,1),
\end{equation}
and we design that $\beta_r$ is an increasing sequence and subsequently $\gamma_r$ is a decreasing one.

Substituting \eqref{eq.defF} into \eqref{eq.updatex0}  yields
\begin{align}\notag
\bx_{r+1}\!=&\arg\min_{\bx\in\mathcal{X}}\!\big\langle\!\nabla f(\bx_r)
\!+\!\beta_rJ^{\T}(\bx_r)\!\left[g(\bx_r)\!+\!\frac{(1-\tau)\blambda_r}{\beta_r}\right]_+\!\!\!,
\\ & \qquad\qquad \bx-\bx_r\!\big\rangle\!+\!\frac{1}{2\alpha_r}\|\bx-\bx_r\|^2,\label{eq.updatex}
\end{align}
where $J(\bx)$ denotes the Jacobian matrix of function $g(\cdot)$ at point $\bx$, $r$ is the index of iterations, $\alpha_r$ is the step-size of the minimization step.

Regarding the update of the dual variable, it is dependent on the functional constraints satisfaction.  Let
\begin{equation}
\mathcal{S}_r\bydef\left\{i|g_i(\bx_r)+\frac{(1-\tau)[\blambda_r]_i}{\beta_r}>0\right\}, \label{eq.defs}
\end{equation}
where $g_i(\bx)$ denotes the $i$th constraint, and notation $[\bx]_i$ denotes the $i$th entry of vector $\bx$. Then, it is obvious that $g_i(\bx_r)\le 0, i\in\overbar{\mathcal{S}}_r$.
Substituting \eqref{eq.defF} into \eqref{eq.updatel0} results in the following two cases of updating dual variable $\blambda$, which are
\begin{subequations}\label{eq.gpda2}
\begin{align}\notag
[\blambda_{r+1}]_i=&{}\arg\max_{\blambda_i\ge 0}\langle g_i(\bx_{r+1}),\blambda_i-[\blambda_r]_i\rangle-\frac{\gamma_r}{2}\|\blambda_i\|^2
\\
&-\frac{1-\tau}{2\beta_r}\|\blambda_i-[\blambda_r]_i\|^2, \forall i\in\mathcal{S}_r
\\\notag
[\blambda_{r+1}]_i=&{}\arg\max_{\blambda_i\ge 0}\langle -\frac{(1-\tau)}{\beta_r}[\blambda_r]_i,\blambda_i-[\blambda_r]_i\rangle-\frac{\gamma_r}{2}\|\blambda_i\|^2
\\
&-\frac{1-\tau}{2\beta_r}\|\blambda_i-[\blambda_r]_i\|^2,\forall i\in\overbar{\mathcal{S}}_r.
\end{align}
\end{subequations}

It can be seen that the update of primal variable $\bx$ in the minimization step is a standard one, which is optimizing the linearized function $F_{\beta_r}$ at point $(\bx_r,\blambda_r)$ with a proximal term. The perturbed dual update is the key innovation, where the perturbation parameter $\gamma_r$ adds the negative curvature to the maximization problem so that the dual update is well-behaved and easy to analyze. As $\gamma_r$ shrinks, the maximization step reduces to the classic update of the Lagrangian multiplier \cite{bertsekas99,boyd2004convex}.

Note that GDPA is a single-loop algorithm and its updating rules \eqref{eq.gpda},\eqref{eq.updatex},\eqref{eq.gpda2} can be written equivalently in the following form, i.e.,
\begin{subequations}\label{eq.alg}
\begin{align}\notag
\bx_{r+1}=& {} \mathcal{P}_{\mathcal{X}}\big(\bx_r-\alpha_r \big(\nabla f(\bx_r)
\\
&+J^{\T}(\bx_r)\left[(1-\tau)\blambda_r+\beta_rg(\bx_r)\right]_+\big)\big),\label{eq.primal}
\\
[\blambda_{r+1}]_i = & {} \begin{cases}\mathcal{P}_{\ge0}\left((1\!-\!\tau)[\blambda_r]_i\!+\!\beta_r g_i(\bx_{r+1})\right), i\in\mathcal{S}_r
\\
0, i\in\overbar{\mathcal{S}}_r
\end{cases} \label{eq.dualupated}
\end{align}
\end{subequations}
where $\mathcal{P}_{\mathcal{X}}$ denotes the projection of iterates to the feasible set, $\mathcal{P}_{\ge 0}\bydef[]_+$ is the component-wise nonnegative projection operator, and $\beta_r$ serves as the step-size of the maximization step.

\remark  If subproblem in \eqref{eq.updatex} is an unconstrained one, such as  a machine learning model parametrized by a neural network, then $\mathcal{X}=\mathbb{R}^d$, where $\bx$ denotes the weights. In this case, the update \eqref{eq.primal} in GDPA reduces to $\bx_{r+1}=\bx_r-\alpha_r (\nabla f(\bx_r)+J^{\T}(\bx_r)[(1-\tau)\blambda_r+\beta_rg(\bx_r)]_+)$.

\section{Theoretical Guarantees}

Before showing the theoretical convergence rate result of GDPA, we first make the following blanket assumption for problem \eqref{eq.org}.
\subsection{Assumptions}
\begin{assumption}\label{ass.lif}
(Lipschitz continuity of function $f(\bx)$)
We assume that $f(\bx)$ is smooth and has gradient Lipschitz continuity with constant $L_f$, i.e.,
$\|\nabla f(\bx) - \nabla f(\by)\|\le L_f\|\bx-\by\|,\forall \bx,\by\in\mathcal{X}\subseteq\mathbb{R}^d$.
\end{assumption}
\begin{assumption}\label{ass.lig}
(Lipschitz continuity of function $g(\bx)$)
function $g(\bx)$ has function Lipschitz continuity with constant $L_g$, i.e., $\|g(\bx)-g(\by)\|\le L_g\|\bx-\by\|,\forall \bx,\by\in\mathcal{X}\subseteq\mathbb{R}^d$, and the Jacobian function of $g(\bx)$ is Lipschitz continuous with constant $L_J$, i.e., $\|J(\bx)-J(\by)\|\le L_J\|\bx-\by\|,\forall \bx,\by\in\mathcal{X}\subseteq\mathbb{R}^d$.
\end{assumption}
\begin{assumption}\label{ass.bfg}
(Boundedness of function $f(\bx)$)
Further, we assume that the lower bound of $f(\bx)$ is $f^{\star}$, i.e.,  $\min_{\bx\in\mathcal{X}}f(\bx)>f^{\star}>-\infty$, and the upper bound of the size of the gradient of $f(\bx)$ with respect to variable $\bx$ is $M$, i.e., $\|\nabla_{\bx} f(\bx)\|\le M,\forall \bx\in\mathcal{X}\subseteq\mathbb{R}^d$.
\end{assumption}
\begin{assumption}\label{ass.bd} (Boundedness of function $g(\bx)$)
Assume that the size of function $g_+(\bx)$ is upper bounded by $G$, i.e., $\|g_+(\bx)\|^2\le G,\forall \bx\in\mathcal{X}\subseteq\mathbb{R}^d$, and the Jacobian function of $g(\bx)$ is upper bounded by $U_J$, i.e., $\|J(\bx)\|\le U_J,\forall \bx\in\mathcal{X}\subseteq\mathbb{R}^d$.
\end{assumption}

All the above assumptions are based on the functions themselves and standard in analyzing the convergence of algorithms. Alternatively, we can assume compactness  of feasible set $\mathcal{X}$ as follows.
\begin{assumption}\label{ass.compact}
Assume that $\mathcal{X}$ is convex and compact.
\end{assumption}
Previous works \cite{sahin2019inexact,li21d,ma20d,lin2019inexact} also assume the compactness of feasible sets, which imply \asref{ass.lif} to \asref{ass.bd} directly for smooth functions.

Due to the nonconvexity of the constraints, we need the following regularity condition to ensuring the feasibility of solutions with the functional nonconvex constraints.

\begin{assumption}\label{ass.re}(Regularity condition) We assume that there exists a constant $\sigma>0$  such that
\begin{equation}\label{eq.re}
\sigma \|g_+(\bx)\| \le \textrm{dist}\left(J^{\T}(\bx) g_+(\bx),-\mathcal{N}_{\mathcal{X}}(\bx)\right)
\end{equation}
where $\textrm{dist}(\bx,\mathcal{X})\bydef\min_{\by\in\mathcal{X}}\|\bx-\by\|$ and $\mathcal{N}_{\mathcal{X}}(\bx)$ denotes the normal cone of feasible set $\mathcal{X}$ at point $\bx$.
\end{assumption}

\remark When $\mathcal{X}=\mathbb{R}^d$, condition \eqref{eq.re} reduces to $\sigma \|g_+(\bx)\| \le \|J^{\T}(\bx) g_+(\bx)\|$.

The regularity qualification is required in solving constrained optimization problems, e.g., Slater's condition, linear independence constraint qualification (LICQ), and so on, which is different from the unconstrained problems or the constrained one with closed-form projection operators. Condition \eqref{eq.re} is a standard one and has been adopted to analyzing the convergence of iterative algorithms with nonconvex functional constrained problems \cite{sahin2019inexact,li21d,lin2019inexact}.

\vspace{-5pt}
\subsection{Theoretical Guarantees}

\noindent{\bf Convergence Rate}. Next, we use the following measure to quantify the optimality of the iterates generated by GDPA:
\begin{equation}\label{eq.defopt}
\mathcal{G}(\bx,\blambda) \bydef\left[\begin{matrix} \frac{1}{\alpha}\left[\bx - \textrm{proj}_{\mathcal{X}}(\bx-\alpha\nabla_{\bx}\mathcal{L}(\bx,\blambda))\right]\\  \frac{1}{\beta}\left[\blambda - \textrm{proj}_{\ge 0}(\blambda+\beta\nabla_{\blambda}\mathcal{L}(\bx,\blambda))\right]\end{matrix}\right],
\end{equation}
where $\alpha,\beta>0$. This optimality gap has been widely used in the theoretical analysis of nonconvex algorithms for solving constrained optimization problems \cite{hong2016convergence} and nonconvex minmax problems \cite{hibsa}. Besides, we need feasibility and slackness conditions in quantifying constraint satisfaction. Together with \eqref{eq.defopt},  approximate stationarity conditions are given as follows.
\begin{definition}\label{de.st}($\epsilon$-approximate Stationary Points)  A point $\bx$ is called an $\epsilon$-approximate stationary  point of problem \eqref{eq.org} if there is a $\blambda\ge0$ such that
\begin{equation}\label{eq.optst}
\|\mathcal{G}(\bx,\blambda)\|\le \epsilon.
\end{equation}
\end{definition}
Then, we provide the main theorem of the convergence rate of GDPA as follows.

\begin{theorem}\label{th.main1}
Suppose that \asref{ass.lif}-\asref{ass.bd}  (or \asref{ass.compact} ) and \asref{ass.re} hold and iterates $\{\bx_r,\blambda_r,\forall r\ge0\}$ are generated by GDPA.  When the step-sizes are chosen as
\begin{equation}\label{eq.spcon}
\alpha_r\sim\gamma_r\sim\frac{1}{\beta_r}\sim\mathcal{O}\left(\frac{1}{r^{\frac{1}{3}}}\right), \quad \gamma_r\beta_r=\tau,
\end{equation}
and $\max\{1-\sigma/\sqrt{66U^2_J+\sigma^2},1/2\}<\tau<1$, there exist constants $K_1,K_2,K_3$ such that
\begin{subequations}
\begin{align}
&\|\mathcal{G}(\overbar{\bx}_{T(\epsilon)},\overbar{\blambda}_{T(\epsilon)})\|^2\le \frac{K_1}{T(\epsilon)^{\frac{2}{3}}},
\\
&\|g_+(\overbar{\bx}_{T(\epsilon)})\|^2\le\frac{K_2}{T(\epsilon)^{\frac{2}{3}}},
\\
& \sum^m_{i=1}\left|[\overbar{\blambda}_{T(\epsilon)}]_i g_i(\overbar{\bx}_{T(\epsilon)})\right|\le\frac{K_3}{T(\epsilon)^{\frac{1}{3}}}
\end{align}
\end{subequations}
where
\begin{multline}\label{eq.gapc}
T(\epsilon)\bydef\min\{r| \|g_+(\bx_{r+1})\|^2\le\epsilon, r> 1\},
\end{multline}
 and outputs $\overbar{\bx}_{T(\epsilon)}\bydef(\sum^{T(\epsilon)}_{r=1}1/\beta_r)^{-1}\bx_r/\beta_{r}$, $\overbar{\blambda}_{T(\epsilon)}\bydef(\sum^{T(\epsilon)}_{r=1}1/\beta_r)^{-1}\blambda_r/\beta_r$.
\end{theorem}
Due to the space limit, all the detailed proofs in this paper are relegated to the supplemental material.

\remark When $\mathcal{X}=\mathbb{R}^d$, then according to \eqref{eq.lag} and \eqref{eq.defopt} it is concluded immediately that $\|\mathcal{G}(\bx,\blambda)\|\le\epsilon$ implies $\|\nabla f(\bx)+\sum_i\blambda_i\nabla g_i(\bx)\|\le\epsilon$.

\noindent{\bf KKT Points} Based on the notation of quasi-Nash equilibrium points, we have obtained the convergence rate of GDPA to $\epsilon$-approximate stationary points under with the constraint satisfaction and slackness condition. While in the classic constrained optimization theory, the convergence results are established on the KKT conditions \cite{bertsekas99}, which are given as follows.
\begin{definition}\label{de.kkt}($\epsilon$-approximate KKT)  A point $\bx$ is called an $\epsilon$-approximate KKT  point of problem \eqref{eq.org} if there is a $\blambda\ge0$ such that
\begin{subequations}\label{eq.kkt}
\begin{align}
&\textrm{dist}\left(\nabla f(\bx)+ \sum^m_{i=1}\blambda_i\nabla g_i(\bx), -\mathcal{N}_{\mathcal{X}}(\bx)\right)\le \epsilon, \label{eq.kkts1}
 \\
& \|g_+(\bx)\|\le\epsilon, \label{eq.kktc1}
 \\
 & \sum^m_{i=1}|\blambda_i g_i(\bx)|\le\epsilon.\label{eq.kktsl1}
\end{align}
\end{subequations}
\end{definition}

When $\epsilon=0$, then $\|\mathcal{G}(\bx,\blambda)\|=0$ or more precisely $\|1/\alpha [\bx - \textrm{proj}_{\mathcal{X}}(\bx-\alpha \nabla_{\bx}\mathcal{L}(\bx,\blambda))]\|=0$ implies ${\textrm{dist}}\left(\nabla f(\bx)+ \sum^m_{i=1}\blambda_i\nabla g_i(\bx), -\mathcal{N}_{\mathcal{X}}(\bx)\right)=0$. In this work, we also provide the following proposition to show the relation between the approximate stationary points \eqref{eq.optst} and \eqref{eq.kkts1}.
\begin{proposition}\label{pro.eq}
When $\mathcal{X}$ is convex, the stationarity of the approximate saddle points defined by \eqref{eq.optst} is a sufficient condition for \eqref{eq.kkts1}, namely if a point $(\bx,\blambda)$  satisfies $\|1/\alpha [\bx - \textrm{proj}_{\mathcal{X}}(\bx-\alpha \nabla_{\bx}\mathcal{L}(\bx,\blambda))]\|\le\epsilon, \alpha>0$, then it also satisfies \eqref{eq.kkts1}.
\end{proposition}
To show the above result, classic Farkas Lemma, (i.e., Lemma 12.4 in \cite{nocedal2006numerical}) is not applicable since the definition of the separating hyperplane is built on the exact stationary points, i.e., the case where $\epsilon=0$. However, here we need a notation of approximate stationary points of being able to quantify the convergence rate. Therefore, we give the following variant of approximate Farkas lemma, which bridges the connection between \eqref{eq.optst} and \eqref{eq.kkts1}.
\def\ck{{\mathcal{K}}}
\def\cf{\emph{\textmd{F}}}
\begin{lemma}\label{le.farkas}
(Approximate Farkas Lemma). Let the cone $\ck$ be defined as $\ck=\{\bB\by+\bC\bw|\by\ge 0\}$ where $\bB\in\mathbb{R}^{d\times m},\bC\in\mathbb{R}^{d\times p}, \by\in\mathbb{R}^m,\bw\in\mathbb{R}^p$. Given any vector $\bg\in\mathbb{R}^d$ and $0<\epsilon<1$, we have either  $(\bg+\Delta)\in\ck$ where $\Delta\in\ck$ and $\|\Delta\|\le\epsilon$ or that there exists a $\bd\in\mathbb{R}^d$ where $\|\bd\|=1$ satisfying
\begin{subequations}\label{eq.pf}
\begin{align}
&\bg^{\T}\bd<-\epsilon,\label{eq.pfa}
\\
&\bB^{\T}\bd\ge0,\label{eq.pfb}
\\
&\bC^{\T}\bd=0, \label{eq.pfc}
\end{align}
\end{subequations}
but not both.
\end{lemma}
\remark The main difference between \leref{le.farkas} and the classic one is that the size of separating hyperplane $\bd$ is bounded, otherwise, the error tolerance defined based on the inner product, i.e., \eqref{eq.pfa}, is meaningless.

 Combining \prref{pro.eq}, \thref{th.main1} claims that the proposed GDPA can find an approximate KKT point at a rate of $\mathcal{O}(1/\epsilon^3)$.
\begin{corollary}\label{th.main2}
Suppose that \asref{ass.lif}-\asref{ass.bd}  (or \asref{ass.compact} ) and \asref{ass.re} hold and iterates $\{\bx_r,\blambda_r,\forall r\ge0\}$ are generated by GDPA.  When the step-sizes are chosen as $\alpha_r\sim1/\beta_r\sim\mathcal{O}(1/r^{1/3})$ and $\max\{1-\sigma/\sqrt{66U^2_J+\sigma^2},1/2\}<\gamma_r\beta_r=\tau<1$, then the outputs of GDPA $\overbar{\bx}_{T(\epsilon)},\overbar{\blambda}_{T(\epsilon)}$ converge to an $\epsilon$-approximate KKT point satisfying
\begin{align}\notag
&\textrm{dist}\left(\!\nabla f(\overbar{\bx}_{T(\epsilon)})\!+\!\sum^m_{i=1}[\overbar{\blambda}_{T(\epsilon)}]_i\nabla g_i(\overbar{\bx}_{T(\epsilon)}), \!-\mathcal{N}_{\mathcal{X}}(\bx)\right)\le \epsilon,
\\
&\|g_+(\overbar{\bx}_{T(\epsilon)})\|\le\epsilon, \quad \sum^m_{i=1}|[\blambda_{T(\epsilon)}]_i g_i(\overbar{\bx}_{T(\epsilon)})|\le\epsilon,\label{eq.optst2}
\end{align}
in the number of $\mathcal{O}(1/\epsilon^3)$ iterations.
\end{corollary}

\subsection{Convergence Analysis}\label{sec.key}

The following theoretical results showcase the main ideas of measuring the convergence rate of GDPA to KKT points, where the key step is to derive the upper bound of the size of dual variable $\blambda$.

\noindent{\bf Perturbation}. The first result is that we found that after adding the perturbation the optimality gap defined in \eqref{eq.defopt} will be upper bounded the successive difference between primal and dual variables plus a perturbed term $\gamma^2_r\|\blambda_{r+1}\|^2$.
\begin{lemma}\label{le.opt}
Suppose that \asref{ass.lif}--\asref{ass.bfg} hold. If the iterates $\{\bx_r,\blambda_r,\forall r\}$ are generated by GDPA, where the step-sizes are chosen according to \eqref{eq.spcon}, then we have
\begin{align}
\notag
&\|\mathcal{G} (\bx_r, \blambda_r)\|^2
\\\notag
\le & 4\left(\left(\frac{3}{\alpha_r}+2L_g\right)^2+2U^2_JL^2_g\beta^2_r\right)\|\bx_{r+1}-\bx_r\|^2
\\
 +& 4\left(\frac{(3+2\tau)^2}{\beta^2_r}\!+\!3\right)\!\!\|\blambda_{r+1}\!-\!\blambda_r\|^2\!+\!16\gamma^2_r\|\blambda_{r+1}\|^2.\label{eq.ccov1}
\end{align}
\end{lemma}
Under \asref{ass.lif} to \asref{ass.bfg}, we can know an upper bound of the first term on the right-hand side (RHS) of \eqref{eq.ccov1} by applying gradient Lipchitz continuity of $F_{\beta}(\bx,\blambda)$, and an upper  bound of the second term on RHS of \eqref{eq.ccov1} by quantifying the strong concavity of $F_{\beta}(\bx,\blambda)$. Also, these upper bounds can be written as the difference between $F_{\beta_{r+1}}(\bx_{r+1},\blambda_{r+1})$ and $F_{\beta_r}(\bx_r,\blambda_r)$ in general. The major challenge is to get an upper bound of the last term in RHS of  \eqref{eq.ccov1}. From \eqref{eq.dualupated}, it is not hard to show that the upper bound of $\blambda$ is $\mathcal{O}(\beta^2_r)$ under \asref{ass.bd}. However, it is implied from \eqref{eq.deftau} that the last term on RHS of \eqref{eq.ccov1} is only upper bounded by a constant, giving rise to a constant error in the optimality gap.

\noindent{\bf Boundedness of Dual Variable}. The main idea of having a sharper upper bound of dual variable $\blambda$ or sum of the last term on RHS of \eqref{eq.ccov1} is using the regularity condition \eqref{eq.re}, which provides certain contraction property on the size of $\blambda$ up to some terms. The detailed claim is given as follows.
\begin{lemma}\label{le.rege}
Suppose that \asref{ass.lif}-\asref{ass.bd}  (or \asref{ass.compact} ) and \asref{ass.re} hold. Let the active set at the $r$th iteration be
\begin{equation}
 \mathcal{A}_{r}\bydef\{i|g_i(\bx_{r})>0\}.
\end{equation}
If the iterates $\{\bx_r,\blambda_r,\forall r\}$ are generated by GDPA, where the step-sizes are chosen according to \eqref{eq.spcon}, then we have
\begin{multline}\label{eq.itofl}
\sigma^2\|\blambda_{r+1}^g\|^2\le  16M^2+ \frac{36}{\alpha^2_r}\|\bx_{r+1}-\bx_r\|^2
\\+ 64U^2_J\|\blambda_{r+1}-\blambda_r\|^2
+64(1-\tau)^2U^2_J\|\blambda_r\|^2
\end{multline}
where $\blambda^g_{r+1}$ is a vector whose $i$th entry is
\begin{equation}\label{eq.defl}
 [\blambda^g_{r+1}]_i = \begin{cases} [\blambda_{r+1}]_i, & \textrm{if}\quad i\in\mathcal{A}_{r+1}\cap\mathcal{S}_r; \\ 0, & \textrm{otherwise}.\end{cases}.
\end{equation}
\end{lemma}
Due to the nonnegativity of $\blambda_r$ and $\tau\in(0,1)$, we know that when $g_i(\bx_{r+1})\le 0$, we have the following contraction property from \eqref{eq.dualupated}
\begin{equation}\label{eq.shk}
[\blambda_{r+1}]_i\le(1-\tau)[\blambda_r]_i, \quad \forall i\in\overbar{\mathcal{A}}_{r+1}\cap\mathcal{S}_r.
\end{equation}
Combing \eqref{eq.itofl} and \eqref{eq.shk}, we can have a contraction property of the size of the dual variable with some additional terms when $\tau$ is close to 1, which demonstrates the importance and reason of adding the perturbation term in the dual update. Based on the fact, we are able to have the following encouraging result.
\begin{lemma}\label{le.sumd}
Suppose that \asref{ass.lif}-\asref{ass.bd}  (or \asref{ass.compact} ) and \asref{ass.re} hold. If the iterates $\{\bx_r,\blambda_r,\forall r\}$ are generated by GDPA, where the step-sizes are chosen according to \eqref{eq.spcon}, then we have
\begin{equation}\label{eq.bdsum}
\sum^T_{r=1}\alpha^2_r\|\blambda_r\|^2\sim\mathcal{O}(T^{1/3}),
\end{equation}
where $T$ denotes the total number of iterations.
\end{lemma}

Note that $\alpha_r\sim\gamma_r$,  from \eqref{eq.spcon} we have  $\sum^T_{r=1}\alpha^2_r\|\blambda_r\|^2\sim\mathcal{O}(T^{1/3})$. If we only use \asref{ass.bd} to get the upper bound of $\gamma^2_r\|\blambda_r\|^2$, we will have $\sum^T_{r=1}\alpha^2_r\|\blambda_r\|^2\sim\mathcal{O}(T)$, resulting in failure to show the convergence GDPA to KKT points.

\begin{figure*}[t]
\centering
\vspace{-0.2cm}
\subfigure[Stationarity of the primal variable]{
\label{fig:aconi}
\includegraphics[width=.38\linewidth]{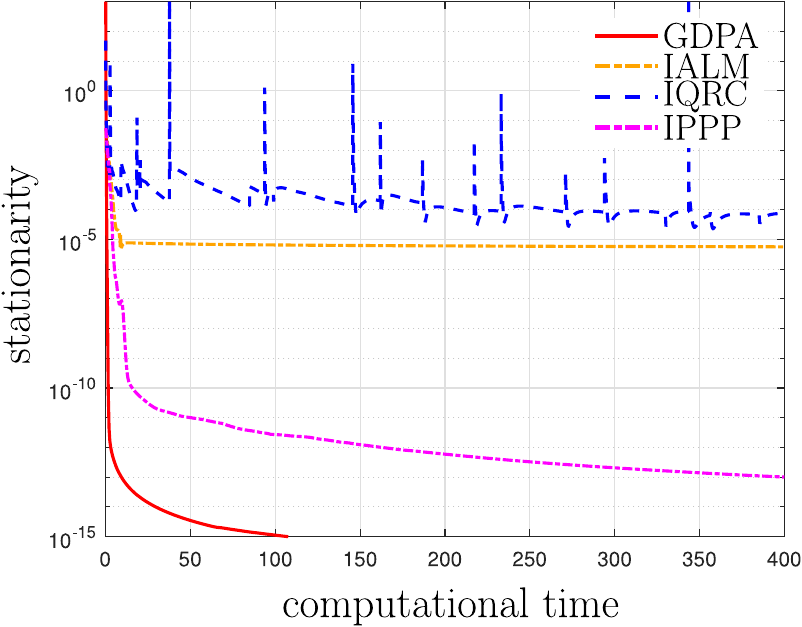}}
\hspace{0.4in}
\subfigure[{Feasibility}]{
\label{fig:acont}
\includegraphics[width=.38\linewidth]{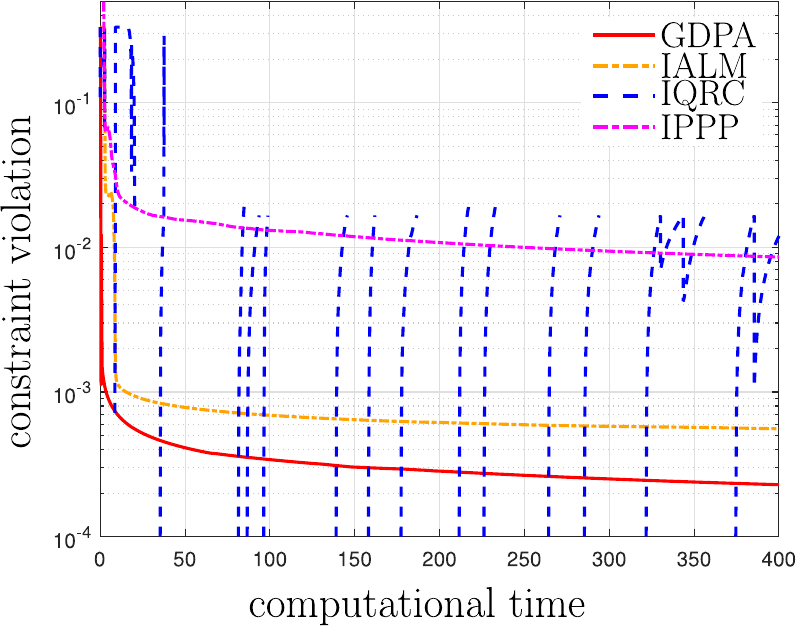}}
\caption{\footnotesize Computational time comparison of GDPA, IALM, IQRC, IPPP.}
\vspace{-4pt}
\label{fig:acon}
\end{figure*}
\begin{figure*}[htp]
\centering
\vspace{-0.1cm}
\subfigure[Stationarity of the primal variable]{
\label{fig:coni}
\includegraphics[width=.38\linewidth]{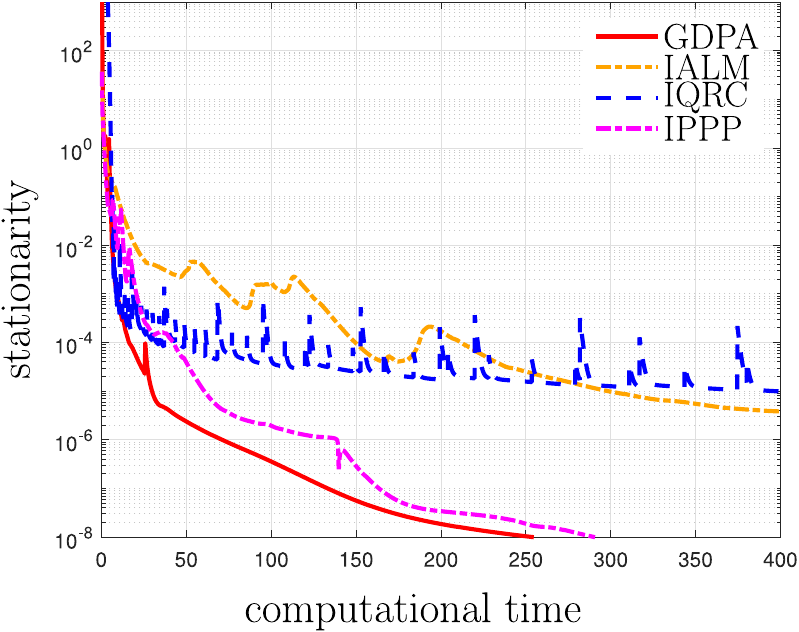}}
\hspace{0.4in}
\subfigure[{Feasibility}]{
\label{fig:cont}
\includegraphics[width=.38\linewidth]{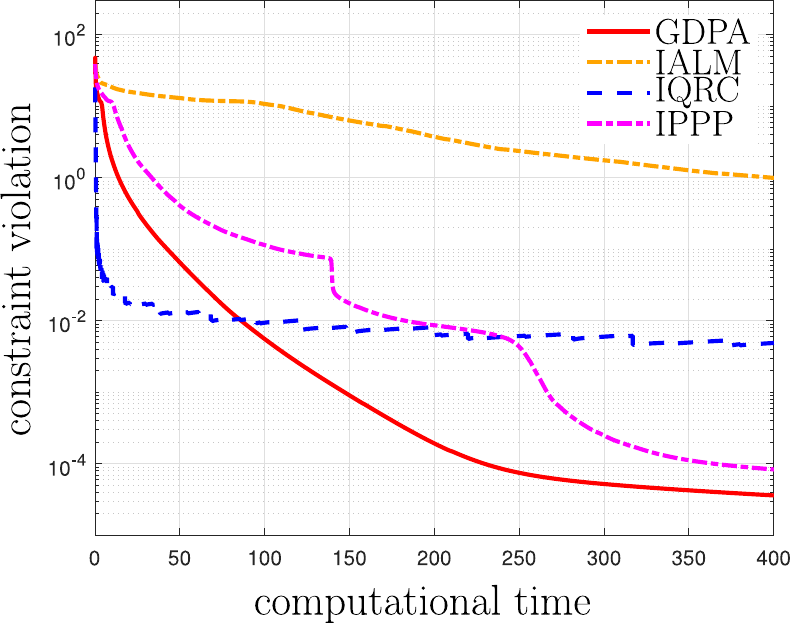}}
\caption{\footnotesize Computational time comparison of GDPA, IALM, IQRC, IPPP.}
\label{fig:opt}
\vspace{-8pt}
\end{figure*}

\section{Discussion}
\subsection{Regularity Condition}
The regularity condition \eqref{eq.re} has been proven to hold for many practical problems, such as a certain kind of mNPC problems \cite{lin2019inexact}, clustering and basis pursuit problems \cite{sahin2019inexact} under some mild conditions on initializations. Also, the variant of this condition for functional linear equality constraints holds automatically for either feasible set $\mathcal{X}$ is a ball constraint or a compact polyhedral one \cite{li21d}.

\vspace{-5pt}
\subsection{Comparison with Existing Works}
\noindent{\bf Penalty based Method}. Comparing with the penalty-based methods \cite{pdd,lin2019inexact}, primal-dual methods can deal with the multiple nonconvex constraints more flexibly, since the dual variable takes the variety of constraints automatically. Numerically, the primal-dual type of algorithms converges, in general, faster than the penalty based method. But theoretical analysis of the penalty-based method is easier and more accessible due to the simplicity of the algorithms. Besides, the regularity condition \eqref{eq.re} used in this work is the same as \cite{lin2019inexact}, and the convergence rate achieved by IPPP is $\mathcal{O}(1/\epsilon^3)$.

\noindent{\bf Primal Based Methods}. IQRC \cite{ma20d} is designed based on a uniform regularity condition, which is different from \eqref{eq.re} and  not easily verified. The inner loop of IQRC is realized by an accelerated gradient descent method, where the overall iteration complexity is $\mathcal{O}(1/\epsilon^3)$ for solving problem \eqref{eq.org}. The convergence analysis is neat due to the design of the algorithm, which focuses on optimizing the constraints first and then switching to optimize the objective function after the constraints are satisfied.

\noindent{\bf Inexact Augmented Lagrangian Method}. IALM type of algorithms, e.g., \cite{sahin2019inexact, li21d}  is the most popular one in solving problems with equality nonconvex constraints, where the convergence rate of IALM to find the KKT points of problem \eqref{eq.org} is  $\mathcal{O}(1/\epsilon^4)$  \cite{sahin2019inexact} or $\mathcal{O}(1/\epsilon^3)$ \cite{li21d} under an almost identical regularity condition as \eqref{eq.re}. The main issue of this class of algorithms is that they use a rather small step-size for the dual variable update, and the performance of IALM varies significantly based on the different problems due to the nested structure in terms of loops.

Overall, all these algorithms are double-loop or triple-loop ones, and they all use the accelerated methods to solve the inner loop problems so that a fast theoretical convergence rate is obtained. As there are more than one inner loops in these algorithms, their implementations will involve tunings of multiple hyperparameters, while GDPA only needs to adjust the two step-sizes and one hyperparameter (which is quite insensitive to the convergence behavior of GDPA numerically). Our theoretical analysis of GDPA is also unique, which is established on a new way of showing the boundedness of dual variables.

\vspace{-5pt}
\section{Numerical Experiments}
\vspace{-5pt}

\noindent{\bf mNPC problem}. We compare the convergence performance of the proposed GDPA algorithm with other existing ones, including IALM \cite{sahin2019inexact,li21d}, IQRC \cite{ma20d}, and IPPP \cite{lin2019inexact}, on a mNPC problem. We divide the MNIST dataset \cite{lecun1998gradient} into $m+1$ parts according to the classes of handwritten digits, and consider identifying digit 1 as the prioritized learning class while the other $m$ digits as the secondary ones and each digit as one task. The loss functions of constraints are $\sum_i\phi_j(\bx_{(i)})\le r_j$, where $\phi_j(\bx_{(i)})$ is just a sigmoid function $1/(1+\exp(z))$ as used in \cite{lin2019inexact,ma20d}, $z=(\bx_{(j)}-\bx_{(i)})^{\T}\xi_j,j\neq i$, and the objective loss function is  $\frac{\lambda}{2}\sum_{i}\|\bx_{(i)}\|+\sum_i\phi_1(\bx_{(i)})$ for $\lambda>0$. Also, the input dimension is $784$ (i.e., an image size of $28\times 28$). Since the classification problem on the MNIST dataset is not a hard one in general, we add some random noise at each pixel on the images, where each entry of the noise follows the \emph{i.i.d} Gaussian distribution, so that the constraints are not easily satisfied.  More detailed settings of this numerical experiment can be found in the \secref{sec.addn}. It can be observed in \figref{fig:acon}  that our proposed GDPA converges faster in orders of magnitude compared with the other benchmarks in terms of computational time.  All the IAML type of algorithms with being able to deal with nonconvex constraints is developed based on the functional equality nonconvex constraints, so we add a nonnegative slack variable to reformulate the inequality constrained problem \eqref{eq.org} as an equality one. More detailed settings of this numerical experiment can be found in the \secref{sec.addn}.

\noindent{\bf Neural Nets Training with Budget Constraints}. We also test these algorithms on a training problem with some accuracy budget for fair learning problems, i.e.,
\vspace{-3pt}
\begin{equation}\notag
\min_{\bW}\; \ell_1(\bW), \quad\textrm{subject to}\quad \ell_i(\bW)\le 1, i=2,\ldots, m+1
\end{equation}
where we also use the MNIST dataset and again split the dataset as $m+1$ parts based on the class of digits, $\ell_1(\bW)$ denotes the loss of training the neural net on digit 1, $\ell_i(\bW)$ for digits $2,3,4,5,6$ with $m=5$. The goal of this problem is to train the neural nets on a prioritized dataset with limited accuracy loss on the other ones.

The results are shown in \figref{fig:opt}. It can be observed that the proposed GDPA converges faster than the rest of the methods in terms of computing time, showing the computational efficiency of single-loop algorithms compared with the double-loop or triple-loop ones. Also, we can see that IQRC shows a faster convergence rate in terms of constraint satisfaction since it first optimizes the constraints and then back to minimize the objective function values if the constraint violation has been achieved to a small predefined error. It is worth noting that in this case IALM performs worse than IPPP.

\vspace{-10pt}
\section{Concluding Remark}
\vspace{-6pt}

In this work, we proposed the first single-loop algorithm for solving general nonconvex optimization problems with functional non-convex constraints. Under a mild regularization condition, we show that the proposed GDPA is able to converge to KKT points of this class of non-convex problems at a rate of  $\mathcal{O}(1/\epsilon^3)$. To the best of knowledge, this is the first theoretical result that a single-loop gradient primal-dual (i.e., gradient descent and ascent) algorithm can solve the nonconvex functional constrained problems with the same provable convergence rate guarantees as the double- and/or triple-loop algorithms.

\section{Acknowledgement}
The author would like to thank Professor Qihang Lin for providing feedback on this work \cite{lu2022single}.

\bibliographystyle{abbrvnat}
\bibliographystyle{icml2022}
\bibliographystyle{IEEEbib}

\bibliography{refs}

%

\appendix
\onecolumn

\section{Preliminaries}
Before showing the detailed derivations of the lemmas and theorems, we first list the following notations and inequalities used in the proofs.
\subsection{Notations}
The following facts would help understand the relation between the satisfaction of the functional constraints and the perturbed augmented Lagrangian function.

\begin{enumerate}

\item Expression of $F_{\beta_r}(\bx,\blambda)$

From \eqref{eq.defF}, we observe that  when $\left[g(\bx)+\frac{(1-\tau)\blambda}{\beta_r}\right]_+>0$, then
\begin{equation}\label{eq.fb2}
F_{\beta_r}(\bx,\blambda)= f(\bx)+\langle (1-\tau)\blambda, g(\bx)\rangle +\frac{\beta_r}{2}\|g(\bx)\|^2.
\end{equation}
When $\left[g(\bx)+\frac{(1-\tau)\blambda}{\beta_r}\right]_+\le0$, then we have
\begin{equation}\label{eq.fbl}
F_{\beta_r}(\bx,\blambda)= f(\bx)-\frac{\|(1-\tau)\blambda\|^2}{2\beta_r}.
\end{equation}

Towards this end,  recall
\begin{equation}\notag
\mathcal{S}_r\bydef\left\{i|\left[g(\bx_r)+\frac{(1-\tau)\blambda_r}{\beta_r}\right]_i>0\right\}, 
\end{equation}
then we have
\begin{equation}\label{eq.fbo}
F_{\beta_r}(\bx,\blambda)=f(\bx)+\sum_{i\in\mathcal{S}_r}\left((1-\tau)\blambda_i g_i(\bx)+\frac{\beta_r}{2}g^2_i(\bx)\right) - \sum_{j\in\overbar{\mathcal{S}}_r}\frac{(1-\tau)^2\blambda^2_j}{2\beta_r}
\end{equation}
where $[\bx]_i$ denotes the $i$th entry of vector $\bx$, $\overbar{\mathcal{S}}_r$ represents the complement of set $\mathcal{S}_r$.

From \eqref{eq.defF}, we know that the gradient of $F_{\beta_r}(\bx,\blambda)$ with respect to $\bx$ is 
\begin{align}
\nabla_{\bx} F_{\beta_r}(\bx,\blambda)\bydef&\nabla f(\bx)+J^{\T}(\bx)\left[(1-\tau)\blambda+\beta_rg(\bx)\right]_+. \label{eq.as}
\end{align}

\item Dual variable: based on the constraints satisfaction, we split the corresponding dual variables as two parts, i.e., $\blambda'_r$ and $\bblambda_r$:
\begin{equation}\label{eq.dre}
 [\blambda'_{r}]_i = \begin{cases} [\blambda_{r}]_i, & \textrm{if}\quad i\in\mathcal{S}_r; \\ 0, & \textrm{otherwise}.\end{cases},\quad
 [\bblambda_r]_i = \begin{cases} [\blambda_r]_i, & \textrm{if}\quad i\in\overbar{\mathcal{S}}_r; \\ 0, & \textrm{otherwise}.\end{cases},
\end{equation}
and we have $|\mathcal{S}_r|+|\overbar{\mathcal{S}}_r|=m$. Also, we define
\begin{equation}
 [\blambda''_{r+1}]_i \bydef \begin{cases} [\blambda_{r+1}]_i, & \textrm{if}\quad i\in\mathcal{S}_r; \\ 0, & \textrm{otherwise}.\end{cases}, \quad  [\bblambda'_{r+1}]_i \bydef \begin{cases} [\blambda_{r+1}]_i, & \textrm{if}\quad i\in\overbar{\mathcal{S}}_r; \\ 0, & \textrm{otherwise}.\end{cases}.
\end{equation}

Dual variable $\blambda^g_{r+1}$ and corresponding $\blambda'^g_{r}$ are defined based on the constraints satisfaction at $\bx_r$ and $\bx_{r+1}$:
\begin{equation}
 [\blambda^g_{r+1}]_i = \begin{cases} [\blambda_{r+1}]_i, & \textrm{if}\quad i\in\mathcal{A}_{r+1}\cap\mathcal{S}_r; \\ 0, & \textrm{otherwise}.\end{cases},\quad
 [\blambda'^g_r]_i = \begin{cases} [\blambda_r]_i, & \textrm{if}\quad i\in\mathcal{A}_{r+1}\cap\mathcal{S}_r; \\ 0, & \textrm{otherwise}.\end{cases}.
\end{equation}

Dual variable $\overbar{\blambda}'^g_{r+1}$ and corresponding $\overbar{\blambda}'^g_{r}$ are defined based on the constraints satisfaction at $\bx_r$ and  $\bx_{r+1}$.
\begin{equation}
 [\overbar{\blambda}'^g_{r+1}]_i = \begin{cases} [\blambda_{r+1}]_i, & \textrm{if}\quad i\in\mathcal{A}_{r+1}\cap\overbar{\mathcal{S}}_r; \\ 0, & \textrm{otherwise}.\end{cases},\quad
  [\overbar{\blambda}''^g_{r}]_i = \begin{cases} [\blambda_{r}]_i, & \textrm{if}\quad i\in\mathcal{A}_{r+1}\cap\overbar{\mathcal{S}}_r; \\ 0, & \textrm{otherwise}.\end{cases}.
\end{equation}

\item Functional constraints:

Similarly, we define functional constraints $g_1(\bx_{r+1})$ and $g'(\bx_{r+1})$ based on the constraints satisfaction as well, i.e.,
\begin{equation}
 [g_1(\bx_{r+1})]_i = \begin{cases} g_i(\bx_{r+1}), & \textrm{if}\quad i\in \mathcal{S}_r; \\ 0, & \textrm{otherwise}.\end{cases},\quad
 [g'(\bx_{r+1})]_i = \begin{cases} g_i(\bx_{r+1}), & \textrm{if}\quad i\in\mathcal{A}_{r+1}\cap\mathcal{S}_r; \\ 0, & \textrm{otherwise}.\end{cases}\label{eq.deffgp}
 \end{equation}

Also, functional constraints $\overbar{g}(\bx_{r+1})$ and $\overbar{g}'(\bx_{r})$ are defined as follows:
\begin{equation}
 [\overbar{g}(\bx_{r+1})]_i = \begin{cases} g_i(\bx_{r+1}), & \textrm{if}\quad i\in\mathcal{A}_{r+1}\cap\overbar{\mathcal{S}}_r; \\ 0, & \textrm{otherwise}.\end{cases},\quad
 [\overbar{g}'(\bx_{r})]_i = \begin{cases} g_i(\bx_{r}), & \textrm{if}\quad i\in\mathcal{A}_{r+1}\cap\overbar{\mathcal{S}}_r; \\ 0, & \textrm{otherwise}.\end{cases}.
\end{equation}
\end{enumerate}

\subsection{Inequalities}

\begin{enumerate}

\item Quadrilateral identity:
\begin{equation}\label{eq.quad}
\left\langle \bv_{r+1},\blambda_{r+1}-\blambda_{r}\right\rangle=\frac{1}{2}\left(\|\blambda_{r+1}-\blambda_r\|^2+\|\bv_{r+1}\|^2-\|\blambda_r-\blambda_{r-1}\|^2\right)
\end{equation}
where
\begin{equation}
\bv_{r+1}\bydef\blambda_{r+1}-\blambda_r-(\blambda_r-\blambda_{r-1}).
\end{equation}

\item  Young's inequality with parameter $\theta>0$:
\begin{equation}
\langle \bx,\by\rangle\le\frac{1}{2\theta}\|\bx\|^2 + \frac{\theta}{2}\|\by\|^2,\quad \forall\bx,\by.
\end{equation}

\end{enumerate}

\subsection{Relation Among the Lemmas in the Proof}

During the theorem proving process, we will first show the descent lemma and then quantify the maximum ascent after the dual update \leref{le.des} . Next, from the dual update, we find the recursion of successive difference between two iterates \leref{le.recur} so that we can construct a potential function that shows the possible progress after one round of updating primal and dual variables \leref{le.cpo}. Giving the upper bound of the dual variable is a preliminary step shown in \leref{le.bdl} for functional constraints satisfaction. The key step is shown in \leref{le.sumd} that sharpens the bound of measuring the size of a sum of dual variables by using the regularity condition \leref{le.rege}, which plays the important role of characterizing the convergence of GDPA. At the same time, the optimality criterion is derived in terms of the successive difference of variables and the size of dual variables \leref{le.opt} so that we can use the obtained potential function to evaluate the descent achieved by GDPA. Combing these results together leads to the main theorem (\thref{th.main1}) of quantifying the convergence rate of GDPA to KKT points.

\section{On the Descent and Ascent of Potential Function}
In this section, we will provide detailed proofs of the convergence rate of GDPA. First, we give the following descent lemma that quantifies the decrease of the objective value after performing one round of GDPA update.
\vspace{-10pt}
\subsection{Descent Lemma}
\begin{mdframed}[backgroundcolor=gray!10,topline=false,
	rightline=false,
	leftline=false,
	bottomline=false]
\begin{lemma}\label{le.des}
Suppose that \asref{ass.lif}--\asref{ass.bfg} hold. If the iterates $\{\bx_r,\blambda_r,\forall r\}$ are generated by GDPA and the step-sizes $\alpha_r$ and $\beta_r$ satisfy
\begin{center}
\fbox{
\begin{minipage}{0.68\columnwidth}
\begin{equation}\label{eq.conal}
\frac{1}{\alpha_r}\ge L_f+(1-\tau)\|\blambda_r\|L_J+\beta_r U_J L_g,
\end{equation}
\end{minipage}
}
\end{center}
then, we have
\begin{align}\notag
&F_{\beta_{r+1}}(\bx_{r+1}, \blambda_{r+1})-F_{\beta_r}(\bx_r, \blambda_r)
\\\notag
\le & \frac{1-\tau}{2\beta_{r-1}}\|\blambda_r-\blambda_{r-1}\|^2-\left(\frac{1}{2\alpha_r}-\frac{3\beta_{r-1} L^2_g}{2}\right)\|\bx_{r+1}-\bx_r\|^2+\beta_{r-1}L^2_g\|\bx_{r}-\bx_{r-1}\|^2
\\\notag
&+\frac{1-\tau}{\beta_{r-1}}\|\blambda_{r+1}-\blambda_r\|^2
+(\beta_{r+1}-\beta_r)\|g_1(\bx_{r+1})\|^2
+\left(\frac{1}{\beta_r}-\frac{1}{\beta_{r+1}}\right)(1-\tau)^2\|\blambda_{r+1}\|^2
\\
&+\frac{(1-\tau)\gamma_r}{2}\|\blambda_{r+1}\|^2-\frac{(1-\tau)\gamma_{r-1}}{2}\|\blambda_r\|^2+(1-\tau)\left(\frac{\gamma_{r-1}-\gamma_r}{2}\right)\|\blambda_{r+1}\|^2.\label{eq.descentlemma}
\end{align}
\end{lemma}
\end{mdframed}
\begin{proof}
The proof mainly include two parts: 1) quantify the decrease of $F_{\beta_r}$ when primal variable $\bx$ is updated; 2) measure the ascent of $F_{\beta_r}$ after dual variable $\blambda$ is updated and $\beta_r$ is changed to $\beta_{r+1}$.

\noindent{\bf Part 1 Update of $\bx$ [$(\bx_{r},\blambda_r)\to(\bx_{r+1},\blambda_{r})$]}.

According to the gradient Lipschitz continuity shown in \asref{ass.lif}--\asref{ass.bfg} and \eqref{eq.as}, we have
\begin{equation}
F_{\beta_r}(\bx_{r+1},\blambda_r)\le F_{\beta_r}(\bx_r,\blambda_r)+\langle\nabla F_{\beta_r}(\bx_r,\blambda_r),\bx_{r+1}-\bx_r\rangle
+\frac{L_f+(1-\tau)\|\blambda_r\|L_J+\beta_r U_JL_g}{2}\|\bx_{r+1}-\bx_r\|^2\label{eq.des1}
\end{equation}
where we use the fact that the gradient Lipschitz constant is $L_f$ when $\left[g_i(\bx_r)+\frac{(1-\tau)[\blambda_r]_i}{\beta_r}\right]_+\le0, i\in\mathcal{S}_r$ and $L_f+(1-\tau)\|\blambda_r\|L_J+\beta_rU_JL_g$ when $\left[g_i(\bx_r)+\frac{(1-\tau)[\blambda_r]_i}{\beta_r}\right]_+>0, i\in\overbar{\mathcal{S}}_r$.

From the optimality condition of subproblem-$\bx$, we have
\begin{equation}
\langle\bx_{r+1}-(\bx_r-\alpha_r\nabla F_{\beta_r}(\bx_r,\blambda_r)), \bx_r-\bx_{r+1}\rangle\ge 0.\label{eq.optx}
\end{equation}
Therefore, substituting \eqref{eq.optx} into \eqref{eq.des1} gives
\begin{equation}\label{eq.descentx}
F_{\beta_r}(\bx_{r+1},\blambda_r)-F_{\beta_r}(\bx_r,\blambda_r)\le -\left(\frac{1}{\alpha_r}-\frac{L_f+(1-\tau)\|\blambda_r\|L_J+\beta_r U_J L_g}{2}\right)\|\bx_{r+1}-\bx_r\|^2.
\end{equation}

It can be seen that when
\begin{equation}
\frac{1}{\alpha_r}\ge L_f+(1-\tau)\|\blambda_r\|L_J+\beta_r U_J L_g,
\end{equation}
there is at least a decrease of $F_{\beta_r}$ in terms of $1/(2\alpha_r)\|\bx_{r+1}-\bx_r\|^2$ after the update of variable $\bx$.

\noindent{\bf Part 2 Update of $\blambda$ [$(\bx_{r+1},\blambda_r)\to(\bx_{r+1},\blambda_{r+1})$]}. There are two sub-cases in this part when $\blambda$ is udpated. To be more specific, we split partition $\sum_{i\in\mathcal{S}_r}(1-\tau)\blambda_i g_i(\bx) - \sum_{j\in\overbar{\mathcal{S}}_r}\frac{(1-\tau)^2\blambda^2_j}{2\beta_r}$ in function $F_{\beta_r} $ \eqref{eq.fbo} as the following two parts, i.e.,
\begin{equation}
h_1(\bx,\blambda_r)\bydef \sum_{i\in\mathcal{S}_r}[\blambda_r]_i g_i(\bx),
\quad
h_2(\blambda_r)\bydef  - \sum_{i\in\overbar{\mathcal{S}}_r}\frac{(1-\tau)[\blambda_r]^2_i}{2\beta_r}. \label{defh1h2}
\end{equation}

Then, from \eqref{eq.fbo}, we can have
\begin{equation}
F_{\beta_r}(\bx,\blambda_r)=f(\bx)+\sum_{i\in\mathcal{S}_r}\frac{\beta_r}{2}g^2_i(\bx) +(1-\tau) (\underbrace{h_1(\bx,\blambda_r)+h_2(\blambda_r)}_{\bydef h(\bx,\blambda_r)}).
\end{equation}

Let $h'(\bx,\blambda)\bydef h(\bx,\blambda)-\mathbbm{1}(\blambda)$, where $\mathbbm{1}(\blambda)$ denotes the indicator function. Let $\xi$ denote the subgradient of $\mathbbm{1}(\blambda)$. Since function $h'(\bx,\blambda)$ is concave with respect to $\blambda$, we have
\begin{align}
\notag
&h'(\bx_{r+1},\blambda_{r+1})-h'(\bx_{r+1},\blambda_r)
\\
\le & \langle\nabla_{\blambda} h(\bx_{r+1},\blambda_r),\blambda_{r+1}-\blambda_r\rangle-\langle\xi_r,\blambda_{r+1}-\blambda_r\rangle
\\\notag
=&\langle \nabla_{\blambda} h(\bx_{r+1},\blambda_r),\blambda_{r+1}-\blambda_r)\rangle
-\langle\xi_{r+1},\blambda_{r+1}-\blambda_r\rangle-\langle\xi_r-\xi_{r+1},\blambda_{r+1}-\blambda_r\rangle
\\\notag
\mathop{=}\limits^{(a)}&\frac{1-\tau}{\beta_r}\|\blambda_{r+1}-\blambda_r\|^2
+\gamma_r\langle \blambda_{r+1}, \blambda_{r+1}-\blambda_r\rangle+\langle\xi_{r+1}-\xi_r, \blambda_{r+1}-\blambda_r\rangle
\\\notag
\mathop{=}\limits^{(b)}&\frac{1-\tau}{\beta_r}\|\blambda_{r+1}-\blambda_r\|^2+\gamma_{r-1}\langle \blambda_{r},\blambda_{r+1}-\blambda_r\rangle
+\langle  \nabla_{\blambda} h(\bx_{r+1},\blambda_r) -\nabla_{\blambda} h(\bx_{r},\blambda_{r-1}), \blambda_{r+1}-\blambda_r\rangle
\\\notag
&-\frac{1-\tau}{\beta_{r-1}}\langle \bv_{r+1}, \blambda_{r+1}-\blambda_r\rangle+\left(\frac{1-\tau}{\beta_{r-1}}-\frac{1-\tau}{\beta_{r}}\right)\|\blambda_{r+1}-\blambda_r\|^2
\\\notag
\mathop{\le}\limits^{(c)}&\frac{1-\tau}{\beta_r}\|\blambda_{r+1}-\blambda_{r}\|^2+\!\!\!\!\!\!\sum_{i\in\mathcal{S}_{r+1}, i\in\mathcal{S}_r}\!\!\!  \frac{\beta_{r-1}}{2(1-\tau)} \|g_i(\bx_{r+1})-g_i(\bx_r)\|^2+\frac{1-\tau}{2\beta_{r-1}}\|[\blambda_{r+1}]_i-[\blambda_r]_i\|^2
\\\notag
& +\sum_{ i\in\mathcal{S}_{r+1}, i\in\overbar{\mathcal{S}}_r} \frac{\beta_{r-1}}{1-\tau}\|g_i(\bx_{r+1})-g_i(\bx_{r})\|^2+ \frac{\beta_{r-1}}{1-\tau}\|g_i(\bx_{r})-g_i(\bx_{r-1})\|^2+\frac{1-\tau}{2\beta_{r-1}}\|[\blambda_{r+1}]_i-[\blambda_r]_i\|^2
\\
&+\left(\frac{1-\tau}{\beta_{r-1}}-\frac{1-\tau}{\beta_{r}}\right)\|\blambda_{r+1}-\blambda_r\|^2+\gamma_{r-1}\langle \blambda_{r},\blambda_{r+1}-\blambda_r\rangle
-\frac{1-\tau}{\beta_{r-1}}\langle \bv_{r+1}, \blambda_{r+1}-\blambda_r\rangle\label{eq.lipgh}
\end{align}
where in $(a)$ we use the optimality condition of $\blambda$-problem, i.e.,
\begin{equation}
\xi_{r+1}- \nabla_{\blambda} h(\bx_{r+1},\blambda_r)+\frac{1-\tau}{\beta_r}(\blambda_{r+1}-\blambda_r)+\gamma_r\blambda_{r+1}=0,
\end{equation}
and in $(b)$ we substitute
\begin{multline}
\langle\xi_{r+1}-\xi_{r}, \blambda_{r+1}-\blambda_r\rangle=\langle \nabla_{\blambda} h(\bx_{r+1},\blambda_r) -\nabla_{\blambda} h(\bx_{r},\blambda_{r-1}), \blambda_{r+1}-\blambda_r\rangle+\left(\frac{1-\tau}{\beta_{r-1}}-\frac{1-\tau}{\beta_{r}}\right)\|\blambda_{r+1}-\blambda_r\|^2
\\
-\frac{1-\tau}{\beta_{r-1}}\langle \underbrace{(\blambda_{r+1}-\blambda_r)-(\blambda_r-\blambda_{r-1})}_{\bydef \bv_{r+1}}, \blambda_{r+1}-\blambda_r\rangle
-\langle\gamma_r \blambda_{r+1}-\gamma_{r-1} \blambda_r, \blambda_{r+1}-\blambda_r\rangle,
\end{multline}
and $(c)$ is true because
\begin{subequations}
\begin{align}\notag
&\quad\;\langle \nabla_{\blambda} h(\bx_{r+1},\blambda_r) -\nabla_{\blambda} h(\bx_{r},\blambda_{r-1}), \blambda_{r+1}-\blambda_r\rangle
\\
&= \sum_{i\in\mathcal{S}_{r}, i\in\mathcal{S}_{r-1}} \left(g_i(\bx_{r+1})- g_i(\bx^r)\right) ([\blambda_{r+1}]_i-[\blambda_r]_i)
\\
&\quad\quad  + \sum_{i\in\mathcal{S}_{r}, i\in\overbar{\mathcal{S}}_{r-1}} \left(g_i(\bx^{r+1})- -\frac{1-\tau}{\beta_{r-1}}[\blambda_{r-1}]_i\right) ([\blambda_{r+1}]_i-[\blambda_r]_i) \label{eq.err2}
\\
&\quad\quad  +\sum_{i\in\overbar{\mathcal{S}}_{r},i\in\mathcal{S}_{r-1} } \left(-\frac{1-\tau}{\beta_r}[\blambda_r]_i- g_i(\bx^r)\right) ([\blambda_{r+1}]_i-[\blambda_r]_i)\label{eq.err3}
\\
&\quad\quad  + \sum_{i\in\overbar{\mathcal{S}}_{r},i\in\overbar{\mathcal{S}}_{r-1}}\left (-\frac{1-\tau}{\beta_r}[\blambda_r]_i- -\frac{1-\tau}{\beta_{r-1}}[\blambda_{r-1}]_i\right)([\blambda_{r+1}]_i-[\blambda_r]_i)\label{eq.err4}
\end{align}
\end{subequations}
where we decompose $\nabla_{\blambda} h(\bx_{r+1},\blambda_r) -\nabla_{\blambda} h(\bx_{r},\blambda_{r-1})$ by definition  \eqref{defh1h2} and each of the term can be bounded as follows:  1) applying the gradient Lipschitz continuity and Young's inequality, we have
\begin{align}\notag
&\quad\sum_{i\in\mathcal{S}_{r}, i\in\mathcal{S}_{r-1}} (g_i(\bx^{r+1})- g_i(\bx^r)) ([\blambda_{r+1}]_i-[\blambda_r]_i)
\\
&\le\!\!\!\sum_{i\in\mathcal{S}_{r}, i\in\mathcal{S}_{r-1}}\!\!\!\frac{\beta_{r-1} }{2(1-\tau)}\|g_i(\bx^{r+1})-g_i(\bx^r)\|^2+\frac{1-\tau}{2\beta_{r-1}}\|[\blambda_{r+1}]_i-[\blambda_r]_i\|^2;\label{eq.t1}
\end{align}
2) term \eqref{eq.err2} can be decomposed as  
\begin{align}\notag
&\quad \left(g_i(\bx_{r+1})- -\frac{1-\tau}{\beta_{r-1}}[\blambda_{r-1}]_i\right) ([\blambda_{r+1}]_i-[\blambda_r]_i)
\\\notag
&\mathop{=}\limits^{(c.1)} (g_i(\bx_{r+1})-g_i(\bx_r)([\blambda_{r+1}]_i-[\blambda_r]_i)+
(g_i(\bx_{r})- 0)([\blambda_{r+1}]_i - [\blambda_{r}]_i)
\\
&\quad +\left(-\frac{1-\tau}{\beta_{r-1}}[\blambda_{r}]_i- -\frac{1-\tau}{\beta_{r-1}}[\blambda_{r-1}]_i\right) ([\blambda_{r+1}]_i-[\blambda_r]_i)
 \\\notag
 & \mathop{\le}\limits^{(c.2)}  \frac{\beta_{r-1}}{1-\tau} \|g_i(\bx_{r+1})-g_i(\bx_{r})\|^2+\frac{1-\tau}{4\beta_{r-1}}\|[\blambda_{r+1}]_i-[\blambda_r]_i\|^2+ \frac{\beta_{r-1}}{1-\tau} \|g_i(\bx_{r})-g_i(\bx_{r-1})\|^2
 \\
&\quad+\frac{1-\tau}{4\beta_{r-1}}\|[\blambda_{r+1}]_i-[\blambda_r]_i\|^2
 +\frac{1-\tau}{\beta_{r-1}}\left(-[\blambda_{r}]_i- (- [\blambda_{r-1}]_i)\right) ([\blambda_{r+1}]_i-[\blambda_r]_i)
\\
 & \mathop{\le}\limits^{(c.3)} \notag \frac{\beta_{r-1}}{1-\tau}\|g_i(\bx_{r+1})-g_i(\bx_{r})\|^2+ \frac{\beta_{r-1}}{1-\tau}\|g_i(\bx_{r})-g_i(\bx_{r-1})\|^2+\frac{1-\tau}{2\beta_{r-1}}\|[\blambda_{r+1}]_i-[\blambda_r]_i\|^2
\\
&\quad  +\frac{1-\tau}{2\beta_{r-1}}\|[\bv_{r+1}]_i\|^2-\frac{1-\tau}{2\beta_{r-1}}\|[\blambda_r]_i-[\blambda_{r-1}]_i\|^2\label{eq.t2}
\end{align}
where $(c.1)$ is true because $i\in\overbar{\mathcal{S}}_{r-1}$ implies $g_i(\bx_{r-1})\le0$ and $[\blambda_{r}]_i=0$ according to \eqref{eq.dualupated}, $(c.2)$ holds due to the fact that $i\in\mathcal{S}_{r}$ implies $g_i(\bx_{r})>-(1-\tau)[\blambda_r]_i/\beta_r$, leading to  $g_i(\bx_{r})>0$, and in $(c.3)$ we apply
\begin{align}\notag
& \left\langle -([\blambda_{r}]_i-[\blambda_{r-1}]_i), [\bv_{r+1}]_i + [\blambda_{r}]_i-[\blambda_{r-1}]_i\right\rangle
\\
\le& \left \langle -([\blambda_{r}]_i-[\blambda_{r-1}]_i), [\bv_{r+1}]_i\right\rangle- \|[\blambda_r]_i-[\blambda_{r-1}]_i\|^2
\\
\le&  \frac{1}{2}\|[\blambda_r]_i-[\blambda_{r-1}]_i\|^2+\frac{1}{2}\|[\bv_{r+1}]_i\|^2- \|[\blambda_r]_i-[\blambda_{r-1}]_i\|^2
\\
\le&  \frac{1}{2}\|[\bv_{r+1}]_i\|^2-\frac{1}{2}\left\|[\blambda_r]_i-[\blambda_{r-1}]_i\right\|^2;
\end{align}
3) regarding term \eqref{eq.err3},  we can have the following two facts: first, $i\in\overbar{\mathcal{S}}_r$ implies $g_i(\bx_{r})+\frac{1-\tau}{\beta_r}[\blambda_r]_i\le0$ according to \eqref{eq.defs}; second, $i\in\overbar{\mathcal{S}}_{r}$ also gives $\blambda_{r+1}= 0$ according to \eqref{eq.dualupated}. Combining both of them term \eqref{eq.err3} can be further upper bounded by
\begin{equation}
\sum_{i\in\overbar{\mathcal{S}}_{r},i\in\mathcal{S}_{r-1} } \left(-\frac{1-\tau}{\beta_r}[\blambda_r]_i - g_i(\bx^r)\right) ([\blambda_{r+1}]_i-[\blambda_r]_i)\le0; \label{eq.t5}
\end{equation}
 and  4) the fourth term is not greater than $0$ due to $[\blambda_{r+1}]_i=0$ and $[\blambda_{r}]_i=0$ in this case.

Note that
\begin{align}\notag
 &\sum_{i\in\mathcal{S}_{r+1}, i\in\mathcal{S}_r \cup i\in\mathcal{S}_{r+1}, i\in\overbar{\mathcal{S}}_r}  \!\!\!\!\frac{\beta_{r-1}}{2(1-\tau)}  \|g_i(\bx_{r+1})-g_i(\bx_r)\|^2+\frac{\beta_{r-1}}{1-\tau} \|g_i(\bx_{r+1})-g_i(\bx_r)\|^2+\frac{1-\tau}{2\beta_{r-1}}\|[\blambda_{r+1}]_i-[\blambda_r]_i\|^2
 \\\notag
 &\quad\quad \quad \quad \quad  +\sum_{ i\in\mathcal{S}_{r+1}, i\in\overbar{\mathcal{S}}_r} \frac{\beta_{r-1}}{1-\tau}\|g_i(\bx_{r})-g_i(\bx_{r-1})\|^2
 \\
 &\le\frac{3\beta_{r-1} L^2_g}{2(1-\tau)} \|\bx_{r+1}-\bx_r\|^2+\frac{\beta_{r-1} L^2_g}{1-\tau} \|\bx_{r}-\bx_{r-1}\|^2+\frac{1-\tau}{2\beta_{r-1}}\|\blambda_{r+1}-\blambda_r\|^2.
\end{align}
Therefore, combing the results \eqref{eq.t1}, \eqref{eq.t2}, \eqref{eq.t5}, we can further bound $h'(\bx_{r+1},\blambda_{r+1})-h'(\bx_{r+1},\blambda_r)$ as follows.
\begin{align}\notag
&\quad h'(\bx_{r+1},\blambda_{r+1})-h'(\bx_{r+1},\blambda_r)
\\\notag
&\mathop{\le}\limits^{(a)}\frac{ 1-\tau }{ \beta_{r}}\|\blambda_{r+1}-\blambda_{r}\|^2+  \frac{3\beta_{r-1} L^2_g}{2(1-\tau)}\|\bx_{r+1}-\bx_r\|^2+  \frac{\beta_{r-1} L^2_g}{1-\tau}\|\bx_{r}-\bx_{r-1}\|^2+\frac{1-\tau}{2\beta_{r-1}}\|\blambda_{r+1}-\blambda_r\|^2
\\\notag
&\quad +\sum_{ i\in\mathcal{S}_{r+1}, i\in\overbar{\mathcal{S}}_r} \frac{1-\tau}{2\beta_{r-1}}\|[\bv_{r+1}]_i\|^2-\frac{1-\tau}{2\beta_{r-1}}\|[\blambda_r]_i-[\blambda_{r-1}]_i\|^2
\\\notag
&\quad+\left(\frac{1-\tau}{\beta_{r-1}}-\frac{1-\tau}{\beta_{r}}\right)\|\blambda_{r+1}-\blambda_r\|^2+\gamma_{r-1}\langle \blambda_{r},\blambda_{r+1}-\blambda_r\rangle
\\
&\quad-\frac{1-\tau}{2\beta_{r-1}}\|\blambda_{r+1}-\blambda_r\|^2-\frac{1-\tau}{2\beta_{r-1}}\|\bv_{r+1}\|^2+\frac{1-\tau}{2\beta_{r-1}}\|\blambda_{r}-\blambda_{r-1}\|^2
\\\notag
&=\frac{ 1-\tau }{ \beta_{r-1}}\|\blambda_{r+1}-\blambda_{r}\|^2+  \frac{3\beta_{r-1} L^2_g}{2(1-\tau)}\|\bx_{r+1}-\bx_r\|^2+  \frac{\beta_{r-1} L^2_g}{1-\tau}\|\bx_{r}-\bx_{r-1}\|^2
\\
&\quad+\frac{1-\tau}{2\beta_{r-1}}\|\blambda_{r}-\blambda_{r-1}\|^2+\gamma_{r-1}\langle \blambda_{r},\blambda_{r+1}-\blambda_r\rangle
\\\notag
&\mathop{\le}\limits^{(b)} \frac{1-\tau}{2\beta_{r-1}}\|\blambda_r-\blambda_{r-1}\|^2+  \frac{3\beta_{r-1} L^2_g}{2(1-\tau)}\|\bx_{r+1}-\bx_r\|^2+  \frac{\beta_{r-1} L^2_g}{1-\tau}\|\bx_{r}-\bx_{r-1}\|^2-\left(\frac{\gamma_{r-1}}{2}-\frac{1-\tau}{\beta_{r-1}}\right)\|\blambda_{r+1}-\blambda_r\|^2
\\
&\quad+\frac{\gamma_r}{2}\|\blambda_{r+1}\|^2-\frac{\gamma_{r-1}}{2}\|\blambda_r\|^2+\frac{\gamma_{r-1}-\gamma_r}{2}\|\blambda_{r+1}\|^2\label{eq.descenty}
\end{align}
where $(a)$ holds by applying the Lipschitz continuity and quadrilateral identity \eqref{eq.quad} and in $(b)$ we use the following fact
\begin{align}
\gamma_{r-1}\langle \blambda_r,\blambda_{r+1}-\blambda_r\rangle
&=\frac{\gamma_{r-1}}{2}\left(\|\blambda_{r+1}\|^2-\|\blambda_r\|^2-\|\blambda_{r+1}-\blambda_r\|^2\right)
\\
&=\frac{\gamma_r}{2}\|\blambda_{r+1}\|^2-\frac{\gamma_{r-1}}{2}\|\blambda_r\|^2-\frac{\gamma_{r-1}}{2}\|\blambda_{r+1}-\blambda_r\|^2+\left(\frac{\gamma_{r-1}-\gamma_r}{2}\right)\|\blambda_{r+1}\|^2.\label{eq.gl}
\end{align}

Note that from \eqref{eq.fbo} we have
\begin{equation}\label{eq.bec}
F_{\beta_{r+1}}(\bx_{r+1},\blambda_{r+1})\le F_{\beta_r}(\bx_{r+1},\blambda_{r+1})+(\beta_{r+1}-\beta_r)\|g(\bx_{r+1})\|^2
+\left(\frac{1}{\beta_r}-\frac{1}{\beta_{r+1}}\right)(1-\tau)^2\|\blambda_{r+1}\|^2
\end{equation}
where we use the following fact
\begin{equation}
\sum_{i\in\mathcal{S}_r}\frac{\beta_{r+1}}{2}g^2_i(\bx_{r+1})-\frac{\beta_r}{2}g^2_i(\bx_{r+1})=\frac{\beta_{r+1}-\beta_r}{2}\|g_1(\bx_{r+1})\|^2.
\end{equation}
Combining \eqref{eq.descenty}and  \eqref{eq.bec}, we can obtain the changes of $F_{\beta_r}$ after one round update of $\bx$ and $\blambda$, i.e.,
\begin{align}\notag
&F_{\beta_{r+1}}(\bx_{r+1}, \blambda_{r+1})-F_{\beta_r}(\bx_r, \blambda_r)
\\\notag
\le&\frac{(1-\tau)^2}{2\beta_{r-1}}\|\blambda_r-\blambda_{r-1}\|^2-\left(\frac{1}{2\alpha_r}-\frac{3\beta_{r-1} L^2_g}{2}\right)\|\bx_{r+1}-\bx_r\|^2 + \beta_{r-1}L^2_g\|\bx_{r}-\bx_{r-1}\|^2
\\\notag
&+\frac{(1-\tau)^2}{\beta_{r-1}}\|\blambda_{r+1}-\blambda_r\|^2
+(\beta_{r+1}-\beta_r)\|g_1(\bx_{r+1})\|^2
+\left(\frac{1}{\beta_r}-\frac{1}{\beta_{r+1}}\right)(1-\tau)^2\|\blambda_{r+1}\|^2
\\
&+\frac{(1-\tau)\gamma_r}{2}\|\blambda_{r+1}\|^2-\frac{(1-\tau)\gamma_{r-1}}{2}\|\blambda_r\|^2+(1-\tau)\left(\frac{\gamma_{r-1}-\gamma_r}{2}\right)\|\blambda_{r+1}\|^2,
\end{align}
which gives the desired result as $0<\tau<1$.
\end{proof}

\subsection{Recursion of the Size of the Difference between Two Successive Dual Variables}
\begin{mdframed}[backgroundcolor=gray!10,
topline=false,rightline=false,leftline=false,bottomline=false]
\begin{lemma}\label{le.recur}
Suppose that \asref{ass.lif}--\asref{ass.lig} hold. If the iterates $\{\bx_r,\blambda_r,\forall r\}$ are generated by GDPA, then we have
\begin{align}\notag
&\frac{2\vartheta(1-\tau)}{\beta_{r-1}\tau}\|\blambda_{r+1}-\blambda_{r}\|^2-\frac{2\vartheta}{\beta_{r}}\left(\frac{\gamma_{r-1}}{\gamma_r}-1\right)\|\blambda_{r+1}\|^2
\\\notag
\le&\frac{2\vartheta(1-\tau)}{\beta_{r-2}\tau}\|\blambda_r-\blambda_{r-1}\|^2-\frac{2\vartheta}{\beta_{r-1}}\left(\frac{\gamma_{r-2}}{\gamma_{r-1}}-1\right)\|\blambda_{r}\|^2+\frac{2\vartheta}{\beta_r}\left(\frac{\gamma_{r-2}}{\gamma_{r-1}}-\frac{\gamma_{r-1}}{\gamma_r}\right)\|\blambda_r\|^2
-\frac{2\vartheta}{\beta_r}\|\blambda_{r+1}-\blambda_r\|^2
\\\notag
&+\frac{6\vartheta L^2_{g}}{\gamma_{r-1}(1-\tau)}\|\bx_{r+1}-\bx_r\|^2+\frac{4\vartheta L^2_{g}}{\gamma_{r-1}(1-\tau)}\|\bx_{r}-\bx_{r-1}\|^2
\\
&+\frac{6\vartheta}{\tau}\left(\frac{2-\tau}{\beta_{r-1}}-\frac{2-\tau}{\beta_{r}}\right)\|\blambda_{r+1}-\blambda_r\|^2 +\left(\frac{2\vartheta}{\beta_{r-1}}-\frac{2\vartheta}{\beta_r}\right)\left(\frac{\gamma_{r-2}}{\gamma_{r-1}}-1\right)\|\blambda_{r}\|^2.\label{eq.diffiter}
\end{align}
\end{lemma}
\vspace{-5pt}
\end{mdframed}

\begin{proof}
 From the optimality condition of $\blambda$-problem at the $r+1$th iteration, we have
\begin{equation}\label{eq.yopt1}
-\langle \nabla_{\blambda} h(\bx_{r+1},\blambda_r) -\frac{1-\tau}{\beta_r}(\blambda_{r+1}-\blambda_r)-\gamma_r\blambda_{r+1}, \blambda_{r+1}-\blambda\rangle\le 0, \quad\forall \blambda\ge 0.
\end{equation}
Similarly, from the optimality condition of $\blambda$-problem at the $r$th iteration, we have
\begin{equation}\label{eq.yopt2}
-\langle \nabla_{\blambda} h(\bx_{r},\blambda_{r-1})-\frac{1-\tau}{\beta_{r-1}}(\blambda_{r}-\blambda_{r-1})-\gamma_{r-1}\blambda_{r}, \blambda-\blambda_r\rangle\ge 0, \quad \forall \blambda\ge 0.
\end{equation}

Plugging in $\blambda=\blambda_r$ in \eqref{eq.yopt1}, $\blambda=\blambda_{r+1}$ in \eqref{eq.yopt2} and combining them together, we can get
\begin{align}\notag
&\frac{1-\tau}{\beta_{r-1}}\langle \bv_{r+1}, \blambda_{r+1}-\blambda_r\rangle + \langle\gamma_r\blambda_{r+1}-\gamma_{r-1}\blambda_r, \blambda_{r+1}-\blambda_r\rangle
\\
\le &  \langle  \nabla_{\blambda} h(\bx_{r+1},\blambda_r)  - \nabla_{\blambda} h(\bx_{r},\blambda_{r-1}),\blambda_{r+1}-\blambda_r\rangle
+\left(\frac{1}{\beta_{r-1}}-\frac{1}{\beta_{r}}\right)\|\blambda_{r+1}-\blambda_r\|^2. \label{eq.recur}
\end{align}

In the following, we will use this inequality to analyze the recurrence of the size of the difference between two consecutive iterates. First, we have
\begin{align}
\notag
&\langle\gamma_r\blambda_{r+1}-\gamma_{r-1}\blambda_r, \blambda_{r+1}-\blambda_r\rangle
\\
=&\langle\gamma_r\blambda_{r+1}-\gamma_r\blambda_r+\gamma_r\blambda_r-\gamma_{r-1}\blambda_r, \blambda_{r+1}-\blambda_r\rangle
\\\notag
=&\gamma_r\|\blambda_{r+1}-\blambda_r\|^2+(\gamma_r-\gamma_{r-1})\langle \blambda_r, \blambda_{r+1}-\blambda_r\rangle
\\\notag
=&\gamma_r\|\blambda_{r+1}-\blambda_r\|^2+\frac{\gamma_r-\gamma_{r-1}}{2}\left(\|\blambda_{r+1}\|^2-\|\blambda_r\|^2-\|\blambda_{r+1}-\blambda_r\|^2\right)
\\
=&\frac{\gamma_{r}+\gamma_{r-1}}{2}\|\blambda_{r+1}-\blambda_r\|^2-\frac{\gamma_{r-1}-\gamma_r}{2}\left(\|\blambda_{r+1}\|^2-\|\blambda_r\|^2\right),\label{eq.gammar}
\end{align}
and quadrilateral identity
\begin{equation}\label{eq.eqrel}
\langle \bv_{r+1}, \blambda_{r+1}-\blambda_r\rangle=\frac{1}{2}\left(\|\blambda_{r+1}-\blambda_r\|^2+\|\bv_{r+1}\|^2-\|\blambda_r-\blambda_{r-1}\|^2\right).
\end{equation}

Next, substituting \eqref{eq.gammar}  and \eqref{eq.eqrel} into \eqref{eq.recur}, we have
\begin{align}\notag
&\quad\;\frac{1-\tau}{2\beta_{r-1}}\|\blambda_{r+1}-\blambda_r\|^2-\frac{\gamma_{r-1}-\gamma_r}{2}\|\blambda_{r+1}\|^2
\\\notag
&\le  \frac{1-\tau}{2\beta_{r-1}}\|\blambda_r-\blambda_{r-1}\|^2-\frac{1-\tau}{2\beta_{r-1}}\|\bv_{r+1}\|^2-\frac{\gamma_{r-1}-\gamma_r}{2}\|\blambda_r\|^2-\frac{\gamma_r+\gamma_{r-1}}{2}\|\blambda_{r+1}-\blambda_r\|^2
\\
&\quad+\langle \nabla_{\blambda} h(\bx_{r+1},\blambda_r)  - \nabla_{\blambda} h(\bx_{r},\blambda_{r-1}), \blambda_{r+1}-\blambda_r\rangle+\left(\frac{1}{\beta_{r-1}}-\frac{1}{\beta_{r}}\right)\|\blambda_{r+1}-\blambda_r\|^2
\\\notag
&\mathop{\le}\limits^{(a)}\frac{1-\tau}{2\beta_{r-1}}\|\blambda_r-\blambda_{r-1}\|^2-\frac{1-\tau}{2\beta_{r-1}}\|\bv_{r+1}\|^2-\gamma_r\|\blambda_{r+1}-\blambda_r\|^2-\frac{\gamma_{r-1}-\gamma_r}{2}\|\blambda_r\|^2
\\
&\quad+\langle \nabla_{\blambda} h(\bx_{r+1},\blambda_r)  - \nabla_{\blambda} h(\bx_{r},\blambda_{r-1}), \blambda_{r+1}-\blambda_r\rangle+\left(\frac{1}{\beta_{r-1}}-\frac{1}{\beta_{r}}\right)\|\blambda_{r+1}-\blambda_r\|^2\label{eq.stre}
\\\notag
&\mathop{\le}\limits^{(b)}\frac{1-\tau}{2\beta_{r-1}}\|\blambda_r-\blambda_{r-1}\|^2-\gamma_r\|\blambda_{r+1}-\blambda_r\|^2-\frac{\gamma_{r-1}-\gamma_r}{2}\|\blambda_r\|^2+\frac{1-\tau}{2\beta_{r-1}}\|\blambda_{r+1}-\blambda_r\|^2
\\
&\quad +  \frac{3\beta_{r-1} L^2_g}{2(1-\tau)}\|\bx_{r+1}-\bx_r\|^2+  \frac{\beta_{r-1} L^2_g}{1-\tau}\|\bx_{r}-\bx_{r-1}\|^2 +\left(\frac{2-\tau}{\beta_{r-1}}-\frac{2-\tau}{\beta_{r}}\right)\|\blambda_{r+1}-\blambda_r\|^2
\\\notag
&\mathop{\le}\limits^{(c)}\frac{1-\tau}{2\beta_{r-1}}\|\blambda_r-\blambda_{r-1}\|^2-\frac{\gamma_r}{2}\|\blambda_{r+1}-\blambda_r\|^2-\frac{\gamma_{r-1}-\gamma_r}{2}\|\blambda_r\|^2
\\
&\quad +  \frac{3\beta_{r-1} L^2_g}{2(1-\tau)}\|\bx_{r+1}-\bx_r\|^2+  \frac{\beta_{r-1} L^2_g}{1-\tau}\|\bx_{r}-\bx_{r-1}\|^2 +\frac{3}{2}\left(\frac{2-\tau}{\beta_{r-1}}-\frac{2-\tau}{\beta_{r}}\right)\|\blambda_{r+1}-\blambda_r\|^2
\end{align}
where $(a)$ is true because $0<\gamma_r<\gamma_{r-1}$, in $(b)$ we use \eqref{eq.lipgh}, and in $(c)$ we choose $\tau>1/2$ such that $\gamma_r/2\ge (1-\tau)/(2\beta_{r-1})$.

Multiplying by 4$\vartheta$ and dividing by $\tau$ on the both sides of the above equation, we can get
\begin{align}\notag
&\quad\;\frac{2\vartheta(1-\tau)}{\beta_{r-1}\tau}\|\blambda_{r+1}-\blambda_{r}\|^2-\frac{2\vartheta}{\beta_r}\left(\frac{\gamma_{r-1}}{\gamma_r}-1\right)\|\blambda_{r+1}\|^2
\\\notag
&\le \frac{2\vartheta(1-\tau)}{\beta_{r-1}\tau}\|\blambda_r-\blambda_{r-1}\|^2-\frac{2\vartheta}{\beta_r}\left(\frac{\gamma_{r-1}}{\gamma_r}-1\right)\|\blambda_{r}\|^2
-\frac{2\vartheta}{\beta_r}\|\blambda_{r+1}-\blambda_r\|^2+\frac{6\vartheta L^2_{g}}{\gamma_{r-1}(1-\tau)}\|\bx_{r+1}-\bx_r\|^2
\\
& \quad +\frac{4\vartheta L^2_{g}}{\gamma_{r-1}(1-\tau)}\|\bx_{r}-\bx_{r-1}\|^2+\frac{6\vartheta}{\tau}\left(\frac{2-\tau}{\beta_{r-1}}-\frac{2-\tau}{\beta_{r}}\right)\|\blambda_{r+1}-\blambda_r\|^2\label{eq.cre}
\\\notag
&\le\frac{2\vartheta(1-\tau)}{\beta_{r-1}\tau}\|\blambda_r-\blambda_{r-1}\|^2-\frac{2\vartheta}{\beta_{r-1}}\left(\frac{\gamma_{r-2}}{\gamma_{r-1}}-1\right)\|\blambda_{r}\|^2+\frac{2\vartheta}{\beta_r}\left(\frac{\gamma_{r-2}}{\gamma_{r-1}}-\frac{\gamma_{r-1}}{\gamma_r}\right)\|\blambda_r\|^2
-\frac{2\vartheta}{\beta_r}\|\blambda_{r+1}-\blambda_r\|^2
\\
&\quad+\frac{6\vartheta L^2_{g}}{\gamma_{r-1}(1-\tau)}\|\bx_{r+1}-\bx_r\|^2+\frac{4\vartheta L^2_{g}}{\gamma_{r-1}(1-\tau)}\|\bx_{r}-\bx_{r-1}\|^2+\frac{6\vartheta}{\tau}\left(\frac{2-\tau}{\beta_{r-1}}-\frac{2-\tau}{\beta_{r}}\right)\|\blambda_{r+1}-\blambda_r\|^2
\\\notag
&\mathop{\le}\limits^{(a)} \frac{2\vartheta(1-\tau)}{\beta_{r-2}\tau}\|\blambda_r-\blambda_{r-1}\|^2-\frac{2\vartheta}{\beta_{r-1}}\left(\frac{\gamma_{r-2}}{\gamma_{r-1}}-1\right)\|\blambda_{r}\|^2+\frac{2\vartheta}{\beta_r}\left(\frac{\gamma_{r-2}}{\gamma_{r-1}}-\frac{\gamma_{r-1}}{\gamma_r}\right)\|\blambda_r\|^2
-\frac{2\vartheta}{\beta_r}\|\blambda_{r+1}-\blambda_r\|^2
\\\notag
&\quad+\frac{6\vartheta L^2_{g}}{\gamma_{r-1}(1-\tau)}\|\bx_{r+1}-\bx_r\|^2+\frac{4\vartheta L^2_{g}}{\gamma_{r-1}(1-\tau)}\|\bx_{r}-\bx_{r-1}\|^2
\\
&\quad+\frac{6\vartheta}{\tau}\left(\frac{2-\tau}{\beta_{r-1}}-\frac{2-\tau}{\beta_{r}}\right)\|\blambda_{r+1}-\blambda_r\|^2 +\left(\frac{2\vartheta}{\beta_{r-1}}-\frac{2\vartheta}{\beta_r}\right)\left(\frac{\gamma_{r-2}}{\gamma_{r-1}}-1\right)\|\blambda_{r}\|^2\label{eq.diffiter1}
\end{align}
where in $(a)$ we use $1/\beta_{r-1}<1/\beta_{r-2}$.
 \end{proof}
\subsection{Construction of a Potential Function}
\begin{mdframed}[backgroundcolor=gray!10,topline=false,
	rightline=false,
	leftline=false,
	bottomline=false]
\begin{lemma}\label{le.cpo}
Suppose that \asref{ass.lif}--\asref{ass.bfg} hold and the iterates $\{\bx_r,\blambda_r,\forall r\}$ are generated by GDPA. If the step-sizes satisfy
\begin{center}
\fbox{
\begin{minipage}{0.68\columnwidth}
\begin{equation}\label{eq.stc1}
\frac{1}{4\alpha_r}>\frac{5\beta_{r} L^2_g}{2}+\frac{2\vartheta L^2_{g}}{\gamma_r}\left(\frac{1}{\tau}+\frac{2}{1-\tau}\right),
\end{equation}
\end{minipage}
}
\end{center}
and constant $\vartheta>1$, then, we have
\begin{align} \notag
 \frac{1}{2\beta_{r}} &\|\blambda_{r+1}-\blambda_r\|^2 +\frac{1}{4\alpha_r}\|\bx_{r+1}-\bx_r\|^2
\le 
\mathcal{P}_r-\mathcal{P}_{r+1}+\frac{2\vartheta}{\beta_r}\left(\frac{\gamma_{r-2}}{\gamma_{r-1}}-\frac{\gamma_{r-1}}{\gamma_r}\right)\|\blambda_r\|^2
\\\notag
&+\left(\frac{2\vartheta}{\beta_{r-1}}-\frac{2\vartheta}{\beta_r}\right)\left(\frac{\gamma_{r-2}}{\gamma_{r-1}}-1\right)\|\blambda_{r}\|^2
+\frac{6\vartheta(3-2\tau)}{\tau}\left(\frac{1}{\beta_{r-1}}-\frac{1}{\beta_{r}}\right)\|\blambda_{r+1}-\blambda_r\|^2
\\
&+(\beta_{r+1}-\beta_r)\|g_1(\bx_{r+1})\|^2
+\left(\frac{1}{\beta_r}-\frac{1}{\beta_{r+1}}\right)(1-\tau)^2\|\blambda_{r+1}\|^2
+(1-\tau)\left(\frac{\gamma_{r-1}-\gamma_r}{2}\right)\|\blambda_{r+1}\|^2 \label{eq.recc}
\end{align}
where the potential function is defined as
\begin{multline}\label{eq.defp}
\mathcal{P}_{r}\bydef F_{\beta_r}(\bx_r,\blambda_r)+\left(\frac{1-\tau}{2\beta_{r-1}}+\frac{2\vartheta(1-\tau)}{\beta_{r-2}\tau}\right)\|\blambda_{r}-\blambda_{r-1}\|^2-\frac{(1-\tau)\gamma_{r-1}}{2}\|\blambda_r\|^2
\\
-\frac{2\vartheta}{\beta_{r-1}}\left(\frac{\gamma_{r-2}}{\gamma_{r-1}}-1\right)\|\blambda_r\|^2
+\beta_{r-1}L^2_g\|\bx_{r}-\bx_{r-1}\|^2+\frac{4\vartheta L^2_{g}}{\gamma_{r-1}(1-\tau)}\|\bx_{r}-\bx_{r-1}\|^2.
\end{multline}
\end{lemma}
\end{mdframed}

\begin{proof}
 Combining \eqref{eq.descentlemma} in \leref{le.des} and \eqref{eq.diffiter} in \leref{le.recur}, we have
\begin{align}\notag
&\quad\;F_{\beta_{r+1}}(\bx_{r+1},\blambda_{r+1})+\left(\frac{1-\tau}{2\beta_r}+\frac{2\vartheta(1-\tau)}{\beta_{r-1}\tau}\right)\|\blambda_{r+1}-\blambda_r\|^2-\frac{(1-\tau)\gamma_r}{2}\|\blambda_{r+1}\|^2-\frac{2\vartheta}{\beta_r}\left(\frac{\gamma_{r-1}}{\gamma_r}-1\right)\|\blambda_{r+1}\|^2
\\ \notag
&\quad\;+\beta_{r}L^2_g\|\bx_{r+1}-\bx_{r}\|^2+\frac{4\vartheta L^2_{g}}{\gamma_{r}(1-\tau)}\|\bx_{r+1}-\bx_{r}\|^2
\\\notag
&\le  F_{\beta_r}(\bx_r,\blambda_r)+\left(\frac{1-\tau}{2\beta_{r-1}}+\frac{2\vartheta(1-\tau)}{\beta_{r-2}\tau}\right)\|\blambda_{r}-\blambda_{r-1}\|^2-\frac{(1-\tau)\gamma_{r-1}}{2}\|\blambda_r\|^2-\frac{2\vartheta}{\beta_{r-1}}\left(\frac{\gamma_{r-2}}{\gamma_{r-1}}-1\right)\|\blambda_r\|^2
\\ \notag
&\quad +\beta_{r-1}L^2_g\|\bx_{r}-\bx_{r-1}\|^2+\frac{4\vartheta L^2_{g}}{\gamma_{r-1}(1-\tau)}\|\bx_{r}-\bx_{r-1}\|^2
\\\notag
&\quad-\left(\frac{2\vartheta}{\beta_r}-\left(\frac{1-\tau}{2\beta_{r}}+\frac{1-\tau}{\beta_{r-1}}+\frac{6\vartheta}{\tau}\left(\frac{2-\tau}{\beta_{r-1}}-\frac{2-\tau}{\beta_{r}}\right)\right)\right)\|\blambda_{r+1}-\blambda_r\|^2
\\\notag
&\quad-\left(\frac{1}{2\alpha_r}-\left(\frac{5\beta_{r} L^2_g}{2}+\frac{2\vartheta L^2_{g}}{\gamma_r\tau}+\frac{4\vartheta L^2_{g}}{\gamma_{r}(1-\tau)}\right)\right)\|\bx_{r+1}-\bx_r\|^2
\\\notag
&\quad+\frac{2\vartheta}{\beta_r}\left(\frac{\gamma_{r-2}}{\gamma_{r-1}}-\frac{\gamma_{r-1}}{\gamma_r}\right)\|\blambda_r\|^2
+\left(\frac{2\vartheta}{\beta_{r-1}}-\frac{2\vartheta}{\beta_r}\right)\left(\frac{\gamma_{r-2}}{\gamma_{r-1}}-1\right)\|\blambda_{r}\|^2
\\
&\quad+(\beta_{r+1}-\beta_r)\|g_1(\bx_{r+1})\|^2
+\left(\frac{1}{\beta_r}-\frac{1}{\beta_{r+1}}\right)(1-\tau)^2\|\blambda_{r+1}\|^2
+(1-\tau)\left(\frac{\gamma_{r-1}-\gamma_r}{2}\right)\|\blambda_{r+1}\|^2.\label{eq.whole}
\end{align}

It can be easily checked that when $\alpha_r$ satisfies \eqref{eq.stc1}, $\tau\ge1/2$ and $\vartheta>1$, we can further rewrite \eqref{eq.whole} as
\begin{align}\notag
\mathcal{P}_{r+1}& \le  \mathcal{P}_{r} - \frac{1}{4\alpha_r} \|\bx_{r+1}-\bx_r\|^2
-\frac{1}{2\beta_{r}}\|\blambda_{r+1}-\blambda_r\|^2
+\frac{2\vartheta}{\beta_r}\left(\frac{\gamma_{r-2}}{\gamma_{r-1}}-\frac{\gamma_{r-1}}{\gamma_r}\right)\|\blambda_r\|^2
\\\notag
&\quad +\left(\frac{2\vartheta}{\beta_{r-1}}-\frac{2\vartheta}{\beta_r}\right)\left(\frac{\gamma_{r-2}}{\gamma_{r-1}}-1\right)\|\blambda_{r}\|^2+\frac{6\vartheta(3-2\tau)}{\tau}\left(\frac{1}{\beta_{r-1}}-\frac{1}{\beta_{r}}\right)\|\blambda_{r+1}-\blambda_r\|^2
\\
&\quad +(\beta_{r+1}-\beta_r)\|g_1(\bx_{r+1})\|^2
+\left(\frac{1}{\beta_r}-\frac{1}{\beta_{r+1}}\right)(1-\tau)^2\|\blambda_{r+1}\|^2
+(1-\tau)\left(\frac{\gamma_{r-1}-\gamma_r}{2}\right)\|\blambda_{r+1}\|^2.\label{eq.whole2}
\end{align}
\end{proof}

\section{Optimality Gap}

\subsection{Proof of \leref{le.opt}}
\begin{proof}
Based on the definition of $\mathcal{G}(\bx,\blambda)$ in \eqref{eq.defopt} and the update rule of GDPA in \eqref{eq.alg}, we can have the upper bound of $\|\mathcal{G}(\bx_r, \blambda_r)\|$ in terms of $\|\bx_{r+1}-\bx_r\|$,  $\|\blambda_{r+1}-\blambda_r\|$, and $\|\blambda_{r+1}\|$.

\noindent{\bf First Step.} Notice that $[\beta_rg_i(\bx_r)+(1-\tau)\blambda_i]_+=0,\forall i\in \overbar{\mathcal{S}}_r$ and $\blambda_i\ge0$, meaning that $\beta_rg_i(\bx_r)+(1-\tau)\blambda_i\le 0$, i.e., $\beta_rg_i(\bx_r)\le -(1-\tau)\blambda_i$. So, we have $\blambda_i+ \beta_rg_i(\bx_r)\le\tau\blambda_i, \forall i\in \overbar{\mathcal{S}}_r$. Due to non-expansiveness of the projection operator, we have
\begin{equation}\label{eq.lam}
\|\bblambda'_{r+1}-\textrm{proj}_{\ge 0}\left(\bblambda_r+\beta_rg_2(\bx_r)\right)\|\le\tau\|\bblambda_{r}\|\mathop{=}\limits^{(a)}\tau\|\bblambda'_{r+1}-\bblambda_{r}\|,
\end{equation}
where
\begin{equation}
 [g_2(\bx_{r})]_i = \begin{cases} g_i(\bx_{r}), & \textrm{if}\quad i\in \overbar{\mathcal{S}}_r; \\ 0, & \textrm{otherwise}.\end{cases},
 \end{equation}
and in $(a)$ we use the fact $\bblambda'_{r+1}=0$.

\noindent{\bf Second Step.}  Then, we can have
\begin{align}
\notag
&\quad\;\|\mathcal{G}(\bx_r, \blambda_r)\|
\\\notag
&\le\frac{1}{\alpha_r}\|\bx_{r+1}-\bx_r\|+\frac{1}{\alpha_r}\left\|\bx_{r+1}-\textrm{proj}_{\mathcal{X}}(\bx_r-\alpha_r\nabla_{\bx} \mathcal{L}(\bx_r, \blambda_r))\right\|
\\
&\quad+\frac{1}{\beta_r}\|\blambda_{r+1}-\blambda_r\|+\frac{1}{\beta_r}\|\blambda_{r+1}-\textrm{proj}_{\ge 0}(\blambda_r+ \beta_r\nabla_{\blambda} \mathcal{L}(\bx_r, \blambda_r))\|
\\\notag
&\mathop{\le}\limits^{(a)}\frac{1}{\alpha_r}\|\bx_{r+1}-\bx_r\|
\\\notag
&\quad+\frac{1}{\alpha_r}\Bigg\|\textrm{proj}_{\mathcal{X}}\left(\bx_{r+1}-\alpha_r(\nabla_{\bx} f(\bx_{r})+J^{\T}(\bx_r)[(1-\tau)\blambda_r+\beta_rg(\bx_r)]_++\frac{1}{\alpha_r}(\bx_{r+1}-\bx_r))\right)
\\\notag
&\quad-\textrm{proj}_{\mathcal{X}}(\bx_r-\alpha_r(\nabla_{\bx} f(\bx_r)+J^{\T}(\bx_r)\blambda_r)\Bigg\|
\\\notag
&\quad+\frac{1}{\beta_r}\|\blambda_{r+1}-\blambda_r\|
\\\notag
&\quad+\frac{2}{\beta_r}\left\|\textrm{proj}_{\ge 0}\left(\blambda''_{r+1}+ \beta_r\left( g_1(\bx_{r+1})-\frac{1}{\beta_r}(\blambda''_{r+1}-\blambda'_r)-\gamma_r\blambda''_{r+1}\right)\right)-\textrm{proj}_{\ge 0}\left(\blambda'_r+  \beta_rg_1(\bx_r)\right)\right\|
\\
&\quad+\frac{2}{\beta_r}\left\|\bblambda'_{r+1}-\textrm{proj}_{\ge 0}\left(\bblambda_r+\beta_rg_2(\bx_r)\right)\right\|
\\\notag
&\mathop{\le}\limits^{(b)} \frac{3}{\alpha_r}\|\bx_{r+1}-\bx_r\|+\frac{3+2\tau}{\beta_r}\|\blambda_{r+1}-\blambda_r\|+2\gamma_r\|\blambda_{r+1}\|
+2\|g(\bx_{r+1})- g(\bx_{r})\|
\\\notag
&\quad+U_J\left\|[(1-\tau)\blambda_r+\beta_r g(\bx_r)]_+-\blambda_r\right\|
\\
&\mathop{\le}\limits^{(c)}\left(\frac{3}{\alpha_r}+2L_g\right)\|\bx_{r+1}-\bx_r\|+\frac{3+2\tau}{\beta_r}\|\blambda_{r+1}-\blambda_r\|+2\gamma_r\|\blambda_{r+1}\|
+U_J\left\|[(1-\tau)\blambda_r+\beta_r g(\bx_r)]_+-\blambda_r\right\|
\end{align}
where in $(a)$ we apply the optimality condition of the subproblems, i.e.,
\begin{subequations}
\begin{align}\label{eq.optc}
\bx_{r+1}&=\textrm{proj}_{\mathcal{X}}\left(\bx_{r+1}-\alpha_r\left(\nabla_{\bx} f(\bx_{r})+J^{\T}(\bx_r)[(1-\tau)\blambda_r+\beta_rg(\bx_r)]_++\frac{1}{\alpha_r}(\bx_{r+1}-\bx_r)\right)\right),
\\
\blambda''_{r+1}&=\textrm{proj}_{\ge0}\left(\blambda''_{r+1}+ \beta_r\left(g_1(\bx_{r+1})-\frac{1}{\beta_r}(\blambda''_{r+1}-\blambda'_r)-\gamma_r\blambda''_{r+1}\right)\right),
\end{align}
\end{subequations}
and the fact that $\sqrt{\|\bx\|^2+\|\by\|^2}\le\|\bx\|+\|\by\|$, in $(b)$ we use the triangle inequality, non-expansiveness of the projection operator, \eqref{eq.lam}, and $(c)$ is true due to the Lipschitz continuity.

\noindent{\bf Third Step.}

To quantify term $\left\|[(1-\tau)\blambda_r+\beta_r g(\bx_r)]_+-\blambda_r\right\|$, we first have
\begin{equation}
\left\|[(1-\tau)\bblambda_r+\beta_r g_2(\bx_r)]_+-\bblambda_r\right\|\mathop{=}\limits^{\eqref{eq.defs}}\|\bblambda_r\|=\|\bblambda'_{r+1}-\bblambda_r\|,\label{eq.ca1}
\end{equation}
since $\bblambda'_{r+1}=0$.

Second, for $i\in\mathcal{S}_r$, we have the following two cases:
\begin{enumerate}
\item[1)] when $(1-\tau)[\blambda_r]_i+\beta_rg_i(\bx_{r+1})\ge0$, then $(1-\tau)[\blambda_r]_i+\beta_rg_i(\bx_{r+1})=[\blambda_{r+1}]_i$, so we have
\begin{align}\notag
&\left|[(1-\tau)[\blambda'_r]_i+\beta_r g_i(\bx_r)]_+-[\blambda_r]_i\right|
\\
= & \left|[(1-\tau)[\blambda_r]_i+\beta_r g_i(\bx_{r+1})+\beta_r g_i(\bx_r)-\beta_r g_i(\bx_{r+1})]_+-[\blambda_r]_i\right|
\\
\mathop{\le}\limits^{(a)} & \left|(1-\tau)[\blambda_r]_i+\beta_r g_i(\bx_{r+1})+\beta_r g_i(\bx_r)-\beta_r g_i(\bx_{r+1})-[\blambda_r]_i\right|
\\
\mathop{\le}\limits^{(b)} & \left|\beta_r g_i(\bx_r)-\beta_r g_i(\bx_{r+1})\right|+\left|(1-\tau)[\blambda_r]_i+\beta_r g_i(\bx_{r+1})-[\blambda_r]_i\right|
\\
=& \left|\beta_r g_i(\bx_r)-\beta_r g_i(\bx_{r+1})\right|+\left|[\blambda_{r+1}]_i-[\blambda_r]_i\right|
\end{align}
where in $(a)$ we use non-expansiveness of the projection operator, in $(b)$ we use the triangle inequality;

\item[2)] when $(1-\tau)[\blambda_r]_i+\beta_rg_i(\bx_{r+1})<0$, then $[\blambda_{r+1}]_i=0$, and based on the definition of $\blambda_{r}$ we have $\beta_r(g_i(\bx_{r+1})-g_i(\bx_r))\ge0$, which gives
\begin{equation}\label{eq.bgr}
[(1-\tau)[\blambda_r]_i+\beta_rg_i(\bx_r)-\beta_rg_i(\bx_{r+1})+\beta_rg_i(\bx_{r+1})]_+\le \beta_rg_i(\bx_{r+1})-\beta_rg_i(\bx_{r}),
\end{equation}
so we have
\begin{align}\notag
&\left|[(1-\tau)[\blambda_r]_i+\beta_rg_i(\bx_r)-\beta_rg_i(\bx_{r+1})+\beta_rg_i(\bx_{r+1})]_+-[\blambda_r]_i\right|
\\
\mathop{\le}\limits^{\eqref{eq.bgr}}  & |\beta_rg_i(\bx_{r+1})-\beta_rg_i(\bx_{r})-[\blambda_r]_i|
\\
\mathop{\le}\limits^{(a)}  & |\beta_rg_i(\bx_{r+1})-\beta_rg_i(\bx_{r})+[\blambda_{r+1}]_i-[\blambda_r]_i|
\\
\mathop{\le}\limits^{(b)}  & |\beta_rg_i(\bx_{r+1})-\beta_rg_i(\bx_{r})|+|[\blambda_{r+1}]_i-[\blambda_r]_i|,
\end{align}
where in $(a)$ we use $[\blambda_{r+1}]_i=0$ in this case, in $(b)$ we use the triangle inequality.
\end{enumerate}

Combining the above two cases, we can get
\begin{equation}\label{eq.ca2}
\left\|[(1-\tau)\blambda'_r+\beta_rg_1(\bx_r)]_+-\blambda'_r\right\|
\le \beta_rL_g\|\bx_{r+1}-\bx_r\|+\|\blambda''_{r+1}-\blambda'_r\|
\end{equation}
where we use the Lipschitz continuity of function $g()$.

Combining \eqref{eq.ca1} and \eqref{eq.ca2} gives rise to
\begin{align}\notag
&\left\|[(1-\tau)\blambda_r+\beta_rg(\bx_r)]_+-\blambda_r\right\|^2
\\\notag
=&\left\|[(1-\tau)\blambda'_r+\beta_rg_1(\bx_r)]_+-\blambda'_r\right\|^2+\left\|[(1-\tau)\bblambda_r+\beta_rg_2(\bx_r)]_+-\bblambda_r\right\|^2
\\\notag
\le &\sum_{i\in\mathcal{S}_r}2L^2_g\beta^2_r\|[\bx_{r+1}]_i-[\bx_r]_i\|^2+2\|[\blambda_{r+1}]_i-[\blambda_r]_i\|^2
+\sum_{i\in\overbar{\mathcal{S}}_r}\|[\blambda_{r+1}]_i-[\blambda_r]_i\|^2
\\
\le & 2L^2_g\beta^2_r\|\bx_{r+1}-\bx_r\|^2+3\|\blambda_{r+1}-\blambda_r\|^2
\end{align}
where we use the fact $\|\bx\|^2\le2\|\by\|^2+2\|\bz\|^2$ when $\|\bx\|\le\|\by\|+\|\bz\|$.

Towards this end, we have
\begin{align}\notag
\|\mathcal{G}(\bx_r, \blambda_r)\|^2&\le 4\left(\frac{3}{\alpha_r}+2L_g\right)^2\|\bx_{r+1}-\bx_r\|^2+4\left(\frac{3+2\tau}{\beta_r}\right)^2\|\blambda_{r+1}-\blambda_r\|^2+16\gamma_r^2\|\blambda_{r+1}\|^2
\\
&\quad +4U^2_g\left(2L^2_g\beta^2_r\|\bx_{r+1}-\bx_r\|^2+3\|\blambda_{r+1}-\blambda_r\|^2\right)
\\\notag
&\le   4\left(\left(\frac{3}{\alpha_r}+2L_g\right)^2\!\!\!+2U^2_JL^2_g\beta^2_r\right)\|\bx_{r+1}-\bx_r\|^2
+4\left(\left(\frac{3+2\tau}{\beta_r}\right)^2\!\!\!+3\right)\|\blambda_{r+1}-\blambda_r\|^2+16\gamma_r^2\|\blambda_{r+1}\|^2.
\end{align}
\end{proof}

\subsection{Upper bound of dual variable}

\begin{mdframed}[backgroundcolor=gray!10,topline=false,
	rightline=false,
	leftline=false,
	bottomline=false]
\begin{lemma}\label{le.bdl}
Suppose that \asref{ass.bd} holds. If the iterates $\{\bx_r,\blambda_r,\forall r\}$ are generated by GDPA, where the step-sizes are chosen according to \eqref{eq.spcon}, then we have
\begin{equation}\label{eq.upbl}
\|\blambda_{r+1}\|^2\le   \frac{G \left(1+\frac{1}{\theta'}\right)r^{2/3}}{1-\eta}\sim\mathcal{O}(r^{2/3})
\end{equation}
where $0<\eta<1$ and $\theta'<1/(1-\tau)^2-1$.
\end{lemma}
\end{mdframed}

\begin{proof}
Recall \eqref{eq.defl}. Let $\blambda^g_{r+1}$ be a vector whose the $i$th entry is
\begin{equation}
 [\blambda^g_{r+1}]_i = \begin{cases} [\blambda_{r+1}]_i, & \textrm{if}\quad i\in\mathcal{A}_{r+1}\cap\mathcal{S}_r; \\ 0, & \textrm{otherwise}.\end{cases}
\end{equation}
where $\mathcal{A}_{r+1}\bydef\{i|g_i(\bx_{r+1})>0\}$.

Based on \eqref{eq.dualupated}, we have
\begin{equation}\label{eq.reub}
\|\blambda_{r+1}^g\|^2\le (1+\theta')(1-\tau)^2\|\blambda'^g_r\|^2+\left(1+\frac{1}{\theta'}\right)\beta^2_r\|g_+(\bx_{r+1})\|^2
\end{equation}
where notation $\blambda'^g_r$ means that the $i$th entry of $\blambda'^g_r$ is
\begin{equation}
 [\blambda'^g_r]_i = \begin{cases} [\blambda_r]_i, & \textrm{if}\quad i\in\mathcal{A}_{r+1}\cap\mathcal{S}_r; \\ 0, & \textrm{otherwise}.\end{cases}.
\end{equation}

When \asref{ass.bd} holds, i.e., $\|g_+(\bx_{r+1})\|^2\le G$, then from \eqref{eq.spcon}, we have
\begin{equation}\label{eq.rel}
\|\blambda_{r+1}^g\|^2\le (1+\theta')(1-\tau)^2\|\blambda'^g_r\|^2+\left(1+\frac{1}{\theta'}\right)Gr^{2/3}.
\end{equation}
Here, we can choose $\theta'<1/(1-\tau)^2-1$ so that $(1+\theta')(1-\tau)^2\bydef \eta' <1$. Then, we have
\begin{align}\notag
\|\blambda_{r+1}^g\|^2\le \eta'\|\blambda'^g_r\|^2+ \left(1+\frac{1}{\theta'}\right)G r^{2/3}.
\end{align}
Note that
\begin{subequations}
\begin{align}
[\blambda_{r+1}]_i\le&(1-\tau)[\blambda_r]_i, \forall i\in\overbar{\mathcal{A}}_{r+1}\cap\mathcal{S}_r,
\\
[\blambda_{r+1}]_i=&0, \forall i\in\overbar{\mathcal{S}}_r.\label{eq.shk2}
\end{align}
\end{subequations}
Combining \eqref{eq.rel}, \eqref{eq.shk} and \eqref{eq.shk2}, we have
\begin{align}
\|\blambda_{r+1}\|^2& \mathop{\le}\limits^{(a)} \eta\|\blambda_r\|^2+\left(1+\frac{1}{\theta'}\right)G r^{2/3}\label{eq.bdl}
\\
& \le \eta^2\|\blambda_{r-1}\|^2+\eta\left(1+\frac{1}{\theta'}\right)G r^{2/3} +\left(1+\frac{1}{\theta'}\right)G r^{2/3}
\\
&\le   \frac{G \left(1+\frac{1}{\theta'}\right)r^{2/3}}{1-\eta}
\end{align}
where $\eta\bydef\max\{(1-\tau)^2,\eta'\}$.
\end{proof}

\subsection{Regularity Condition: Proof of \leref{le.rege}}

\begin{proof}
Recall the fact that when $[g(\bx_{r+1})]_i\le 0$, we know  from \eqref{eq.dualupated}  that
\begin{equation}
[\blambda_{r+1}]_i\le(1-\tau)[\blambda_r]_i, \quad \forall i\in\overbar{\mathcal{A}}_{r+1}\cap\mathcal{S}_r.
\end{equation}
Due to the nonnegativity of $\blambda_r$ and $\tau\in(0,1)$. This implies that the size of the dual variable $[\blambda_{r+1}]_i$ is shrunk. Also, note that $[\blambda_{r+1}]_i=0, i\in\overbar{\mathcal{S}}_r$.

Upon the above results, we only need to consider the case, i.e., $[\blambda_{r+1}]_i=(1-\tau)[\blambda_{r}]_i+\beta_{r} g_i(\bx_{r+1}), i\in\mathcal{A}_{r+1}\cap\mathcal{S}_r$.
Here, $i$ denotes the index of these active constraints. Therefore, we only need to consider active constraints at $\bx_{r+1}$ and corresponding $\blambda'^g_r$ and $\blambda^g_{r+1}$, where notation $\blambda'^g_r$ means that the $i$th entry of $\blambda'^g_r$ is
\begin{equation}\label{eq.defl2}
 [\blambda'^g_r]_i = \begin{cases} [\blambda_r]_i, & \textrm{if}\quad i\in\mathcal{A}_{r+1}\cap\mathcal{S}_r; \\ 0, & \textrm{otherwise}.\end{cases}
\end{equation}

Let $\wblambda^g_{r+1}$ be a vector whose the $i$th entry is
\begin{equation}\label{eq.defln1}
 [\wblambda^g_{r+1}]_i = \begin{cases} [\blambda_{r+1}]_i, & \textrm{if}\quad i\in\mathcal{A}_{r+1}; \\ 0, & \textrm{otherwise}.\end{cases}
\end{equation}
and
\begin{equation}\label{eq.defln2}
 [\wblambda'^g_{r}]_i = \begin{cases} [\blambda_{r}]_i, & \textrm{if}\quad i\in\mathcal{A}_{r+1}; \\ 0, & \textrm{otherwise}.\end{cases}.
\end{equation}

By defining the following auxiliary variable
\begin{equation}
\wblambda^g_{r+1}\bydef(1-\tau)(\bx_{r+1})\wblambda'^g_r+\beta_r(\bx_{r+1})g_+(\bx_{r+1}),
\end{equation}
we can obtain 
\begin{align}\notag
&\quad\; \textrm{dist}\left(J^{\T}(\bx_{r+1})\wblambda^g_{r+1},-\mathcal{N}_{\mathcal{X}}(\bx_{r+1})\right)
\\
&=\textrm{dist}\left(\beta_rJ^{\T}(\bx_{r+1})g_+(\bx_{r+1})+(1-\tau)J^{\T}(\bx_{r+1})\wblambda'^g_r,-\mathcal{N}_{\mathcal{X}}(\bx_{r+1})\right)
\\
&\mathop{\ge}\limits^{(a)} \beta_r\textrm{dist}\left(J^{\T}(\bx_{r+1})g_+(\bx_{r+1}),-\mathcal{N}_{\mathcal{X}}(\bx_{r+1})\right)-\textrm{dist}\left((1-\tau)J^{\T}(\bx_{r+1})\wblambda'^g_r,-\mathcal{N}_{\mathcal{X}}(\bx_{r+1})\right)
\\
&\mathop{\ge}\limits^{(b)}  \beta_r\textrm{dist}(J^{\T}(\bx_{r+1})g_+(\bx_{r+1}),-\mathcal{N}_{\mathcal{X}}(\bx_{r+1}))-(1-\tau)U_J\|\blambda_r\|
\\
&\mathop{\ge}\limits^{(c)}  \sigma\beta_r\|g_+(\bx_{r+1})\|-(1-\tau)U_J\|\blambda_r\|\label{eq.disu}
\end{align}
where in $(a)$ we use the inverse triangle inequality, in $(b)$ we use the triangle inequality and convexity of feasible set $\mathcal{X}$, and in $(c)$ we apply the regularity condition \eqref{eq.re}.

Based on \eqref{eq.defln1}, \eqref{eq.defln2}, we have
\begin{align}\label{eq.upofg}
\|\blambda^g_{r+1}\|\le\|\wblambda^g_{r+1}\|\le(1-\tau)\|\wblambda'^g_r\|+\beta_r\|g_+(\bx_{r+1})\|\le(1-\tau)\|\blambda_r\|+\beta_r\|g_+(\bx_{r+1})\|.
\end{align}

Combining \eqref{eq.disu} and \eqref{eq.upofg}, it is sufficient to show $\frac{\sigma}{2}\|\blambda^g_{r+1}\|\le\textrm{dist}(J^{\T}(\bx_{r+1})\wblambda^g_{r+1},-\mathcal{N}_{\mathcal{X}}(\bx_{r+1}))$ if
\begin{align}
\frac{\sigma}{2}\left((1-\tau)\|\blambda_r\|+\beta_r\|g_+(\bx_{r+1})\|\right)\le \sigma\beta_r\|g_+(\bx_{r+1})\|-(1-\tau)U_J\|\blambda_r\|,
\end{align}
which requires
\begin{equation}
\beta_r\ge\frac{2(1-\tau)(\frac{\sigma}{2}+U_J)\|\blambda_r\|}{\sigma}.
\end{equation}

According to the upper bound of $\blambda$ given in \eqref{eq.upbl}, we have when
\begin{center}
\fbox{
\begin{minipage}{0.68\columnwidth}
\begin{equation}\label{eq.reofb1}
\beta_r\ge \underbrace{\frac{2(1-\tau)(\frac{\sigma}{2}+U_J)}{\sigma}\sqrt{\frac{G \left(1+\frac{1}{\theta'}\right)}{1-\eta}}}_{\beta_0}r^{1/3},
\end{equation}
\end{minipage}
}
\end{center}
then,
\begin{equation}
\frac{\sigma}{2}\|\blambda^g_{r+1}\|\le\textrm{dist}\left(J^{\T}(\bx_{r+1})\wblambda^g_{r+1},-\mathcal{N}_{\mathcal{X}}(\bx_{r+1})\right).
\end{equation}

Next, we need further to deal with the case where $i\in\mathcal{A}_{r+1}\cap\overbar{\mathcal{S}}_r$ in the sense that $g_i(\bx_r)\ge 0$ but $g_i(\bx_{r+1})>0$, i.e.,
\begin{align}\notag
&\quad\; \textrm{dist}\left(J^{\T}(\bx_{r+1})\wblambda^g_{r+1},-\mathcal{N}_{\mathcal{X}}(\bx_{r+1})\right)
\\
&\mathop{\le}\limits^{(a)}  \textrm{dist}\left(J^{\T}(\bx_{r+1})\overbar{\blambda}'^g_{r+1},-\mathcal{N}_{\mathcal{X}}(\bx_{r+1}))+\textrm{dist}(J^{\T}(\bx_{r+1})\blambda^g_{r+1},-\mathcal{N}_{\mathcal{X}}(\bx_{r+1})\right)
\\
&\mathop{\le}\limits^{(b)}   (1-\tau)U_J\|\overbar{\blambda}''^g_{r}\|+U_J\beta_r\|\overbar{g}(\bx_{r+1})\|+\textrm{dist}\left(J^{\T}(\bx_{r+1})\blambda^g_{r+1},-\mathcal{N}_{\mathcal{X}}(\bx_{r+1})\right)
\\
&\mathop{\le}\limits^{(c)}   (1-\tau)U_J\|\overbar{\blambda}''^g_{r}\|+U_J\beta_rL_g\|\bx_{r+1}-\bx_{r}\|+\textrm{dist}\left(J^{\T}(\bx_{r+1})\blambda^g_{r+1},-\mathcal{N}_{\mathcal{X}}(\bx_{r+1})\right)\label{eq.reg}
\end{align}
where in $(a)$ we use the triangle inequality and the following definition of $\overbar{\blambda}'^g_{r+1}$, i.e.,
\begin{equation}
 [\overbar{\blambda}'^g_{r+1}]_i = \begin{cases} [\blambda_{r+1}]_i, & \textrm{if}\quad i\in\mathcal{A}_{r+1}\cap\overbar{\mathcal{S}}_r; \\ 0, & \textrm{otherwise}.\end{cases},
\end{equation}
in $(b)$ we define
\begin{align}
 [\overbar{\blambda}''^g_{r}]_i =& \begin{cases} [\blambda_{r}]_i, & \textrm{if}\quad i\in\mathcal{A}_{r+1}\cap\overbar{\mathcal{S}}_r; \\ 0, & \textrm{otherwise}.\end{cases},
 \\
 \overbar{g}(\bx_{r+1}) =& \begin{cases} g_i(\bx_{r+1}), & \textrm{if}\quad i\in\mathcal{A}_{r+1}\cap\overbar{\mathcal{S}}_r; \\ 0, & \textrm{otherwise}.\end{cases},
\end{align}
and use
\begin{equation}
 [\overbar{\blambda}'^g_{r+1}]_i = (1-\tau) [\overbar{\blambda}''^g_{r}]_i +\beta_r\overbar{g}_i(\bx_{r+1}),\quad i\in\mathcal{A}_{r+1}\cap\overbar{\mathcal{S}}_r,
\end{equation}
and $(c)$ is true because we define
\begin{equation}
 \overbar{g}'(\bx_{r}) = \begin{cases} g_i(\bx_{r})\le0, & \textrm{if}\quad i\in\mathcal{A}_{r+1}\cap\overbar{\mathcal{S}}_r; \\ 0, & \textrm{otherwise}.\end{cases},
 \end{equation}
 and apply
\begin{equation}
\left\|[\overbar{g}(\bx_{r+1})]_+-[\overbar{g}'(\bx_r)]_+\right\|\le L_g\|\bx_{r+1}-\bx_r\|.
\end{equation}

Given condition \eqref{eq.reg}, we will show the recursion of $\blambda^g_r$ as follows.

From \eqref{eq.updatex}, we have
\begin{equation}\label{eq.disopt}
\textrm{dist}\left((\nabla f(\bx_r) +J^{\T}(\bx_r)[(1-\tau)\blambda_r+\beta_r g(\bx_{r})]_+ +\frac{1}{\alpha_r}(\bx_{r+1}-\bx_{r}),-\mathcal{N}_{\mathcal{X}}(\bx_{r+1})\right)=0,
\end{equation}
then we re-order some terms  and get
\begin{align}
\notag
&\quad\;\nabla f(\bx_r) +J^{\T}(\bx_r)[(1-\tau)\blambda_r+\beta_r g(\bx_{r})]_+ +\frac{1}{\alpha_r}(\bx_{r+1}-\bx_{r})
\\\notag
&=\nabla f(\bx_r) +J^{\T}(\bx_r)[(1-\tau)\blambda_r+\beta_r g(\bx_{r})]_++J^{\T}(\bx_{r+1})\blambda_{r+1}^g -J^{\T}(\bx_{r+1})\blambda_{r+1}^g
\\
&\quad+J^{\T}(\bx_{r})\blambda_{r+1}-J^{\T}(\bx_{r})\blambda_{r+1}+\frac{1}{\alpha_r}(\bx_{r+1}-\bx_{r})\label{eq.updxu}
\\\notag
&\mathop{=}\limits^{(a)} \nabla f(\bx_r) +J^{\T}(\bx_r)[(1-\tau)\blambda_r+\beta_r g(\bx_{r})]_+-J^{\T}(\bx_{r})\blambda_{r+1}+J^{\T}(\bx_{r+1})\blambda_{r+1}^g
\\
&\quad+J^{\T}(\bx_{r})\blambda_{r+1}-J^{\T}(\bx_{r+1})\blambda'^g_{r}+J^{\T}(\bx_{r+1})\blambda'^g_{r}-J^{\T}(\bx_{r+1})\blambda_{r+1}^g
+\frac{1}{\alpha_r}(\bx_{r+1}-\bx_{r})
\end{align}
where in $(a)$ we add and subtract some same terms.

For simplicity in equation expressions, we define
\begin{align}\notag
&\nabla f(\bx_r) +J^{\T}(\bx_r)[(1-\tau)\blambda_r+\beta_r g(\bx_{r})]_+-J^{\T}(\bx_{r})\blambda_{r+1}
\\
&+J^{\T}(\bx_{r})\blambda_{r+1}-J^{\T}(\bx_{r+1})\blambda'^g_{r}+J^{\T}(\bx_{r+1})\blambda'^g_{r}-J^{\T}(\bx_{r+1})\blambda_{r+1}^g
+\frac{1}{\alpha_r}(\bx_{r+1}-\bx_{r})\bydef \bb.\label{eq.jk}
\end{align}

So, we can write \eqref{eq.disopt} as
\begin{equation}
\textrm{dist}\left(J^{\T}(\bx_{r+1})\blambda_{r+1}^g+\bb,-\mathcal{N}_{\mathcal{X}}(\bx_{r+1})\right)=0.
\end{equation}
Using the inverse triangle inequality, we have
\begin{equation}
0=\textrm{dist}\left(J^{\T}(\bx_{r+1})\blambda_{r+1}^g+\bb,-\mathcal{N}_{\mathcal{X}}(\bx_{r+1})\right)\ge \textrm{dist}\left(J^{\T}(\bx_{r+1})\blambda_{r+1}^g,-\mathcal{N}_{\mathcal{X}}(\bx_{r+1})\right)-\|\bb\|.
\end{equation}
Applying the regularity condition \eqref{eq.reg}, we have
\begin{equation}\label{eq.disre}
\frac{\sigma}{2}\|\blambda_{r+1}^g\| \le \textrm{dist}\left(J^{\T}(\bx_{r+1})\blambda_{r+1}^g,-\mathcal{N}_{\mathcal{X}}(\bx_{r+1})\right)\le\|\bb\|+(1-\tau)U_J\|\overbar{\blambda}''^g_{r}\|+U_J\beta_rL_g\|\bx_{r+1}-\bx_{r}\|.
\end{equation}

Substituting \eqref{eq.jk} into $\|\bb\|$, we have
\begin{align}\notag
\|\bb\|&\mathop{\le}\limits^{(a)} \|\nabla f(\bx_r)\|+\frac{1}{\alpha_r}\|\bx_{r+1}-\bx_{r}\|+\underbrace{(1-\tau)U_J\|\blambda_r\|+ U_J\|\blambda_{r+1}-\blambda_r\|
+L_J\|\blambda_r\|\|\bx_{r+1}-\bx_r\|}_{iii}
\\
&\quad +\underbrace{U_J\beta_rL_g\|\bx_{r+1}-\bx_r\|}_{i}+ \underbrace{U_J\|\blambda_{r+1}-\blambda_r\|}_{ii}
\\
&\le   M+\left(\frac{1}{\alpha_r}+U_J\beta_rL_g+L_J\|\blambda_r\|\right)\|\bx_{r+1}-\bx_{r}\|+ 2U_J\|\blambda_{r+1}-\blambda_r\|+(1-\tau)U_J\|\blambda_r\|\label{eq.bdj}
\end{align}
where in $(a)$ we use \asref{ass.bd}, i.e., $ \|\nabla f(\bx) \|\le M,\forall \bx$, and use the following facts to bound the other terms in \eqref{eq.jk}:

\emph{i}) apply gradient Lipschitz continuity so that we have 
\begin{align}\notag
&\left\|[(1-\tau)\blambda_r+\beta_r g(\bx_{r})]_+-\blambda_{r+1}\right\|
\\
\mathop{=}\limits^{(a)}&\left\|[(1-\tau)\blambda'_r+\beta_r g_1(\bx_{r})]_+-[(1-\tau)\blambda'_r+\beta_r g_1(\bx_{r+1})]_+\right\|
\\
\mathop{\le}\limits^{(b)} & \beta_r L_g\|\bx_{r+1}-\bx_r\|;
\end{align}
where $(a)$ is true because of the definition of $\blambda''_{r+1}$ introduced in \eqref{eq.defs} and \eqref{eq.dualupated}, in $(b)$ we use non-expansiveness of the projection operator and the triangle inequality;

\emph{ii}) use the boundedness of the Jacobian matrix so we can get
\begin{equation}\notag
 \|J^{\T}(\bx_{r+1})\blambda'^g_{r}-J^{\T}(\bx_{r+1})\blambda^g_{r+1}\|
\le  U_J\|\blambda_{r+1}-\blambda_r\|;
\end{equation}

\emph{iii})  we re-organize terms as follows
\begin{align}\notag
&\quad\;J^{\T}(\bx_{r})\blambda_{r+1}-J^{\T}(\bx_{r+1})\blambda'^g_{r}
\\
&=J^{\T}(\bx_{r})\blambda_{r+1}-J^{\T}(\bx_r)\blambda_{r+1}^g+J^{\T}(\bx_r)\blambda_{r+1}^g-J^{\T}(\bx_{r+1})\blambda'^g_{r}
\\\notag
&=J^{\T}(\bx_{r})\blambda_{r+1}-J^{\T}(\bx_r)\blambda_{r+1}^g+J^{\T}(\bx_r)\blambda_{r+1}^g-J^{\T}(\bx_r)\blambda'^g_{r}
+J^{\T}(\bx_r)\blambda'^g_{r}-J^{\T}(\bx_{r+1})\blambda'^g_{r},
\end{align}
which can give us the following inequalities directly
\begin{align}\notag
&\quad\;\|J^{\T}(\bx_{r})\blambda_{r+1}-J^{\T}(\bx_{r+1})\blambda'^g_{r}\|
\\
&\le   \|J^{\T}(\bx_{r})\blambda_{r+1}-J^{\T}(\bx_r)\blambda_{r+1}^g\| + U_J\|\blambda_{r+1}-\blambda_r\|+L_J\|\blambda_r\|\|\bx_{r+1}-\bx_r\|
\\
&\mathop{\le}\limits^{(a.1)}  (1-\tau)U_J\|\blambda_r\|+ U_J\|\blambda_{r+1}-\blambda_r\|+L_J\|\blambda_r\|\|\bx_{r+1}-\bx_r\|
\end{align}
where $(a.1)$ is true because based on the definition of $\blambda^g_{r+1}$, and we also use the upper bound of $\|J^{\T}(\bx_{r})\blambda_{r+1}-J^{\T}(\bx_r)\blambda_{r+1}^g\|$, which is obtained by the following steps:
\begin{align}\notag
&\|J^{\T}(\bx_{r})\blambda_{r+1}-J^{\T}(\bx_r)\blambda_{r+1}^g\|
\\
\mathop{=}\limits^{(a.2)}&U_J\left\|[(1-\tau)\blambda'_r+\beta_rg_1(\bx_{r+1})]_+-[(1-\tau)\blambda'^g_r+\beta_rg'(\bx_{r+1})]_+\right\|
\\
\mathop{\le}\limits^{(a.3)} &(1-\tau)U_J\|\blambda_r\|\label{eq.upj}
\end{align}
where $(a.2)$ is true due to $\bblambda'_{r+1}=0$, the notation $g'$ follows \eqref{eq.defl}, i.e., the $i$th entry of $g'(\bx_{r+1})$ is
\begin{equation}
 [g'(\bx_{r+1})]_i = \begin{cases} g_i(\bx_{r+1}), & \textrm{if}\quad i\in\mathcal{A}_{r+1}\cap\mathcal{S}_r; \\ 0, & \textrm{otherwise}.\end{cases}
\end{equation}
and $(a.3)$ is true since we use the following facts
\begin{enumerate}
\item When $i\in\mathcal{A}_{r+1}\cap\mathcal{S}_r$, then $ [(1-\tau)[\blambda'_r]_i+\beta_r[g_1]_i(\bx_{r+1})]_+-[(1-\tau)[\blambda'^g_r]_i+\beta_rg'_i(\bx_{r+1})]_+ =0$.

\item When  $i\in\overbar{\mathcal{A}}_{r+1}\cap \mathcal{S}_r$, i.e., $[g(\bx_{r+1})]_i\le0$ and $(1-\tau)[\blambda_r]_i+\beta_rg_i(\bx_{r})\ge0$, then we have
\begin{align}\notag
&\left|[(1-\tau)[\blambda'_r]_i+\beta_r[g_1]_i(\bx_{r+1})]_+-[(1-\tau)[\blambda'^g_r]_i+\beta_rg'_i(\bx_{r+1})]_+\right|
\\
=&\left|[(1-\tau)[\blambda'_r]_i+\beta_r[g_1]_i(\bx_{r+1})]_+\right|\le(1-\tau)[\blambda'_r]_i, \forall i\in\overbar{\mathcal{A}}_{r+1}\cap \mathcal{S}_r.
\end{align}
In summary, we have \eqref{eq.upj}.
\end{enumerate}

Substituting \eqref{eq.bdj} into \eqref{eq.disre}, we know that
\begin{align}
\frac{\sigma}{2}\|\blambda_{r+1}^g\|\le
M+\left(\frac{1}{\alpha_r}+2U_J\beta_rL_g+L_J\|\blambda_r\|\right)\|\bx_{r+1}-\bx_{r}\|+ 2U_J\|\blambda_{r+1}-\blambda_r\|+2(1-\tau)U_J\|\blambda_r\|.
\label{eq.bdj2}
\end{align}

Combining the regularity condition \eqref{eq.reg}, we have that
\begin{multline}
\frac{\sigma^2}{4}\|\blambda_{r+1}^g\|^2
\le  4M^2+4\left(\frac{1}{\alpha_r}+2U_J\beta_rL_g+L_J\|\blambda_r\|\right)^2\|\bx_{r+1}-\bx_{r}\|^2
\\
+16U^2_J\|\blambda_{r+1}-\blambda_r\|^2+16(1-\tau)^2U^2_J\|\blambda_r\|^2.\label{eq.itof2}
\end{multline}
Note that
\begin{equation}
\frac{1}{\alpha_r}+2U_J\beta_rL_g+L_J\|\blambda_r\|=\frac{1+2U_JL_g\beta_r\alpha_r+\alpha_rL_J\|\blambda_r\|}{\alpha_r}\mathop{<}\limits^{\eqref{eq.conal}}\frac{3}{\alpha_r}.
\end{equation}
Therefore, we have
\begin{equation}
\sigma^2\|\blambda_{r+1}^g\|^2\le 16M^2+\frac{36}{\alpha^2_r}\|\bx_{r+1}-\bx_{r}\|^2+64U^2_J\|\blambda_{r+1}-\blambda_r\|^2+64(1-\tau)^2U^2_J\|\blambda_r\|^2.
\end{equation}
The proof is complete.
\end{proof}
\subsection{Upper bound of Sum of Dual Variables: proof of \leref{le.sumd}}
\begin{proof}
We will use mathematical induction to obtain the upper bound of $\sum^T_{r=1}\alpha^2_r\|\blambda_r\|^2$.

It is trivial to show the case when $T=1$. Then, we assume that
\begin{equation}\label{eq.min}
\sum^R_{r=1} \alpha^2_r \|\blambda_r\|^2\sim\mathcal{O}(\beta_R),  \forall R\le T,
\end{equation}
i.e., there exists a constant $\Lambda$ such that
\begin{equation}\label{eq.inl}
\sum^R_{r=1} \alpha^2_r \|\blambda_r\|^2\le \Lambda \beta_R, \forall R.
\end{equation}

In the following, we will show that
\begin{equation}\label{eq.inl2}
\sum^R_{r=1}\alpha^2_{r+1}\|\blambda_{r+1}\|^2\le \Lambda\beta_{R}\le\Lambda\beta_{R+1} , \forall R.
\end{equation}

\noindent{\bf Step 1: Upper bound of the size of the difference of two successive primal and dual variables}

Using \eqref{eq.recc} as shown in \leref{le.cpo}, we have
\begin{align} \notag
&\quad\; \frac{1}{2\beta_{r}}  \|\blambda_{r+1}-\blambda_r\|^2 +\frac{1}{4\alpha_r}\|\bx_{r+1}-\bx_r\|^2
\\\notag
&\mathop{\le}\limits^{(a)} \mathcal{P}_r-\mathcal{P}_{r+1}+\frac{6\vartheta(3-2\tau)}{\tau}\left(\frac{1}{\beta_{r-1}}-\frac{1}{\beta_{r}}\right)\|\blambda_{r+1}-\blambda_r\|^2
\\\notag
&\quad +\left(\frac{2\vartheta}{\beta_{r-1}}-\frac{2\vartheta}{\beta_r}\right)\left(\frac{\gamma_{r-2}}{\gamma_{r-1}}-1\right)\|\blambda_{r}\|^2
+\frac{2\vartheta}{\beta_{r-1}}\frac{\gamma_{r-2}-\gamma_{r-1}}{\gamma_{r-1}}\|\blambda_r\|^2
\\
&\quad +(\beta_{r+1}-\beta_r)\|g_1(\bx_{r+1})\|^2
+\left(\frac{1}{\beta_r}-\frac{1}{\beta_{r+1}}\right)(1-\tau)^2\|\blambda_{r+1}\|^2
+(1-\tau)\left(\frac{\gamma_{r-1}-\gamma_r}{2}\right)\|\blambda_{r+1}\|^2
\\\notag
&\mathop{\le}\limits^{(b)} \mathcal{P}_r-\mathcal{P}_{r+1}+\frac{6\vartheta(3-2\tau)}{\tau}\left(\frac{1}{\beta_{r-1}}-\frac{1}{\beta_{r}}\right)\|\blambda_{r+1}-\blambda_r\|^2
  +2\vartheta\left(\frac{1}{\beta_{r-1}}-\frac{1}{\beta_r}\right)\|\blambda_{r}\|^2
+ 2\vartheta \left(\frac{1}{\beta_{r-2}}-\frac{1}{\beta_{r-1}}\right)\|\blambda_r\|^2
\\
&\quad +(\beta_{r+1}-\beta_r)\|g_1(\bx_{r+1})\|^2
+\left(\frac{1}{\beta_r}-\frac{1}{\beta_{r+1}}\right)(1-\tau)^2\|\blambda_{r+1}\|^2
+(1-\tau)\left(\frac{\gamma_{r-1}-\gamma_r}{2}\right)\|\blambda_{r+1}\|^2
\\\notag
&\mathop{\le}\limits^{(c)} \mathcal{P}_r-\mathcal{P}_{r+1}+\left(\frac{6\vartheta(3-2\tau)}{\tau}\left(\frac{1}{\beta_{r-1}}-\frac{1}{\beta_{r}}\right)+4\left(\frac{1}{\beta_{r}}-\frac{1}{\beta_{r+1}}\right)\right)\|\blambda_{r+1}-\blambda_r\|^2
\\\notag
&\quad +\left(2\vartheta\left(\frac{1}{\beta_{r-1}}-\frac{1}{\beta_r}\right)
+ 2\vartheta \left(\frac{1}{\beta_{r-2}}-\frac{1}{\beta_{r-1}}\right)+4\tau^2\left(\frac{1}{\beta_{r}}-\frac{1}{\beta_{r+1}}\right)\right)\|\blambda_r\|^2
\\
&\quad +\left(\left(\frac{1}{\beta_r}-\frac{1}{\beta_{r+1}}\right)(1-\tau)^2
+\frac{(1-\tau)\tau}{2}\left(\frac{1}{\beta_{r-1}}-\frac{1}{\beta_r}\right)\right)\|\blambda_{r+1}\|^2\label{eq.updiffbeta}
\end{align}
where in $(a)$ as $\gamma_r$ is a decreasing sequence we have
\begin{equation}
\left(\frac{\gamma_{r-2}}{\gamma_{r-1}}-\frac{\gamma_{r-1}}{\gamma_r}\right)\le\frac{\gamma_r\gamma_{r-2}}{\gamma^2_{r-1}}-1\le\frac{\gamma_{r-2}}{\gamma_{r-1}}-1\le\frac{\gamma_{r-2}-\gamma_{r-1}}{\gamma_{r-1}},
\end{equation}
and we also use $\beta_r$ is a increasing sequence and $\gamma_{r-2}/\gamma_{r-1}\le2$, $(b)$ is true due to $0<\gamma_r\beta_r=\tau<1, \forall r$, and in $(c)$
from \eqref{eq.fb2} and \eqref{eq.dualupated}, we know that
\begin{equation}
\|g_1(\bx_{r+1})\|^2\le\frac{2}{\beta^2_r}\left(\|\blambda_{r+1}-\blambda_r\|^2+\tau^2\|\blambda_r\|^2\right).
\end{equation}
Thus, we have
\begin{align}
(\beta_{r+1}-\beta_r)\|g_1(\bx_{r+1})\|^2
\le&\frac{2(\beta_{r+1}-\beta_r)}{\beta^2_r}\left(\|\blambda_{r+1}-\blambda_r\|^2+\tau^2\|\blambda_r\|^2\right)
\\
\le&2\left(\frac{1}{\beta_{r}}-\frac{1}{\beta_{r+1}}\right)\frac{\beta_{r+1}}{\beta_r}\left(\|\blambda_{r+1}-\blambda_r\|^2+\tau^2\|\blambda_r\|^2\right)
\\
\le&4\left(\frac{1}{\beta_{r}}-\frac{1}{\beta_{r+1}}\right)\left(\|\blambda_{r+1}-\blambda_r\|^2+\tau^2\|\blambda_r\|^2\right)
\end{align}
where the last inequality is true because $\beta_r$ is an increasing sequence and  chosen as $\mathcal{O}(r^{1/3})$, e.g., $\beta_r=\beta_0r^{1/3}$.

Note that  we choose stepsizes $\beta_r$ as \eqref{eq.spcon}, then
\begin{align}
&\frac{1}{\beta_{r}}-\frac{1}{\beta_{r+1}}\sim  \gamma_r-\gamma_{r+1}\sim \alpha_r-\alpha_{r+1}
\\
\sim&\frac{1}{r^{1/3}}-\frac{1}{(r+1)^{1/3}}\le\frac{(r+1)^{1/3}-r^{1/3}}{r^{2/3}}\le\frac{1}{3r^{4/3}}\mathop{\sim}\limits^{(a)} \mathcal{O}\left(\frac{1}{r^{4/3}}\right)\sim \alpha^4_r\sim\gamma_r^4\sim\frac{1}{\beta^4_r}\label{eq.ssr2}
\end{align}
where in $(a)$ we use the gradient Lipschitz continuity of function $x^{1/3}$, i.e., $(x+y)^{1/3}-x^{1/3}\le y/(3x^{2/3})$ for any $x,y\ge 1$.

Given these results,  we can have if $\beta_r=\beta_0r^{1/3}$, then $1/\beta_{r}-1/\beta_{r+1}\le \frac{1}{\beta_r 3r}$, $1/\beta_{r-1}-1/\beta_{r}\le \frac{1}{\beta_{r-1}3(r-1)}\le\frac{2}{\beta_r3(r-1)}$ and $1/\beta_{r-2}-1/\beta_{r-1}\le \frac{4}{\beta_r3(r-2)}$. It is easy to check that when $r\ge3>8/3$, then $4(1/\beta_{r}-1/\beta_{r+1})<1/(2\beta_r)$ and when $r>8\vartheta(3-2\tau)/\tau +1$, then $(6\vartheta(3-2\tau))/\tau(1/\beta_{r-1}-1/\beta_{r})<1/(2\beta_r)$\footnote{Note that running a constant number of iterations as a warm start for an algorithm will not affect the iteration complexity.}. Therefore, multiplying $4$ on both sides of \eqref{eq.updiffbeta}, we can further simplify it as
\begin{align}
\notag
&\quad\;\frac{1}{\beta_{r}}\|\blambda_{r+1}-\blambda_r\|^2+\frac{1}{\alpha_r}\|\bx_{r+1}-\bx_r\|^2
\\\notag
&\le4(\mathcal{P}_r-\mathcal{P}_{r+1})+ 4\left(2\vartheta\left(\frac{1}{\beta_{r-1}}-\frac{1}{\beta_r}\right)
+ 2\vartheta \left(\frac{1}{\beta_{r-2}}-\frac{1}{\beta_{r-1}}\right)+4\tau^2\left(\frac{1}{\beta_{r}}-\frac{1}{\beta_{r+1}}\right)\right)\|\blambda_r\|^2
\\
&\quad +4\left(\left(\frac{1}{\beta_r}-\frac{1}{\beta_{r+1}}\right)(1-\tau)^2
+\frac{(1-\tau)\tau}{2}\left(\frac{1}{\beta_{r-1}}-\frac{1}{\beta_r}\right)\right)\|\blambda_{r+1}\|^2
 \\
&\mathop{\le}\limits^{(a)}4(\mathcal{P}_r-\mathcal{P}_{r+1})+4\left(\frac{1}{\beta_03r^{4/3}}+\frac{(3/2)^{4/3}}{\beta_06r^{4/3}}\right)\|\blambda_{r+1}\|^2
+\left(8\vartheta\frac{\frac{3}{2}^{4/3}+2^{4/3}}{\beta_03r^{4/3}}+\frac{4}{\beta_03r^{4/3}}\right)\|\blambda_{r}\|^2
\\
&\mathop{\le}\limits^{(b)}4(\mathcal{P}_r-\mathcal{P}_{r+1})+\frac{2.5}{\beta_0r^{4/3}}\|\blambda_{r+1}\|^2+\frac{13\vartheta}{\beta_0r^{4/3}}\|\blambda_{r}\|^2\label{eq.updiffbeta2}
\end{align}
where in $(a)$ we use \eqref{eq.ssr2} again and $r\ge3$, and $(b)$ is true by a direct numerical calculation.

\noindent{\bf Step 2: Recursion of upper bound of dual variable}

Multiplying $\alpha_r$ on both sides of \eqref{eq.itofl} shown in \leref{le.rege}, we have
\begin{align}\notag
&\quad\;\alpha_r\sigma^2\|\blambda_{r+1}^g\|^2
\\
&\le 16\alpha_rM^2+\frac{36}{\alpha_r}\|\bx_{r+1}-\bx_r\|^2 + \alpha_r64U^2_J\|\blambda_{r+1}-\blambda_r\|^2
+64\alpha_r(1-\tau)^2U^2_J\|\blambda_r\|^2
\\
&\mathop{\le}\limits^{\eqref
{eq.conal}} 16\alpha_rM^2+\max\left\{36,\frac{64U_J}{L_g}\right\}\left(\frac{1}{\alpha_r}\|\bx_{r+1}-\bx_r\|^2 + \frac{1}{\beta_r}\|\blambda_{r+1}-\blambda_r\|^2\right)
+64\alpha_r(1-\tau)^2U^2_J\|\blambda_r\|^2.\label{eq.sl}
\end{align}

Then, substituting \eqref{eq.updiffbeta2} into \eqref{eq.sl} gives
\begin{align}\notag
&\alpha_{r+1}\sigma^2\|\blambda_{r+1}^g\|^2
\le\alpha_r\sigma^2\|\blambda_{r+1}^g\|^2
\\
\le& 16\alpha_rM^2
+\max\left\{36,\frac{64U_J}{L_g}\right\}\bigg(4(\mathcal{P}_r-\mathcal{P}_{r+1})+\frac{2.5}{\beta_0r^{4/3}}\|\blambda_{r+1}\|^2+\frac{13\vartheta}{\beta_0r^{4/3}}\|\blambda_{r}\|^2
\bigg)+64\alpha_r(1-\tau)^2U^2_J\|\blambda_r\|^2 \label{eq.tel}
\end{align}
where we use the fact that $\alpha_r$ is a decreasing sequence.

We further multiply another $\alpha_r$ on both sides of  \eqref{eq.tel} and can get
\begin{multline}
\alpha^2_{r+1}\sigma^2\|\blambda_{r+1}^g\|^2
\le 16M^2 \alpha^2_r +64(1-\tau)^2U^2_J \alpha^2_r \|\blambda_r\|^2
\\
+ \alpha_r \max\left\{36,\frac{64U_J}{L_g}\right\}\left(4(\mathcal{P}_r-\mathcal{P}_{r+1})+\frac{2.5}{\beta_0r^{4/3}}\|\blambda_{r+1}\|^2+\frac{13\vartheta}{\beta_0r^{4/3}}\|\blambda_{r}\|^2
\right) \label{eq.tel2}
\end{multline}
where we use the fact that $\alpha_r$ is a decreasing sequence again.

\noindent{\bf Step 3: Telescoping sum of dual variables}

Recall the potential function: from \eqref{eq.defF} and \eqref{eq.defp}, we know
\begin{align}\notag
\mathcal{P}_r  = & f(\bx_r) +\frac{\beta_r}{2}\left\|\left[g(\bx_r)+\frac{(1-\tau)\blambda_r}{\beta_r}\right]_+\right\|^2-\frac{\|(1-\tau)\blambda_r\|^2}{2\beta_r}
\\\notag
&+(1-\tau)\left(\frac{1}{2\beta_{r-1}}+\frac{2\vartheta}{\beta_{r-2}\tau}\right)\|\blambda_{r}-\blambda_{r-1}\|^2
-\frac{(1-\tau)\gamma_{r-1}}{2}\|\blambda_r\|^2-\frac{2\vartheta}{\beta_{r-1}}\left(\frac{\gamma_{r-2}}{\gamma_{r-1}}-1\right)\|\blambda_r\|^2
\\
&+\beta_{r-1}L^2_g\|\bx_{r}-\bx_{r-1}\|^2+\frac{4\vartheta L^2_{g}}{\gamma_{r-1}(1-\tau)}\|\bx_{r}-\bx_{r-1}\|^2.\label{eq.dep}
\end{align}

Applying the telescoping sum of \eqref{eq.tel2} over $r=1,\ldots, R$, we have from \eqref{eq.dep}
\begin{align}\notag
&\quad\;\sigma^2\sum^R_{r=1} \alpha^2_{r+1} \|\blambda_{r+1}^g\|^2
\\\notag
&\le 
16M^2\sum^R_{r=1} \alpha^2_r
+\max\left\{36,\frac{64U_J}{L_g}\right\} 4\alpha_1\Delta_F
\\\notag
&\quad+\max\left\{36,\frac{64U_J}{L_g}\right\} 4\alpha_1 \left(\frac{(1-\tau)^2}{\beta_{R+1}}\|\blambda_{R+1}\|^2
+\frac{(1-\tau)\gamma_{R}}{2}\|\blambda_{R+1}\|^2+\frac{2\vartheta}{\beta_{R}}\left(\frac{\gamma_{R-1}}{\gamma_{R}}-1\right)\|\blambda_{R+1}\|^2\right)
\\\notag
&\quad+ \max\left\{36,\frac{64U_J}{L_g}\right\}\sum^R_{r=1}\frac{2.5}{\beta_0 r^{4/3}}\alpha_r\|\blambda_{r+1}\|^2+\max\left\{36,\frac{64U_J}{L_g}\right\}\sum^R_{r=1}\frac{13\vartheta}{\beta_0 r^{4/3}}\alpha_r\|\blambda_{r}\|^2
\\
&\quad+64(1-\tau)^2U^2_J\sum^R_{r=1} \alpha^2_r \|\blambda_r\|^2 \label{eq.sumofp}
\\\notag
&\mathop{\le}\limits^{(a)}   4M^2\sum^R_{r=1} \alpha^2_r
+4\max\left\{36,\frac{64U_J}{L_g}\right\} \alpha_1 \Delta_F+ 4\max\left\{36,\frac{64U_J}{L_g}\right\} \alpha_1 \frac{4G \left(1+\frac{1}{\theta'}\right)R^{1/3}}{(1-\eta)\beta_0}
\\\notag
&\quad+\max\left\{36,\frac{64U_J}{L_g}\right\}\sum^R_{r=1}\frac{2.5}{\beta_0 r^{4/3}}\alpha_r\|\blambda_{r+1}\|^2+\max\left\{36,\frac{64U_J}{L_g}\right\}\sum^R_{r=1}\frac{13\vartheta}{\beta_0 r^{4/3}}\alpha_r\|\blambda_{r}\|^2
\\
&\quad+64(1-\tau)^2U^2_J\sum^R_{r=1} \alpha^2_r \|\blambda_r\|^2\label{eq.lup}
\end{align}
where in $(a)$ we have that
\begin{multline}
\Delta_F\bydef f(\bx_1) +\frac{\beta_1}{2}\left\|\left[g(\bx_1)+\frac{(1-\tau)\blambda_1}{\beta_1}\right]_+\right\|^2+(1-\tau)\left(\frac{1}{2\beta_{1}}+\frac{2\vartheta}{\beta_{0}\tau}\right)\|\blambda_{1}-\blambda_{0}\|^2
\\
+\left(\beta_{0}L^2_g+\frac{4\vartheta L^2_{g}}{\gamma_{0}(1-\tau)}\right)\|\bx_{1}-\bx_{0}\|^2-f^{\star}
\end{multline}
and use the fact that
\begin{equation}\label{eq.stepo2}
\frac{2}{\beta_{r}}\left(\frac{\gamma_{r-1}}{\gamma_{r}}-1\right)=\frac{2\gamma_r}{\tau}\left(\frac{\gamma_{r-1}}{\gamma_{r}}-1\right) = \frac{2}{\tau}\left(\gamma_{r-1} - \gamma_r\right)=\frac{2}{\beta_{r-1}}-\frac{2}{\beta_r},
\end{equation}
so that we can have
\begin{align}\notag
&\quad\;\frac{(1-\tau)^2}{\beta_{R+1}}\|\blambda_{R+1}\|^2
+\frac{(1-\tau)\gamma_{R}}{2}\|\blambda_{R+1}\|^2+\frac{2\vartheta}{\beta_{R}}\left(\frac{\gamma_{R-1}}{\gamma_{R}}-1\right)\|\blambda_{R+1}\|^2
\\
&\le \frac{(1-\tau)^2}{\beta_{R+1}}\|\blambda_{R+1}\|^2
+\frac{(1-\tau)\tau}{2\beta_R}\|\blambda_{R+1}\|^2+\left(\frac{2\vartheta}{\beta_{R-1}}-\frac{2\vartheta}{\beta_R}\right)\|\blambda_{R+1}\|^2
\\
&\le   \left(\frac{1-\frac{3}{2}\tau+\frac{\tau^2}{2}}{\beta_{R}}+\frac{2\vartheta}{\beta_{R-1}}-\frac{2\vartheta}{\beta_R}\right)\|\blambda_{R+1}\|^2\label{eq.pte0}
\\
&\le    \frac{2}{\beta_{R-1}}\|\blambda_{R+1}\|^2\label{eq.pte}
\\
& \mathop{\le}\limits^{\eqref{eq.upbl}}   \frac{2G \left(1+\frac{1}{\theta'}\right)R^{2/3}}{(1-\eta)\beta_{R-1}}\le\frac{4G \left(1+\frac{1}{\theta'}\right)R^{1/3}}{(1-\eta)\beta_0}.
\end{align}

Then, we can obtain
\begin{align}\notag
&\quad\;\sigma^2\sum^R_{r=1}\alpha^2_{r+1}\|\blambda_{r+1}^g\|^2
\\\notag
&\mathop{\le}\limits^{(a)}  R^{1/3}\left( \frac{16M^2}{\beta^2_0U_JL_g}+ 4\max\left\{36,\frac{64U_J}{L_g}\right\}\alpha_1\frac{4G \left(1+\frac{1}{\theta'}\right)}{(1-\eta)\beta_0}\right)
+4\max\left\{36,\frac{64U_J}{L_g}\right\}\alpha_1\Delta_F
\\\notag
&\quad+\max\left\{36,\frac{64U_J}{L_g}\right\}\sum^R_{r=1}\frac{2.5}{\beta_0r^{4/3}}\alpha_r\|\blambda_{r+1}\|^2+\max\left\{36,\frac{64U_J}{L_g}\right\}\sum^R_{r=1}\frac{13\vartheta}{\beta_0r^{4/3}}\alpha_r\|\blambda_{r}\|^2
\\
&\quad+64(1-\tau)^2U^2_J\sum^R_{r=1}\alpha^2_r\|\blambda_r\|^2
\\\notag
& \mathop{\le}\limits^{(b)}  R^{1/3}\left( \frac{16M^2}{\beta^2_0U_JL_g}+ 4\max\left\{36,\frac{64U_J}{L_g}\right\}\alpha_1\frac{4G \left(1+\frac{1}{\theta'}\right)}{(1-\eta)\beta_0}\right)
+4\max\left\{36,\frac{64U_J}{L_g}\right\}\alpha_1\Delta_F
\\\notag
&\quad+\max\left\{36,\frac{64U_J}{L_g}\right\}\sum^R_{r=1}\frac{5}{\beta_0r^{4/3}}\alpha_r\|\blambda_{r}\|^2+\max\left\{36,\frac{64U_J}{L_g}\right\}\sum^R_{r=1}\frac{13\vartheta}{\beta_0r^{4/3}}\alpha_r\|\blambda_{r}\|^2
\\
&\quad+64(1-\tau)^2U^2_J\sum^R_{r=1}\alpha^2_r\|\blambda_r\|^2+\max\left\{36,\frac{64U_J}{L_g}\right\}\frac{5}{\beta_0R^{4/3}}\alpha_R\|\blambda_{R+1}\|^2
\\\notag
& \mathop{\le}\limits^{(c)}  R^{1/3}\left( \frac{16M^2}{\beta^2_0U_JL_g}+ 4\max\left\{36,\frac{64U_J}{L_g}\right\}\alpha_1\frac{4G \left(1+\frac{1}{\theta'}\right)}{(1-\eta)\beta_0}\right)
+4\max\left\{36,\frac{64U_J}{L_g}\right\}\alpha_1\Delta_F
\\
&\quad+66(1-\tau)^2U^2_J\sum^R_{r=1}\alpha^2_r\|\blambda_r\|^2+\max\left\{36,\frac{64U_J}{L_g}\right\}\frac{5}{\beta_0R^{4/3}}\alpha_R\|\blambda_{R+1}\|^2
\label{eq.bdd}
\end{align}
where in $(a)$ we use \eqref{eq.lup} and
\begin{equation}
16M^2\sum^R_{r=1}\alpha^2_r\mathop{\le}\limits^{\eqref{eq.conal}} 16M^2\sum^R_{r=1}\frac{1}{\beta^2_rU_JL_g}\le\frac{16M^2R^{1/3}}{\beta^2_0U_JL_g},
\end{equation}
$(b)$ follows as
\begin{equation}
\alpha_r=\frac{\alpha_r}{\alpha_{r+1}}\alpha_{r+1}\le\frac{(r+1)^{1/3}}{r^{1/3}}\alpha_{r+1}\le 2^{1/3}\alpha_{r+1}< 2\alpha_{r+1},\forall r\ge1,
\end{equation}
and $(c)$ is true when
\begin{equation}
\alpha_r\ge\max\left\{36,\frac{64U_J}{L_g}\right\}\frac{13\vartheta}{(1-\tau)^2U^2_J}\frac{1}{\beta_0r^{4/3}},\vartheta\ge 1, 
\end{equation}
i.e.,
\begin{center}
\fbox{
\begin{minipage}{0.8\columnwidth}
\begin{equation}\label{eq.rea1}
\alpha_r\ge\max\left\{36,\frac{64U_J}{L_g}\right\} \frac{13\vartheta}{(1-\tau)^2U^2_J} \frac{1}{\beta_0(r-1)^{4/3}},\forall r\ge 2,
\end{equation}
\end{minipage}
}
\end{center}
then, we can have
\begin{equation}
\max\left\{36,\frac{64U_J}{L_g}\right\}\frac{\max\{13\vartheta,5\}}{\beta_0r^{4/3}}\le(1-\tau)^2U^2_J\alpha_r.
\end{equation}

For the case  when $[g(\bx_{r+1})]_i\le 0$, we have $[\blambda_{r+1}]_i\le(1-\tau)[\blambda_r]_i$. So, we have $\alpha^2_{r+1}[\blambda_{r+1}]^2_i<\alpha^2_r[\blambda_{r+1}]^2_i\le(1-\tau)^2
\alpha^2_r[\blambda_r]^2_i$, i.e.,
\begin{equation}\label{eq.sofs}
\sigma^2\alpha^2_{r+1}[\blambda_{r+1}]^2_i < \sigma^2(1-\tau)^2\alpha^2_r[\blambda_r]^2_i,\quad \textrm{if} \quad g_i(\bx_{r+1})\le 0.
\end{equation}

Then, from \eqref{eq.bdd}, we can obtain
\begin{align}\notag
 \sigma^2 \sum^R_{r=1}\alpha^2_{r+1}\|\blambda_{r+1}\|^2
\le & R^{1/3}\left( \frac{16M^2}{\beta^2_0U_JL_g}+ 4\max\left\{36,\frac{64U_J}{L_g}\right\}\alpha_1\frac{4G \left(1+\frac{1}{\theta'}\right)}{(1-\eta)\beta_0}\right)
+4\max\left\{36,\frac{64U_J}{L_g}\right\}\alpha_1\Delta_F
\\
&+(66U^2_J+\sigma^2)(1-\tau)^2\sum^R_{r=1}\alpha^2_r\|\blambda_r\|^2+\max\left\{36,\frac{64U_J}{L_g}\right\}\frac{5}{\beta_0R^{4/3}}\alpha_R\|\blambda_{R+1}\|^2
.\label{eq.upl}
\end{align}
Dividing $\sigma^2$ on both sides of \eqref{eq.upl} results in the following recursion
\begin{align}\notag
\sum^R_{r=1}\alpha^2_{r+1}\|\blambda_{r+1}\|^2
\le & \frac{R^{1/3}}{\sigma^2}\left( \frac{16M^2}{\beta^2_0U_JL_g}+ 4\max\left\{36,\frac{64U_J}{L_g}\right\}\alpha_1\frac{4G \left(1+\frac{1}{\theta'}\right)}{(1-\eta)\beta_0}\right)
+\frac{4}{\sigma^2}\max\left\{36,\frac{64U_J}{L_g}\right\}\alpha_1\Delta_F
\\
&+\frac{5}{\sigma^2}\max\left\{36,\frac{64U_J}{L_g}\right\}\frac{\alpha_0}{\beta_0 R}\frac{G \left(1+\frac{1}{\theta'}\right)}{1-\eta}+\frac{(66U^2_J+\sigma^2)(1-\tau)^2}{\sigma^2}\sum^R_{r=1}\alpha^2_r\|\blambda_r\|^2,\label{eq.upl2}
\end{align}
which follows from \eqref{eq.upbl}.

Then, when
\begin{equation}
\frac{(66U^2_J+\sigma^2)(1-\tau)^2}{\sigma^2}\bydef \rho <1,
\end{equation}
i.e.,
\begin{center}
\fbox{
\begin{minipage}{0.4\columnwidth}
\begin{equation}\label{eq.retau}
\tau>1-\frac{\sigma}{\sqrt{66U^2_J+\sigma^2}},
\end{equation}
\end{minipage}
}
\end{center}
then, we have
\begin{multline}\notag
\sum^R_{r=1}\alpha^2_{r+1}\|\blambda_{r+1}\|^2
\le \rho \Lambda\beta_0R^{1/3}
+\frac{4\alpha_1\Delta_F}{\sigma^2}\max\left\{36,\frac{64U_J}{L_g}\right\}+\frac{5}{\sigma^2}\max\left\{36,\frac{64U_J}{L_g}\right\}\frac{\alpha_0}{\beta_0 R}\frac{G \left(1+\frac{1}{\theta'}\right)}{1-\eta}
\\
+\frac{R^{1/3}}{\sigma^2}\left( \frac{16M^2}{\beta^2_0U_JL_g}+ 4\alpha_1\max\left\{36,\frac{64U_J}{L_g}\right\}\frac{4G \left(1+\frac{1}{\theta'}\right)}{(1-\eta)\beta_0}\right).
\end{multline}
Therefore,  when
\begin{equation}
\Lambda\ge\frac{1}{\sigma^2(1-\rho)\beta_0}\left( \frac{16M^2}{\beta^2_0U_JL_g}+ 4\alpha_1\max\left\{36,\frac{64U_J}{L_g}\right\}\frac{4G \left(1+\frac{1}{\theta'}\right)}{(1-\eta)\beta_0}\right)
+\frac{4\alpha_1\Delta_F+5\frac{\alpha_0}{\beta_0 R}\frac{G \left(1+\frac{1}{\theta'}\right)}{1-\eta}}{R^{1/3}\sigma^2(1-\rho)\beta_0}\max\left\{36,\frac{64U_J}{L_g}\right\},
\end{equation}
we can obtain
\begin{equation}
\sum^R_{r=1}\alpha^2_{r+1}\|\blambda_{r+1}\|^2\le\Lambda\beta_{R},
\end{equation}
which completes the proof of showing \eqref{eq.inl2}.
\end{proof}

\section{Convergence Rate of GDPA to KKT Points}

\begin{proof}
\noindent{\bf Stationarity}.

First, multiplying  $\beta_{r}$ on both sides of \eqref{eq.updiffbeta}, we have
\begin{align}
\|\blambda_{r+1}-\blambda_r\|^2+\frac{\beta_{r}}{\alpha_r}\|\bx_{r+1}-\bx_r\|^2
\le4\beta_r(\mathcal{P}_r-\mathcal{P}_{r+1})+\frac{2.5}{r}\|\blambda_{r+1}\|^2+\frac{13\vartheta}{r}\|\blambda_{r}\|^2.\label{eq.ccov2}
\end{align}
From \eqref{eq.conal}, we choose step-size $\alpha_r$ as $\alpha_{0, 1}/(\alpha_{0,2}+\alpha_{0,3}r^{1/3} )$. According to \eqref{eq.upbl}, we can see that if $\alpha_{0,3}\ge(1-\tau)\sqrt{ \frac{G \left(1+\frac{1}{\theta'}\right)}{1-\eta}}L_J+\beta_0U_JL_g$, then  it is easy to obtain a constant as the lower bound of $\alpha_r\beta_r$, which satisfies \eqref{eq.conal}. For example, when $\alpha_{0,2}\ge L_f$ and $\alpha_{0,1}<1$, then $\alpha_r\ge1/(\alpha_{0,2}+\alpha_{0,3}r^{1/3})$ and we can get
\begin{center}
\fbox{
\begin{minipage}{0.68\columnwidth}
\begin{equation}
\alpha_r\beta_r=\frac{\alpha_{0,1}\beta_0r^{1/3}}{\alpha_{0,2}+\alpha_{0,3}r^{1/3}}\ge \frac{\alpha_{0,1}\beta_0}{\alpha_{0,2}+\alpha_{0,3}}\bydef I_1.\label{eq.ab}
\end{equation}
\end{minipage}
}
\end{center}

So, we can have equation \eqref{eq.ccov2} re-written as
\begin{align}\notag
&\min\left\{1,I_1\right\}\left(\|\blambda_{r+1}-\blambda_r\|^2+\frac{1}{\alpha^2_r}\|\bx_{r+1}-\bx_r\|^2\right)
\\\notag
\le&\|\blambda_{r+1}-\blambda_r\|^2+\frac{I_1}{\alpha^2_r}\|\bx_{r+1}-\bx_r\|^2
\\
\le&4\beta_r(\mathcal{P}_r-\mathcal{P}_{r+1})+\frac{2.5}{r}\|\blambda_{r+1}\|^2+\frac{13\vartheta}{r}\|\blambda_{r}\|^2.\label{eq.ccov3}
\end{align}

From \eqref{eq.ccov1}, we know
\begin{align}\notag
&\quad\;\|\mathcal{G}(\bx_r, \blambda_r)\|^2
\\
&\mathop{\le}\limits^{(a)}   4\left(2\left(\frac{3}{\alpha_r}\right)^2+8L^2_g+2U^2_JL^2_g\beta^2_r\right)\|\bx_{r+1}-\bx_r\|^2
+4\left(\left(\frac{3+2\tau}{\beta_0}\right)^2+3\right)\|\blambda_{r+1}-\blambda_r\|^2+16\gamma_r^2\|\blambda_{r+1}\|^2
\\
&\mathop{\le}\limits^{\eqref{eq.conal}} \frac{112}{\alpha^2_r}\|\bx_{r+1}-\bx_r\|^2
+4\left(\frac{25}{\beta^2_0}+3\right)\|\blambda_{r+1}-\blambda_r\|^2+16\gamma_r^2\|\blambda_{r+1}\|^2
\\
&\le \max\left\{112,4\left(\frac{25}{\beta^2_0}+3\right)\right\}\left(\|\blambda_{r+1}-\blambda_r\|^2+\frac{1}{\alpha^2_r}\|\bx_{r+1}-\bx_r\|^2\right)+16\gamma_r^2\|\blambda_{r+1}\|^2 \label{eq.ccov4}
\end{align}
where $(a)$ is true due to the increase of sequence $\beta_r$, and we require $\beta_0\ge1/U_J$.

Combining \eqref{eq.ccov2} and \eqref{eq.ccov4},  we know that there exists a constant $C$ such that
\begin{align}
\|\mathcal{G}(\bx_r, \blambda_r)\|^2
\le  \underbrace{\frac{\max\left\{112,4\left(\frac{25}{\beta^2_0}+3\right)\right\}}{\min\left\{1,I_1\right\}}}_{\bydef C}\left(4\beta_r(\mathcal{P}_r-\mathcal{P}_{r+1})+\frac{2.5}{r}\|\blambda_{r+1}\|^2+\frac{13\vartheta}{r}\|\blambda_{r}\|^2\right)+16\gamma_r^2\|\blambda_{r+1}\|^2 .\label{eq.ccov5}
\end{align}

Applying the telescoping sum, we have that
\begin{align}\notag
&\quad\;\sum^T_{r=1}\frac{1}{\beta_r}\|\mathcal{G}(\bx_r, \blambda_r)\|^2
\\\notag
& \mathop{\le}\limits^{\eqref{eq.ccov5} } 4C \sum^T_{r=1} \left(\mathcal{P}_r-\mathcal{P}_{r+1}\right)+\sum^T_{r=1}\frac{2.5C}{r\beta_r}\|\blambda_{r+1}\|^2
+\sum^T_{r=1}\frac{13C\vartheta}{r\beta_r}\|\blambda_{r}\|^2+16\sum^T_{r=1}\gamma_r^2\|\blambda_{r+1}\|^2
\\\notag
&\mathop{\le}\limits^{\eqref{eq.dep}\eqref{eq.sumofp}\eqref{eq.pte0}} 4C\left(\Delta_F-(1-\tau)\left(\frac{1}{2\beta_{T}}+\frac{2\vartheta}{\beta_{T}\tau}\right)\|\blambda_{T+1}-\blambda_{T}\|^2+\left(\frac{1-\frac{3}{2}\tau+\frac{\tau^2}{2}}{\beta_{T}}+\frac{2}{\beta_{T-1}}-\frac{2}{\beta_T}\right)\|\blambda_{T+1}\|^2\right)
\\\notag
&\quad+\sum^T_{r=1}\frac{2.5C}{r\beta_r}\|\blambda_{r+1}\|^2
+\sum^T_{r=1}\frac{13C\vartheta}{r\beta_r}\|\blambda_{r}\|^2+16\sum^T_{r=1}\frac{\gamma_r^2}{\beta_r}\|\blambda_{r+1}\|^2
\\\notag
& \mathop{\le}\limits^{(a)} 4C\left(\Delta_F-(1-\tau)\left(\frac{1}{2\beta_{T}}+\frac{2\vartheta}{\beta_{T}\tau}\right)\|\blambda_{T+1}-\blambda_{T}\|^2+\left(\frac{1-\frac{3}{2}\tau+\frac{\tau^2}{2}}{\beta_{T}}\right)\|\blambda_{T+1}\|^2+\frac{2}{\beta_0T^{4/3}}\|\blambda_{T+1}\|^2\right)
\\\notag
&\quad+\sum^T_{r=1}\frac{2.5C}{r\beta_r}\|\blambda_{r+1}\|^2
+\sum^T_{r=1}\frac{13C\vartheta}{r\beta_r}\|\blambda_{r}\|^2+16\sum^T_{r=1}\frac{\tau}{\beta^3_r}\|\blambda_{r+1}\|^2
\\\notag
& \mathop{\le}\limits^{(b)} 4C\left(\Delta_F-\frac{2\vartheta}{\beta_{T}\tau}\|\blambda''_{T+1}-\blambda'_{T}\|^2+\left(\frac{1-\frac{3}{2}\tau+\frac{\tau^2}{2}}{\beta_{T}}\right)\|\blambda''_{T+1}\|^2+\frac{2G(1+1/\theta')}{\beta_0T^{2/3}(1-\eta)}\right)
\\\notag
&\quad+\sum^T_{r=1}\frac{2.5C+4\tau}{\beta^3_r}\|\blambda_{r}\|^2+\sum^T_{r=1}\frac{13C\vartheta}{\beta^3_r}\|\blambda_{r}\|^2+\frac{2.5C+16\tau}{\beta^3_T}\|\blambda_{T+1}\|^2
\end{align}
where in $(a)$ we use \eqref{eq.ssr2} so that $1/\beta_{T-1}-1/\beta_T\le1/(T-1)^{4/3}\le 1/(\beta_0T^{4/3}),\forall T\ge2$,
in $(b)$ we use
\begin{equation}
 [\blambda''_{r+1}]_i \bydef \begin{cases} [\blambda_{r+1}]_i, & \textrm{if}\quad i\in\mathcal{S}_r; \\ 0, & \textrm{otherwise}.\end{cases},
\end{equation}
 \eqref{eq.upbl} and \eqref{eq.shk2} as well as the facts that $0<\tau<1$, $\beta_r=\beta_0r^{1/3}$, and $\beta_0$ is a constant.

To further get an upper bound of $\sum^T_{r=1}\frac{1}{\beta_r}\|\mathcal{G}(\bx_r, \blambda_r)\|^2$, we notice that
\begin{equation}
\left(\sum^T_{r=1}\frac{\alpha_r}{\beta^2_r}\|\blambda_{r+1}\|\right)^2\le\left(\sum^T_{r=1}\frac{1}{\beta^4_r}\right)\left(\sum^T_{r=1}\alpha^2_r\|\blambda_{r+1}\|^2\right)\mathop{\le}\limits^{\eqref{eq.inl2}} \frac{1}{\beta^4_0}T^{-1/3}\Lambda\beta_0T^{1/3}=\frac{\Lambda}{\beta^3_0}
\end{equation}
by applying the Cauchy-Schwarz inequality. Then, we can obtain
\begin{align}
\sum^T_{r=1}\frac{1}{\beta_r}\|\mathcal{G}(\bx_r, \blambda_r)\|^2\mathop{\le}\limits^{\eqref{eq.ab}} &4C\left(\Delta_F-\frac{2\vartheta}{\beta_{T}\tau}\|\blambda''_{T+1}-\blambda'_{T}\|^2+\left(\frac{1-\frac{3}{2}\tau+\frac{\tau^2}{2}}{\beta_{T}}\right)\|\blambda''_{T+1}\|^2+\frac{2G(1+1/\theta')}{\beta_0T^{2/3}(1-\eta)}\right)
\\\notag
&+\frac{2.5C+16\tau+13C\vartheta}{I_1}\frac{\Lambda^{1/2}}{\beta^{3/2}_0}+\frac{(2.5C+4\tau)G(1+1/\theta')}{\beta^3_0T^{1/3}(1-\eta)}.
\end{align}

For the case where $[\blambda_{T+1}]_i>0$, i.e., $i\in\mathcal{S}_r$,  we have
\begin{equation}
\|\blambda''_{T+1}\|^2\mathop{\le}\limits^{(a)} \eta\|\blambda'_T\|^2+\left(1+\frac{1}{\theta'}\right)\beta^2_T\|g_+(\bx_{T+1})\|^2
\end{equation}
where $(a)$ follows from \eqref{eq.reub}, \eqref{eq.shk}, and \eqref{eq.shk2}.

Therefore, we have
\begin{align}\notag
\sum^T_{r=1}\frac{1}{\beta_r}\|\mathcal{G}(\bx_r, \blambda_r)\|^2\le &-\frac{2\vartheta}{\beta_{T}\tau}\|\blambda''_{T+1}-\blambda'_{T}\|^2+\left(\frac{1-\frac{3}{2}\tau+\frac{\tau^2}{2}}{\beta_{T}}\right)\left(\eta\|\blambda'_T\|^2+\left(1+\frac{1}{\theta'}\right)\beta^2_T\|g_+(\bx_{T+1})\|^2\right)
\\
&+4C\Delta_F+\frac{2.5C+16\tau+13C\vartheta}{I_1}\frac{\Lambda^{1/2}}{\beta^{3/2}_0}+\frac{(2.5C+4\tau)G(1+1/\theta')}{\beta^3_0T^{1/3}(1-\eta)}+\frac{8CG(1+1/\theta')}{\beta_0T^{2/3}(1-\eta)}.\label{eq.opte}
\end{align}

Before showing the convergence of the optimality gap, we first verify the satisfaction of the functional constraints as follows.

\noindent{\bf Stopping Criterion (Constraints)}

According to \eqref{eq.dualupated} and \eqref{eq.deffgp}, we know
\begin{align}\notag
&\quad\;\|\blambda^g_{r+1}-\blambda'^g_r\|^2
\\
&= \| \beta_rg'(\bx_{r+1})-\tau\blambda'^g_r \|^2
\\
&= \sum_{i\in\mathcal{A}_{r+1}\cap\mathcal{S}_r}\|\beta_rg_i(\bx_{r+1})-\tau[\blambda_r]_i\|^2
\\
&= \sum_{i\in\mathcal{A}_{r+1}}\|\beta_rg_i(\bx_{r+1})-\tau[\blambda_r]_i\|^2-\sum_{i\in\mathcal{A}_{r+1}\cap\overbar{\mathcal{S}}_r}\|\beta_rg_i(\bx_{r+1})-\tau[\blambda_r]_i\|^2,
\end{align}
which gives
\begin{align}
\sum_{i\in\mathcal{A}_{r+1}}\|\beta_rg_i(\bx_{r+1})-\tau[\blambda_r]_i\|^2&\le\|\blambda^g_{r+1}-\blambda'^g_r\|^2 +\sum_{i\in\mathcal{A}_{r+1}\cap\overbar{\mathcal{S}}_r}\|\beta_rg_i(\bx_{r+1})-\tau[\blambda_r]_i\|^2
\\
&\mathop{\le}\limits^{(a)}   \|\blambda_{r+1}-\blambda_r\|^2+\sum_{i\in\mathcal{A}_{r+1}\cap\overbar{\mathcal{S}}_r}2\left\|\beta_r\left(g_i(\bx_{r+1})-g_i(\bx_r)\right)\right\|^2+2\tau^2\|\blambda_r\|^2
\\
&\mathop{\le}\limits^{(b)}   \|\blambda_{r+1}-\blambda_r\|^2+2L^2_g\beta^2_r\|\bx_{r+1}-\bx_r\|^2+2\tau^2\|\blambda_r\|^2\label{eq.recs}
\end{align}
where $(a)$ is true because of $g_i(\bx_r)\le0, g_i(\bx_{r+1})>0,\forall i\in\mathcal{A}_{r+1}\cap\overbar{\mathcal{S}}_r$, in $(b)$ we apply Lipschitz continuity.

Note that
\begin{equation}
\|\beta_r g_+(\bx_{r+1})\|^2-2\beta_r\tau\langle g_+(\bx_{r+1}),\blambda_r\rangle+\tau^2\|\blambda_r\|^2=\| \beta_r g_+(\bx_{r+1})-\tau\blambda_r \|^2\le\sum_{i\in\mathcal{A}_{r+1}}\|\beta_rg_i(\bx_{r+1})-\tau[\blambda_r]_i\|^2.\label{eq.bet}
\end{equation}

Combining \eqref{eq.recs} and \eqref{eq.bet}, we can have
\begin{align}\notag
&\quad\;\|\beta_r g_+(\bx_{r+1})\|^2+\tau^2\|\blambda_r\|^2
\\
&\le   2\beta_r\tau\langle g_+(\bx_{r+1}),\blambda_r\rangle+\|\blambda_{r+1}-\blambda_r\|^2+2L^2_g\|\bx_{r+1}-\bx_r\|^2+2\tau^2\|\blambda_r\|^2
\\
&\le \frac{1}{2}\|\beta_r g_+(\bx_{r+1})\|^2+2\tau^2\|\blambda_r\|^2+\|\blambda_{r+1}-\blambda_r\|^2+2L^2_g\beta^2_r\|\bx_{r+1}-\bx_r\|^2+2\tau^2\|\blambda_r\|^2
\end{align}
where we apply the Young's inequality with parameter $2$ to $\beta_r\tau\langle g_+(\bx_{r+1}),\blambda_r\rangle$.

Then, by algebraic manipulation, we can obtain
\begin{equation}
\|\beta_r g_+(\bx_{r+1})\|^2\le 2\|\blambda_{r+1}-\blambda_r\|^2+4L^2_g\beta^2_r\|\bx_{r+1}-\bx_r\|^2+6\tau^2\|\blambda_r\|^2.
\end{equation}
Dividing $\beta_r$ on both sides of the above inequality gives
\begin{align}\notag
&\quad\;\beta_r\| g_+(\bx_{r+1})\|^2
\\
&=\frac{1}{\beta_r}\|\beta_r g_+(\bx_{r+1})\|^2
\\
&\mathop{\le}\limits^{(a)}   \frac{2}{\beta_r}\|\blambda_{r+1}-\blambda_r\|^2+ 4L^2_g\beta_r\|\bx_{r+1}-\bx_r\|^2+\frac{6\tau^2}{\beta_r}\|\blambda_r\|^2\label{eq.stepb}
\\
&\mathop{\le}\limits^{\eqref{eq.conal}}  \frac{2}{\beta_r}\|\blambda_{r+1}-\blambda_r\|^2+ \frac{4L_g}{U_J\alpha_r}\|\bx_{r+1}-\bx_r\|^2+\frac{6\tau^2}{\beta_r}\|\blambda_r\|^2
\\
&\mathop{\le}\limits^{\eqref{eq.updiffbeta2}} 8\underbrace{\max\left\{1,\frac{2L_g}{U_J}\right\}}_{\bydef C_1}(\mathcal{P}_r-\mathcal{P}_{r+1})+\frac{5\max\left\{1,\frac{2L_g}{U_J}\right\}}{\beta_rr}\|\blambda_{r+1}\|^2+\frac{26\vartheta\max\left\{1,\frac{2L_g}{U_J}\right\}}{\beta_rr}\|\blambda_{r}\|^2+\frac{6\tau^2\|\blambda_r\|^2}{\beta_r}
\end{align}
where in $(a)$ we use \eqref{eq.conal} and choose $\beta_0\ge\sqrt{2L_g/U_J}$ so that $2L^2_g/\beta_r\le\frac{1}{\alpha_r}$.

Applying the telescoping sum, we can have
\begin{align}\notag
&\quad\;\sum^T_{r=1}\beta_r\| g_+(\bx_{r+1})\|^2
\\
&\mathop{\le}\limits^{\eqref{eq.sumofp}\eqref{eq.pte}} 8C_1\left(\Delta_F+\frac{2}{\beta_{T-1}}\|\blambda_{T+1}\|^2\right)+\sum^T_{r=1}\left(\frac{5C_1}{\beta_r}\|\blambda_{r+1}\|^2+\frac{26C_1\vartheta+6\tau^2}{\beta_r}\|\blambda_{r}\|^2\right)
\\
&\le  8C_1\left(\Delta_F+\frac{2}{\beta_{T-1}}\|\blambda_{T+1}\|^2+\frac{5C_1}{\beta_T}\|\blambda_{T+1}\|^2\right)+\sum^T_{r=1}\frac{(5+26\vartheta)C_1+6\tau^2}{\beta_r}\|\blambda_{r}\|^2
\\
&\mathop{\le}\limits^{(a)}   8C_1\left(\Delta_F+\frac{2}{\beta_{T-1}}\|\blambda_{T+1}\|^2+\frac{5C_1}{\beta_T}\|\blambda_{T+1}\|^2\right)+\frac{\Lambda\beta_0((5+26\vartheta)C_1+6\tau^2)}{\alpha_{0,1}I_1}\left(\alpha_{0,2}T^{1/3}+\alpha_{0,3}T^{2/3}\right)
\\\notag
&\mathop{\le}\limits^{(b)}  8C_1\left(\Delta_F+\frac{2}{\beta_{T-1}}\frac{G \left(1+\frac{1}{\theta'}\right)T^{2/3}}{1-\eta}+\frac{5C_1}{\beta_{T}}\frac{G \left(1+\frac{1}{\theta'}\right)T^{2/3}}{1-\eta}\right)
\\
&\quad +\frac{\Lambda\beta_0((5+26\vartheta)C_1+6\tau^2)}{\alpha_{0,1}I_1}\left(\alpha_{0,2}T^{1/3}+\alpha_{0,3}T^{2/3}\right)
\\
&\le   8C_1\left(\Delta_F+\frac{(4+5C_1)G \left(1+\frac{1}{\theta'}\right)T^{1/3}}{\beta_0(1-\eta)}\right)+\frac{\Lambda\beta_0((5+26\vartheta)C_1+6\tau^2)}{\alpha_{0,1}I_1}\left(\alpha_{0,2}T^{1/3}+\alpha_{0,3}T^{2/3}\right)
\end{align}
where in $(a)$ we use
\begin{align}
\sum^T_{r=1}\frac{1}{\beta_r}\|\blambda_{r}\|^2&\mathop{\le}\limits^{\eqref{eq.ab}}\sum^T_{r=1}\frac{\alpha_r}{I_1}\|\blambda_{r}\|^2\le \frac{1}{\alpha_{T}I_1}\sum^T_{r=1}\alpha^2_{r}\|\blambda_{r}\|^2
\\
&\mathop{\le}\limits^{\eqref{eq.inl} }   \Lambda\beta_0T^{1/3}\left(\frac{\alpha_{0,2}+\alpha_{0,3}T^{1/3}}{\alpha_{0,1}I_1}\right)\le\frac{\Lambda\beta_0}{\alpha_{0,1}I_1}\left(\alpha_{0,2}T^{1/3}+\alpha_{0,3}T^{2/3}\right)\sim\mathcal{O}(T^{2/3}),
\end{align}
and in $(b)$ we use \eqref{eq.upbl}.

Therefore, we have
\begin{align}\notag
\left(\sum^T_{r=1}\beta_r\right)^{-1}\sum^T_{r=1}\beta_r\| g_+(\bx_{r+1})\|^2\le &
\frac{8C_1}{\beta_0T^{4/3}}\left(\Delta_F+\frac{(4+5C_1)G \left(1+\frac{1}{\theta'}\right)T^{1/3}}{\beta_0(1-\eta)}\right)
\\
&+\frac{\Lambda\beta_0((5+26\vartheta)C_1+6\tau^2)}{\alpha_{0,1}I_1\beta_0T^{4/3}}\left(\alpha_{0,2}T^{1/3}+\alpha_{0,3}T^{2/3}\right). \label{eq.conc}
\end{align}
It is implied that combining \eqref{eq.conc} and \eqref{eq.spcon} gives $\lim_{T\to\infty}\| g_+(\bx_{T})\|^2=0$, meaning that the constraints are satisfiable.  Based on the definition of $T(\epsilon)$ (we will use $T$ as a shortcut of $T(\epsilon)$ in the following) in \eqref{eq.gapc}, we have
\begin{align}
\| g_+(\bx_{T+1})\|^2&\le\frac{8C_1\left(\Delta_F+\frac{(4+5C_1)G \left(1+\frac{1}{\theta'}\right)T^{1/3}}{\beta_0(1-\eta)}\right)+\frac{\Lambda\beta_0((5+26\vartheta)C_1+6\tau^2)}{\alpha_{0,1}I_1}\left(\alpha_{0,2}T^{1/3}+\alpha_{0,3}T^{2/3}\right)}{\sum^T_{r=1}\beta_r}
\\
&\le\frac{8C_1\Delta_F}{\beta_0T^{4/3}}+\frac{\frac{8C_1(4+5C_1)G \left(1+\frac{1}{\theta'}\right)}{\beta_0(1-\eta)}+\frac{\Lambda\beta_0((5+26\vartheta)C_1+6\tau^2)}{\alpha_{0,1}I_1}\alpha_{0,2}}{\beta_0T}+\frac{\Lambda ((5+26\vartheta)C_1+6\tau^2)\alpha_{0,3}}{\alpha_{0,1}I_1 T^{2/3}}.
\end{align}

For convenience of expression, let
\begin{equation}
C_{1,1}\bydef \frac{8C_1\Delta_F}{\beta_0},\quad C_{1,2}\bydef \frac{8C_1(4+5C_1)G \left(1+\frac{1}{\theta'}\right)}{\beta^2_0(1-\eta)}+\frac{\Lambda\alpha_{0,2}((5+26\vartheta)C_1+6\tau^2)}{\alpha_{0,1}I_1},\quad C_{1,3}\bydef\frac{\Lambda((5+26\vartheta)C_1+6\tau^2)\alpha_{0,3}}{\alpha_{0,1}}.
\end{equation}
Then, when $T\ge\max\{C^{3/2}_{1,1},C^{3}_{1,2}\}$, then we have
\begin{equation}
\| g_+(\bx_{T+1})\|^2\le\frac{2+C_{1,2}}{T^{2/3}}.\label{eq.uofc}
\end{equation}

Next, we will get an upper bound of $-\|\blambda''_{r+1}-\blambda'_r\|^2$ shown in \eqref{eq.opte} as follows:
\begin{align}\notag
\|\blambda''_{r+1}-\blambda'_r\|^2& =\sum_{i\in\mathcal{S}_r}([\blambda_{r+1}]_i-[\blambda_{r}]_i)^2
\\
& = \sum_{i\in\mathcal{S}_r\cap\mathcal{A}_{r+1}}([\blambda_{r+1}]_i-[\blambda_{r}]_i)^2+\sum_{i\in\mathcal{S}_r\cap\overbar{\mathcal{A}}_{r+1}}([\blambda_{r+1}]_i-[\blambda_{r}]_i)^2
\\
& \mathop{=}\limits^{\eqref{eq.shk}}   \sum_{i\in\mathcal{S}_r\cap\mathcal{A}_{r+1}}([\blambda_{r+1}]_i-[\blambda_{r}]_i)^2+\sum_{i\in\mathcal{S}_r\cap\overbar{\mathcal{A}}_{r+1}}\tau^2[\blambda_{r}]^2_i
\\
&  \mathop{=}\limits^{\eqref{eq.dualupated}}   \sum_{i\in\mathcal{S}_r\cap\mathcal{A}_{r+1}}(\beta_rg_i(\bx_{r+1})-\tau[\blambda_{r}]_i)^2+\sum_{i\in\mathcal{S}_r\cap\overbar{\mathcal{A}}_{r+1}}\tau^2[\blambda_{r}]^2_i
\\
&  \mathop{\ge}\limits^{(a)}  \sum_{i\in\mathcal{S}_r\cap\mathcal{A}_{r+1}}\left(\tau^2[\blambda_{r}]^2_i+\beta^2_rg^2_i(\bx_{r+1})\right)+\sum_{i\in\mathcal{S}_r\cap\overbar{\mathcal{A}}_{r+1}}\tau^2[\blambda_{r}]^2_i \\
&  \ge   \sum_{i\in\mathcal{S}_r\cap\mathcal{A}_{r+1}}\beta^2_rg^2_i(\bx_{r+1})+2\sum_{i\in\mathcal{S}_r}\tau^2[\blambda_{r}]^2_i\ge2\tau^2\|\blambda'_r\|^2
\end{align}
where in $(a)$ we use the fact that both $[\blambda_{r}]_i$ and $g_i(\bx_{r+1})$ are non-negative so we have
\begin{align}
(\beta_rg_i(\bx_{r+1})-\tau[\blambda_{r}]_i)^2=\tau^2[\blambda_{r}]^2_i-2\tau[\blambda_{r}]_i\beta_rg_i(\bx_{r+1})+\beta^2_rg^2_i(\bx_{r+1})\ge\tau^2[\blambda_{r}]^2_i+\beta^2_rg^2_i(\bx_{r+1}).
\end{align}

To this end, we can have the upper bound of the optimality gap further based on the satisfaction of the functional constraints, i.e.,
\begin{align}\notag
\sum^T_{r=1}\frac{1}{\beta_r}\|\mathcal{G}(\bx_r, \blambda_r)\|^2\le &-\frac{4\vartheta\tau\|\blambda'_T\|^2}{\beta_{T}}+\left(\frac{1-\frac{3}{2}\tau+\frac{\tau^2}{2}}{\beta_{T}}\right)\eta\|\blambda'_T\|^2+\left(1+\frac{1}{\theta'}\right)\beta_T\|g_+(\bx_{T+1})\|^2
\\
&+4C\Delta_F+\frac{2.5C+16\tau+13C\vartheta}{I_1}\frac{\Lambda^{1/2}}{\beta^{3/2}_0}+\frac{(2.5C+16\tau)G(1+1/\theta')}{\beta^3_0T^{1/3}(1-\eta)}+\frac{8CG(1+1/\theta')}{\beta_0T^{2/3}(1-\eta)}.\!\!
\end{align}

From \eqref{eq.uofc}, we have
\begin{equation}
\left(1+\frac{1}{\theta'}\right)\beta_T\|g_+(\bx_{T+1})\|^2\le\left(1+\frac{1}{\theta'}\right)\frac{2+C_{1,2}}{\beta_0T^{1/3}}.
\end{equation}

When $4\vartheta\tau>\eta(1-3/2\tau+\tau^2/2)$, e.g., $\vartheta=\eta(1+\tau^2/2)/(4\tau)$, $-\frac{4\vartheta\tau\|\blambda'_T\|^2}{\beta_{T}}+\left(\frac{1-\frac{3}{2}\tau+\frac{\tau^2}{2}}{\beta_{T}}\right)\eta\|\blambda'_T\|^2\le0$. Combining \leref{le.cpo}, we can choose constant $\vartheta$ as $\max\{\eta(1+\tau^2/2)/(4\tau),1\}$.

\noindent{\bf Step-sizes Selection:}
GDPA only needs to tune two parameters for ensuring convergence. Regarding the step-size of the dual update, i.e., $\beta_r$, we can choose it as $\beta_r=\beta_0r^{1/3}$ as discussed in the proof.
From \eqref{eq.reofb1}, \eqref{eq.ccov4}, \eqref{eq.stepb}, we require
\begin{equation}
\beta_0\ge\max\left\{ \frac{2(1-\tau)(\frac{\sigma}{2}+U_J)}{\sigma}\sqrt{\frac{G \left(1+\frac{1}{\theta'}\right)}{1-\eta}}, \frac{1}{U_J}, \sqrt{\frac{2L_g}{U_J}}\right\}.
\end{equation}
Then, we can consider $\beta_0$ and choose the step-size of updating primal variable, i.e., $\alpha_r$. Based on the discussion in the above proof, we choose it as $\alpha_r=\frac{\alpha_0}{\alpha_{0,2}+\alpha_{0,3}r^{1/3}}$. From \eqref{eq.conal}, \eqref{eq.stc1}, \eqref{eq.rea1} we require
\begin{equation}
\alpha_{0,3}\ge\max\left\{(1-\tau)\sqrt{ \frac{G \left(1+\frac{1}{\theta'}\right)}{1-\eta}}L_J+\beta_0U_JL_g, 2(1-\tau)L^2_g+(8\vartheta L^2_g)/\tau^2)\beta_0\right\}, \alpha_{0,2}\ge L_f,  \alpha_{0,1}<1,
 \end{equation}
then conditions \eqref{eq.conal} and \eqref{eq.stc1} can be satisfied while resulting in a consistent lower bound of $\alpha_r\beta_r$ discussed in \eqref{eq.ab}.  Note that when $r$ is larger than a constant, \eqref{eq.rea1} is satisfied automatically. To be more specific, let $C_2\bydef13\max\{36,(64U_J)/L_g\}\vartheta/((1-\tau)^2U^2_J)/\beta_0$. Then, when $r\ge(C_2(\alpha_{0,2}+\alpha_{0,3})/\alpha_{0,1})^{3/2}$, it can be easily verified that \eqref{eq.rea1}  holds.

\noindent{\bf Stationarity (for both primal and dual variables):}

Therefore, we can obtain
\begin{multline}\label{eq.sumg}
\sum^T_{r=1}\frac{1}{\beta_r}\|\mathcal{G}(\bx_r, \blambda_r)\|^2\le 4C\Delta_F+\frac{2.5C+16\tau+13C\vartheta}{I_1}\frac{\Lambda^{1/2}}{\beta^{3/2}_0}
\\
+\left(\frac{(2.5C+16\tau)G}{\beta^2_0(1-\eta)}+2+C_{1,2}\right)\frac{(1+1/\theta')}{\beta_0T^{1/3}}+\frac{8CG(1+1/\theta')}{\beta_0T^{2/3}(1-\eta)},
\end{multline}
which leads to
\begin{align}\notag
&\quad\;\|\mathcal{G}(\overbar{\bx}_T, \overbar{\blambda}_T)\|^2
\\
& \le  \left(\sum^T_{r=1}\frac{1}{\beta_r}\right)^{-1}\sum^T_{r=1}\frac{1}{\beta_r}\|\mathcal{G}(\bx_r, \blambda_r)\|^2
\\
& \le   \frac{4C\Delta_F}{\beta_0T^{2/3}}+\frac{2.5C+16\tau+13C\vartheta}{I_1}\frac{\Lambda^{1/2}}{\beta^{5/2}_0T^{2/3}}
+\left(\frac{(2.5C+16\tau)G}{\beta^2_0(1-\eta)}+2+C_{1,2}\right)\frac{(1+1/\theta')}{\beta^2_0T}+\frac{8CG(1+1/\theta')}{\beta^2_0T^{4/3}(1-\eta)}
\\
& \sim  \mathcal{O}\left(\frac{1}{T^{2/3}}\right),
\end{align}
or equivalently $\|\mathcal{G}(\overbar{\bx}_T, \overbar{\blambda}_T)\|\sim\mathcal{O}(T^{-1/3})$.

\noindent{\bf Feasibility (constraint violation):}

From \eqref{eq.conc}, we know
\begin{align}\notag
&\|\beta_T g_+(\overbar{\bx}_{T})\|^2\le\left(\sum^T_{r=1}\frac{1}{\beta_r}\right)^{-1}\sum^T_{r=1}\frac{1}{\beta_r}\|\beta_r g_+(\bx_{r})\|^2
\\
\le &
\frac{8C_1}{\beta_0T^{2/3}}\left(\Delta_F+\frac{(4+5C_1)G \left(1+\frac{1}{\theta'}\right)T^{1/3}}{\beta_0(1-\eta)}\right)
+\frac{\Lambda\beta_0((5+26\vartheta)C_1+6\tau^2)}{\alpha_{0,1}I_1\beta_0T^{2/3}}\left(\alpha_{0,2}T^{1/3}+\alpha_{0,3}T^{2/3}\right)
\\
\le & \frac{\alpha_{0,3}\Lambda\beta_0((5+26\vartheta)C_1+6\tau^2)}{\alpha_{0,1}I_1\beta_0} +  \frac{\alpha_{0,2}\Lambda((5+26\vartheta)C_1+6\tau^2)}{\alpha_{0,1}I_1T^{1/3}}
+\frac{8C_1(4+5C_1)G \left(1+\frac{1}{\theta'}\right)}{\beta^2_0(1-\eta)T^{1/3}}+\frac{8C_1\Delta_F}{\beta_0T^{2/3}},
\end{align}
which directly gives
\begin{multline}
\|g_+(\overbar{\bx}_{T})\|^2
\le\frac{\alpha_{0,3}\Lambda\beta_0((5+26\vartheta)C_1+6\tau^2)}{\alpha_{0,1}I_1\beta^3_0T^{2/3}}
 + \frac{\alpha_{0,2}\Lambda((5+26\vartheta)C_1+6\tau^2)}{\alpha_{0,1}\beta^2_0I_1T}
\\
+\frac{8C_1(4+5C_1)G \left(1+\frac{1}{\theta'}\right)}{\beta^4_0(1-\eta)T}+\frac{8C_1\Delta_F}{\beta^3_0T^{4/3}}
\sim \mathcal{O}\left(\frac{1}{T^{2/3}}\right),
\end{multline}
or equivalently $\|g_+(\overbar{\bx}_{T})\|\sim\mathcal{O}(T^{-1/3})$.

\noindent{\bf Slackness:}

Let $m$ denote the number of constraints or dimension of $\blambda_r$. Then,  when $[\blambda_{r}]_i=0$ then $[\blambda_{r}]_ig_i(\bx_{r})=0$. Otherwise, $[\blambda_{r}]_i>0$, and
\begin{align}\notag
|[\blambda_{r}]_i g_i(\bx_{r})|&\mathop{=}\limits^{\eqref{eq.dualupated}}\left|\frac{1}{\beta_r}\left(\left([\blambda_{r}]_i-[\blambda_{r-1}]_i\right)[\blambda_{r}]_i+\tau[\blambda_{r-1}]_i[\blambda_{r}]_i\right)\right|
\\
&\le \frac{1}{\beta_r}\left(\left|\left([\blambda_{r}]_i-[\blambda_{r-1}]_i\right)[\blambda_{r}]_i\right|+\frac{[\blambda_{r-1}]^2_i}{2}+\frac{[\blambda_{r}]^2_i}{2}\right)
\\
&\le \frac{1}{\beta_r}\left(\frac{([\blambda_{r}]_i-[\blambda_{r-1}]_i)^2}{2}+\frac{[\blambda_{r-1}]_i^2}{2}+[\blambda_{r}]_i^2\right),
\end{align}
so we have
\begin{equation}
\sum^m_{r=1}|[\blambda_{r}]_i g_i(\bx_{r})|\le\frac{1}{\beta_r}\left(\frac{\|\blambda_{r}-\blambda_{r-1}\|^2}{2}+\frac{\|\blambda_{r-1}\|^2}{2}+\|\blambda_r\|^2\right).
\end{equation}
Therefore, we can conclude that
\begin{align}\notag
\sum^m_{i=1}|[\overbar{\blambda}_{T}]_i g_i(\overbar{\bx}_{T})|&\le  
\left(\sum^T_{r=1}\frac{1}{\beta_r}\right)^{-1}\sum^m_{i=1}\frac{1}{\beta_r}|[\blambda_{r}]_i g_i(\bx_{r})|
\\
&\le  \left(\sum^T_{r=1}\frac{1}{\beta_r}\right)^{-1}\sum^T_{r=1}\frac{1}{\beta^2_r}\left(\frac{\|\blambda_{r}-\blambda_{r-1}\|^2}{2}+\frac{\|\blambda_{r-1}\|^2}{2}+\|\blambda_r\|^2\right)
\mathop{\le}\limits^{\eqref{eq.inl2}\eqref{eq.sumg}} \mathcal{O}\left(T^{-1/3}\right).\label{eq.slak}
\end{align}
\end{proof}

\section{Proofs for the Auxiliary Lemmas}
\subsection{Proof of \leref{le.farkas} (Approximate Farkas Lemma)}
\def\ck{\mathcal{K}}
\begin{proof}
We first need to show that the two cases cannot hold simultaneously.  Since $(\bg+\Delta)\in\mathcal{K}$, there exits vectors $\by\ge0$ and $\bw$ such that $\bg+\Delta=\bB\by+\bC\bw$. If there also existed a $\bd$ satisfying \eqref{eq.pf}, then we have by taking inner products that
\begin{align}
-\epsilon&\mathop{>}\limits^{\eqref{eq.pfa}}\bd^{\T}\bg=\bd^{\T}\bB\by+\bd^{\T}\Delta+\bd^{\T}\bC\bw=(\bB^{\T}\bd)^{\T}\by+\bd^{\T}\Delta+(\bC^{\T}\bd)^{\T}\bw
\\
 &\mathop{\ge}\limits^{(a)}\bd^{\T}\Delta\ge-\|\bd\|\|\Delta\|\mathop{\ge}\limits^{(b)}-\epsilon
\end{align}
where in $(a)$ we use $\bC^{\T}\bd=0$, $\bB^{\T}\bd\ge0$, $\by\ge0$, and in $(b)$ we use $\|\Delta\|\le\epsilon$ and $\|\bd\|=1$. Therefore, it is obvious that both cases cannot be holding at once.

Next, we will show that one of the two cases holds, i.e., we show that how to construct a vector $\bd$ such that it satisfies \eqref{eq.pf} in the case that $(\bg+\Delta)\notin\ck$. Let $\bs^{\star}$ be a vector in $\ck$ which is closest to $\bg$, i.e.,
\begin{equation}\label{eq.spro}
\bs^{\star} = \arg\min_{\bz\in\mathcal{K}}\|\bz-(\bg+\Delta)\|^2_2.
\end{equation}
So, we have
\begin{equation}\label{eq.ss}
\langle\bs^{\star}-(\bg+\Delta),\bs-\bs^{\star}\rangle\ge 0, \forall \bs\in\mathcal{K}
\end{equation}
by applying the optimality condition of \eqref{eq.spro}. Note that since $\bs\in\mathcal{K}$ and $\mathcal{K}$ is a cone, $\alpha\bs\in\mathcal{K}$ for all scalars $\alpha\ge0$. Since $\|\alpha\bs^{\star}-(\bg+\Delta)\|^2_2$ is minimized by $\alpha=1$, we have $(\bs^{\star})^{\T}(\bs^{\star}-(\bg+\Delta))=0$. Combing with \eqref{eq.ss}, we have
\begin{equation}\label{eq.opts}
\bs^{\T}(\bs^{\star}-(\bg+\Delta))\ge 0, \forall \bs\in\mathcal{K}.
\end{equation}
Define $\bd\bydef\bs^{\star}-(\bg+\Delta)$ satisfying \eqref{eq.ss}.
From \eqref{eq.opts}, we have that $\bd^{\T}\bs\ge0,\forall \bs\in\mathcal{K}$, so
\begin{equation}\label{eq.dbc}
\bd^{\T}\bB\by+\bC\bw\ge0,\quad\forall\bw,\by\ge0.
\end{equation}
In the following, we will verify that $\bd$ satisfies all the relations shown in \eqref{eq.pf}.

\begin{enumerate}
\item Verify \eqref{eq.pfa}

Note that $\bd\neq 0$ because $(\bg+\Delta)\notin\mathcal{K}$. From the fact that $(\bs^{\star})^{\T}(\bs^{\star}-(\bg+\Delta))=0$, we have
\begin{equation}
\bd^{\T}(\bg+\Delta)=\bd^{\T}(\bs^{\star}-\bd)=(\bs^{\star}-(\bg+\Delta))^{\T}\bs^{\star}-\bd^{\T}\bd=-\|\bd\|^2_2,
\end{equation}
which is
\begin{equation}
\bd^{\T}\bg=-\bd^{\T}\Delta-\|\bd\|^2_2=-(\bs^{\star}-(\bg+\Delta))^{\T}\Delta-\|\bd\|^2_2\mathop{\le}\limits^{(a)}-\|\bd\|^2_2\mathop{<}\limits^{(b)} -\epsilon
\end{equation}
where $(a)$ holds since $\Delta\in\mathcal{K}$ and \eqref{eq.opts}, and $(b)$ is true because $\|\bd\|=1$ and $0<\epsilon<1$.

\item Verify \eqref{eq.pfb} \footnote{Verifying \eqref{eq.pfb} and \eqref{eq.pfc} is identical to the original proof of the Farkas lemma (Please see Lemma 12.4 in \cite{nocedal2006numerical}) since the two relations do not include any perturbation (approximate) term.}

From \eqref{eq.dbc}, we have $(\bC^{\T}\bd)^{\T}\ge 0,\forall \bw$ when $\by=0$, which is true only if $\bC^{\T}\bd=0$.

\item Verify \eqref{eq.pfc}

Similar as the previous case, we have $(\bB\bd)^{\T}\by\ge0,\forall \by\ge0$, which is true only if $\bB^{\T}\bd\ge0$.
\end{enumerate}
In summary, we have shown that the constructed $\bd$ satisfies all the properties, which completes the proof.
\end{proof}

\subsection{Proof of \prref{pro.eq}}

Without of generality, we give the following lemma to show the relation between $\|1/\alpha [\bx - \textrm{proj}_{\mathcal{X}}(\bx-\alpha \nabla_{\bx}\mathcal{L}(\bx,\blambda))]\|\le\epsilon, \alpha>0$ and \eqref{eq.kkts1}.

\begin{mdframed}[backgroundcolor=gray!10,topline=false,
	rightline=false,
	leftline=false,
	bottomline=false]
\begin{lemma}\label{le.prok}
Let feasible set $\mathcal{X}$ be represented by differentiable continuous convex functions  $c_i(\bx),\forall i$ defined by
\begin{equation}\label{eq.defc}
c_i(\bx)=0, \;i\in\mathcal{E},\quad\textrm{and}\quad c_i(\bx)\le0, \; i\in\mathcal{I},
\end{equation}
where $\mathcal{E}$ denotes the set of the indices of equality constraints, and $\mathcal{I}$ denotes the set of the indices of inequality constraints. When $\|\mathcal{G}(\bx^{\star},\blambda^{\star})\|\le\epsilon$, i.e.,
\begin{align}
\left\|\frac{1}{\alpha }\left[\bx^{\star} - \textrm{proj}_{\mathcal{X}}(\bx^{\star}-\alpha \nabla_{\bx}\mathcal{L}(\bx^{\star},\blambda^{\star}))\right] \right\|\le  \epsilon, \quad\forall \alpha>0\label{eq.optkx}
\end{align}
then, we have
\begin{equation}\label{eq.eqe}
\textrm{dist}\left(\nabla\mathcal{L}(\bx^{\star},\blambda^{\star}),-\mathcal{N}_{\mathcal{X}}(\bx^{\star})\right)\le\epsilon'
\end{equation}
where $ \epsilon'$ has a one-to-one correspondence of $\epsilon$ and also $\epsilon'\sim\mathcal{O}(\epsilon)$.
\end{lemma}
\end{mdframed}
\begin{proof}

From Lemma 3 in \cite{lu2020finding}, we know that \eqref{eq.optkx} implies
\begin{equation}\label{eq.copt}
\langle\nabla \mathcal{L}(\bx^{\star}), \bx - \bx^{\star}\rangle\ge -\epsilon', \quad \forall \bx\in\mathcal{X}.
\end{equation}
when $\mathcal{X}$ is bounded and there is a one-to-one correspondence of $ \epsilon'$ and $\epsilon$ in the sense that $\epsilon'\sim\mathcal{O}(\epsilon)$.

 Define a cone
\begin{equation}
\mathcal{K} = \left\{\sum_{i\in\mathcal{A}'(\bx^{\star})}-\bmu_i\nabla c_i(\bx^{\star}), \bmu_i\ge 0\quad \textrm{for}\; i\in\mathcal{A}'(\bx^{\star})\cap\mathcal{I}\right\}
\end{equation}
where $\mathcal{A}'(\bx)$ denotes the set of active constraints at point $\bx^{\star}$.

Let $\bg=\nabla \mathcal{L}(\bx^{\star},\blambda^{\star})$. From \leref{le.farkas} (approximate Farkas lemma), we have that either
\begin{equation}
\nabla \mathcal{L}(\bx^{\star},\blambda^{\star})+\Delta=-\sum_{i\in\mathcal{A}'(\bx^{\star})}\bmu_i\nabla c_i(\bx^{\star}),\quad \bmu_i\ge 0\quad \textrm{for}\; i\in\mathcal{A}'(\bx^{\star})\cap\mathcal{I},
\end{equation}
which is equivalent to
\begin{equation}\label{eq.dis}
\left\|\nabla \mathcal{L}(\bx^{\star},\blambda^{\star})+\sum_{i\in\mathcal{A}'(\bx^{\star})}\bmu_i\nabla c_i(\bx^{\star})\right\|\le\epsilon',\quad \bmu_i\ge 0\quad \textrm{for}\; i\in\mathcal{A}'(\bx^{\star})\cap\mathcal{I}
\end{equation}
or else there is a direction $\bd$ such that $\bd^{\T}\nabla\mathcal{L}(\bx^{\star},\blambda^{\star})<-\epsilon'$ and $\bd\in\cf(\bx^{\star})$ and $\|\bd\|=1$, where
\begin{equation}
\mathcal{F}(\bx)=\left\{\bd\bigg|\begin{cases}\bd^{\T}\nabla c_i(\bx) = 0, & \forall i\in\mathcal{E} \\
\bd^{\T} -\nabla c_i(\bx)\ge 0, & \forall i\in\mathcal{A}'(\bx)\cap\mathcal{I}\end{cases}\right\}
\end{equation}
denotes the cone of linearized feasible directions.

Due to the fact that $c_i(\bx),\forall i$ are convex, we have
\begin{equation}
\sum_{i\in\mathcal{A}'(\bx^{\star})}\bmu_ic_i(\by)\ge \sum_{i\in\mathcal{A}'(\bx^{\star})}\bmu_ic_i(\bx^{\star}) + \sum_{i\in\mathcal{A}'(\bx^{\star})}\bmu_i\nabla^{\T} c_i(\bx^{\star})(\by-\bx^{\star}), \; \forall c_i(\by)\le0, i\in\mathcal{A}'(\bx^{\star})\cap\mathcal{I}.
\end{equation}
Since  $c_i(\by)\le0$, it is obvious that
\begin{equation}\label{eq.a}
-\sum_{i\in\mathcal{A}'(\bx^{\star})}\bmu_ic_i(\bx^{\star}) \ge  \sum_{i\in\mathcal{A}'(\bx^{\star})}\bmu_i\nabla^{\T} c_i(\bx^{\star})(\by-\bx^{\star}), \quad \textrm{for}\; c_i(\by)\le0, i\in\mathcal{A}'(\bx^{\star})\cap\mathcal{I}.
\end{equation}
Also, since $i\in\mathcal{A}'(\bx^{\star})\cap\mathcal{I}$, we have  $c_i(\bx^{\star})=0$, i.e., $\bx^{\star}\notin\textrm{int}(\mathcal{X})$.  Therefore, the above \eqref{eq.a} further yields the following inequality:
\begin{equation}\label{eq.noc}
0 \ge  \sum_{i\in\mathcal{A}'(\bx^{\star})}\bmu_i\nabla^{\T} c_i(\bx^{\star})(\by-\bx^{\star}),   \quad \textrm{for}\; c_i(\by)\le0, i\in\mathcal{A}'(\bx^{\star})\cap\mathcal{I},\bx^{\star}\notin\textrm{int}(\mathcal{X})
\end{equation}

By the definition of the normal cone, we have from \eqref{eq.noc} that
\begin{equation}
\mathcal{N}_{\mathcal{X}}(\bx^{\star})=\begin{cases}\sum_{i\in\mathcal{A}'(\bx^{\star})}\bmu_i\nabla c_i(\bx^{\star}), \quad  & \bmu_i\ge 0, \textrm{for}\; i\in\mathcal{A}'(\bx^{\star})\cap\mathcal{I},\bx^{\star}\notin\textrm{int}(\mathcal{X});\\
\{0\}, & \textrm{otherwise}.\end{cases}
\end{equation}
Combining the facts \eqref{eq.copt} and \eqref{eq.dis} directly gives the result \eqref{eq.eqe}.
\end{proof}

Since we assume that the $\mathcal{X}$ is convex,  directly applying \leref{le.prok} gives  proposition 1.

\vspace{-2pt}
\clearpage
\newpage


\section{Additional Numerical Experiments}\label{sec.addn}
We perform the numerical experiments on a machine with Intel(R) Core(TM) i5-8265U CPU @ 1.60GHz 1.80GHz.

\subsection{More Details of the Experimental Settings in the Main Text}

\noindent{\bf mNPC problem}.   
We take $m=3$, $\lambda=1$, and $r_j=0.1,\forall j$. The initial points of all the algorithms are the same, and $\bx_0$ is randomly generalized, where each entry follows i.i.d. Gaussian distribution $\mathcal{CN}(0,1\times 10^{-3})$. Note that the regularity condition in this problem can be verified, since the sigmoid function is monotonic.
We add i.i.d. Gaussian noise directly to each data with zero mean and unit variance. The initial step-sizes of GDPA are chosen as $\alpha_0=0.1$,  $\beta_0=1\times 10^{-4}$,  and $\tau=0.1$. The initial dual variable of IALM is $\gamma_0=1\times 10^{-2}$, and the initial penalty parameter of IPPP is $1\times 10^{-3}$. For the inner loops of IALM and IPPP, we just use the standard Nesterov's accelerated gradient descent with step-size $0.01$ and momentum parameter $0.1$. For the inner loop of IQRC, the initial step-size is set as $\gamma_0=0.1$, the maximum number of the inner loop is 30, and the predefined feasibility tolerance is 1.

\noindent{\bf Neural nets training with budget constraints}.
The initial step-sizes of GDPA are chosen as $\alpha_0=2\times 10^{-4}$, $\beta_0=2\times 10^{-4}$, and  $\tau=0.1$. The initial dual variable of IALM is $\gamma_0=1\times 10^{-3}$,  and the initial penalty parameter of IPPP is $1\times 10^{-4}$. For the inner loops of IALM and IPPP, we use the standard Nesterov's accelerated gradient descent with step-size $1\times 10^{-4}$ and momentum parameter $0.1$. For the inner loop of  IQRC, the initial step-size is set as $\gamma_0=1\times 10^{-4}$, the maximum number of the inner loop is 30, and the predefined feasibility tolerance is 0.01. The neural net includes two layers, where the hidden layer has 30 neurons, the output dimension of the perception layer is 10, and the activation function is sigmoid.


\subsection{CMDP}

\begin{figure}[h]
\vspace{-0.1cm}
\centering
\subfigure[Objective reward]{
\label{fig:cmdpo}
\includegraphics[width=.43\linewidth]{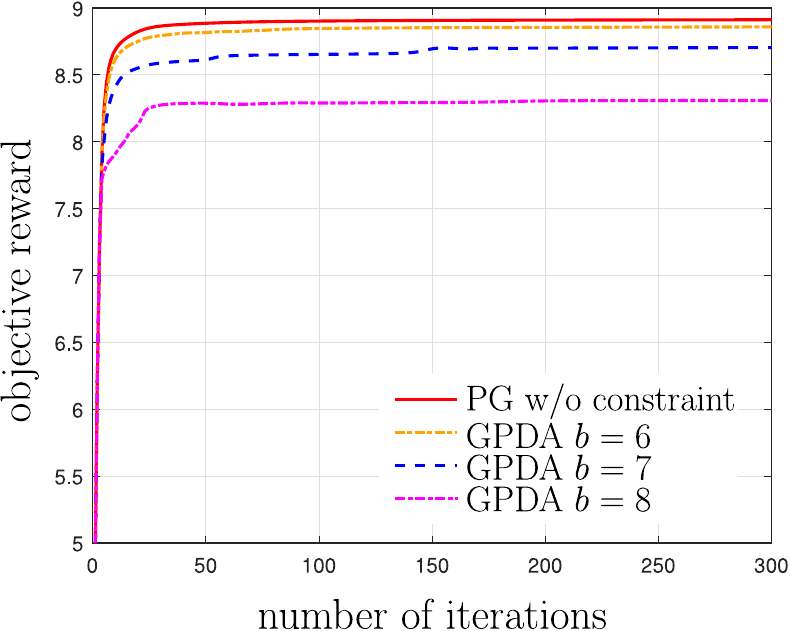}}
\hspace{0.4in}
\subfigure[{Constrained reward}]{
\label{fig:cmdpc}
\includegraphics[width=.43\linewidth]{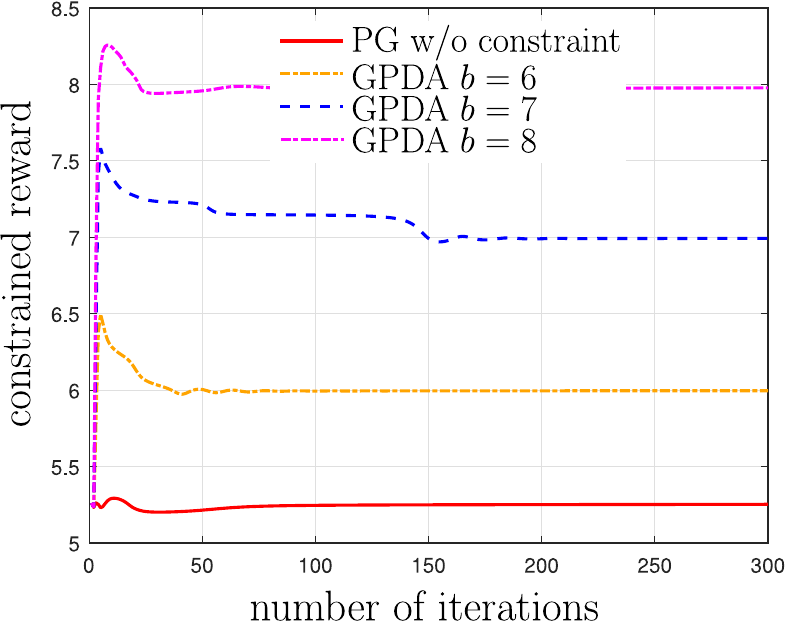}}
\caption{\footnotesize Objective reward v.s. constrained reward achieved by GDPA and PG.}
\label{fig:cmdp}
\end{figure}

We use the code shared in \cite{bhandari2019global} and extend it to CMDP problems for testing the performance of GDPA, where $|\mathcal{S}|=50$, $|\mathcal{A}|=10$ and $\gamma=0.9$. The initial step-sizes of GDPA is $1\times 10^{3}$, $0.5$, and $\tau=0.1$. We set the constraints thresholds as $b=6,7,8$ for three cases. Comparing the classic policy gradient (PG) method, it can be seen in \figref{fig:cmdp} that GDPA can provide the solutions that achieve the predefined constrained rewards while PG fails. Also, it can be observed that if the predefined constrained reward is higher, then the achieved objective rewards will be lower, which makes sense since CMDP is a more complex learning task with multiple optimization objectives than the case without constraints.

\end{document}